\documentclass[11pt,a4paper,reqno]{amsart}

\usepackage{times}
\usepackage{mathrsfs,mathtools}
\usepackage[inline]{enumitem}
\usepackage{stmaryrd}
\usepackage[usenames,dvipsnames]{xcolor}
\usepackage{amsmath}
\usepackage{amsthm}
\usepackage{amssymb}
\usepackage{amsfonts}
\usepackage{subfiles}
\usepackage{dsfont}
\usepackage{hyperref}
\usepackage[T1]{fontenc}
\usepackage[german,english]{babel}
\usepackage{todonotes}
\usepackage{geometry}
\usepackage{xspace,cancel}
\usepackage[capitalize]{cleveref}
\usepackage{etoolbox}

\usepackage[numbers,sort&compress,longnamesfirst]{natbib}

\setlist[enumerate]{label=\textup{(\roman*)},itemindent=.85cm,leftmargin=0.cm}
\newcommand{\oX}{\vphantom{X}\smash{\overline X}}

\makeatletter \@addtoreset{equation}{section}
\def\theequation{\thesection.\arabic{equation}}
\newtheorem{theorem}{Theorem}[section]
\newtheorem*{assumption*}{\assumptionName}
	\providecommand{\assumptionName}{}
	
\newtheorem{corollary}[theorem]{Corollary}
\newtheorem{lemma}[theorem]{Lemma}
\newtheorem{proposition}[theorem]{Proposition}
\newtheorem{definition}[theorem]{Definition}

\theoremstyle{definition}
\newtheorem{remark}[theorem]{Remark}

\newtheorem{example}[theorem]{Example}

\makeatletter
\def\namedlabel#1#2{\begingroup
    #2%
    \def\@currentlabel{#2}%
    \phantomsection\label{#1}\endgroup
}

\numberwithin{equation}{section}

\newcommand{\D}{\ensuremath{\mathbb{D}}}
\newcommand{\E}{\mathbb{E}}
\newcommand{\Fil}{\mathbb{F}}
\newcommand{\FilG}{\mathbb{G}}

\newcommand{\N}{\mathbb{N}}
\newcommand{\Pm}{\mathbb{P}}
\newcommand{\Rp}{\mathbb{R}_+}
\newcommand{\transp}{\top}

\newcommand{\F}{\mathcal{F}}
\newcommand{\G}{\mathcal{G}}
\newcommand{\Pred}{\mathcal{P}}

\newcommand{\cadlag}{c\`adl\`ag\xspace}

\newcommand{\dt}{\ud t}
\newcommand{\ds}{\ud s}
\newcommand{\dx}{\ud x}
\newcommand{\dP}{\ud\Pm}
\newcommand{\extsense}{\xrightarrow{\hspace{0.2cm}\ext\hspace{0.2cm}}}

\newcommand{\mutorthog}
		   {\protect\mathpalette{\protect\independenT}{\perp}}
		   	\def\independenT#1#2{\mathrel{\rlap{$#1#2$}\mkern2mu{#1#2}}}

\newcommand{\inp}{\xrightarrow{\hspace*{0.2cm}\mathbb P\hspace*{0.2cm}}}

\newcommand{\luip}{\xrightarrow{\hspace*{0.2cm}(\lu, \Pm)\hspace*{0.2cm}}}
\newcommand{\mutilde}{\widetilde{\mu}}

\newcommand{\odsdx}{\pair{\omega;\ds,\dx}}

\newcommand{\sip}{\xrightarrow{\hspace*{0.2cm}\left(\J_1(\mathbb R),\Pm\right)\hspace*{0.2cm}}}

\newcommand{\ud}{\ensuremath{\mathrm{d}}}
\newcommand{\weakFil}{\xrightarrow{\hspace*{0.2cm}\w\hspace*{0.2cm}}}

\newcommand{\abs}[1]{\left\vert#1\right\vert}
\newcommand{\borel}[1]{\mathcal{B}\left({#1}\right)}

\newcommand{\Expect}[1]{\mathbb{E}\left[#1\right]}

\newcommand{\Jconv}[1]{\xrightarrow{\hspace{0.2cm}\J_1(\R{#1})\hspace{0.2cm}}}

\newcommand{\inJmulti}[1]{\xrightarrow{\hspace*{0.2cm}\J_1(\mathbb R^{#1})\hspace*{0.2cm}}}

\newcommand{\luil}[1]{\xrightarrow{\hspace{0.2cm}\left(\lu,\mathbb{L}^{#1}\right)\hspace{0.2cm}}}
\newcommand{\norm}[1]{\left\Vert#1\right\Vert}
\newcommand{\pair}[1]{\ensuremath{\! \left( #1 \right) }}
\newcommand{\pqc}[1]{\!\left\langle#1\right\rangle}

\newcommand{\R}[1]{\mathbb{R}^{#1}}

\newcommand{\set}[1]{\ensuremath{\! \left \{ #1 \right \} }}

\newcommand{\sipmulti}[1]{\xrightarrow{\hspace*{0.2cm}\left(\J_1(\R{#1}), \Pm\right)\hspace*{0.2cm}}}

\newcommand{\tnorm}[1]{\left\vert\kern-0.25ex\left\vert\kern-0.25ex\left\vert #1 
    \right\vert\kern-0.25ex\right\vert\kern-0.25ex\right\vert}

\newcommand{\Htwo}[2]{\mathbb{H}^2_{#1}(#2)}
\newcommand{\LGconvmulti}[2]{\xrightarrow{\hspace{0.2cm}\mathbb{L}^{#1}(\Omega,\mathcal{G},\Pm;\R{#2})\hspace{0.2cm}}}

\newcommand{\silmulti}[2]{\xrightarrow{\hspace*{0.2cm}\left(\J_1(\R{#1}), \mathbb{L}^{#2}\right)\hspace*{0.2cm}}}

\newcommand\numberthis{\addtocounter{equation}{1}\tag{\theequation}}

\DeclareMathOperator{\ext}{ext}

\DeclareMathOperator{\J}{J}
\DeclareMathOperator{\lu}{lu}
\DeclareMathOperator{\PUT}{P--UT}

\DeclareMathOperator{\Var}{Var}
\DeclareMathOperator{\w}{w}

\definecolor{red}{rgb}{0.7,0.15,0.15}
\definecolor{green}{rgb}{0,0.5,0}
\definecolor{blue}{rgb}{0,0,0.7}
\hypersetup{	
				colorlinks, linkcolor={red},
				citecolor={green}, urlcolor={blue}
			}

\newtoggle{full} 
\toggletrue{full}

\makeatletter
\newcommand \Dotfill {\leavevmode \leaders \hb@xt@ 6pt{\hss \hss }\hfill \kern \z@}
\makeatother

\makeatletter
\def\@tocline#1#2#3#4#5#6#7{\relax
  \ifnum #1>\c@tocdepth 
  \else
    \par \addpenalty\@secpenalty\addvspace{#2}%
    \begingroup \hyphenpenalty\@M
    \@ifempty{#4}{%
      \@tempdima\csname r@tocindent\number#1\endcsname\relax
    }{%
      \@tempdima#4\relax
    }%
    \parindent\z@ \leftskip#3\relax \advance\leftskip\@tempdima\relax
    \rightskip\@pnumwidth plus4em \parfillskip-\@pnumwidth
    #5\leavevmode\hskip-\@tempdima
      \ifcase #1
       \or\or \hskip 1.65em \or \hskip 3.3em \else \hskip 4.95em \fi%
      #6\nobreak\relax
    \Dotfill
    \hbox to\@pnumwidth{\@tocpagenum{#7}}\par
    \nobreak
    \endgroup
  \fi}
\makeatother

\makeatletter
\def\l@section{\@tocline{1}{0pt}{1pc}{}{\scshape}}
\renewcommand{\tocsection}[3]{%
\indentlabel{\@ifnotempty{#2}{\ignorespaces#1 #2.\hskip 0.7em}}#3}
\def\l@subsection{\@tocline{2}{0pt}{1pc}{5pc}{}}

\def\l@subsubsection{\@tocline{3}{0pt}{1pc}{7pc}{}}

\makeatother

\begin{document}
\setlength\parindent{0pt}

\title[Stability results for martingale representations]{Stability results for martingale representations:\\the general case}

\author[A. Papapantoleon]{Antonis Papapantoleon}
\author[D. Possama\"{i}]{Dylan Possama\"{i}}
\author[A. Saplaouras]{Alexandros Saplaouras}

\address{Department of Mathematics, National Technical University of Athens, Zografou Campus, 15780 Athens, Gree-ce}
\email{papapan@math.ntua.gr}

\address{Department of Industrial Engineering and Operations Research, Columbia University, 500W 120th street, 10027, New York, NY, USA}
\email{dp2917@columbia.edu}

\address{Department of Mathematics, University of Michigan, East Hall, 530 Church Street, Ann Arbor, MI 48109-1043, USA}
\email{asaplaou@umich.edu}

\thanks{We thank Samuel Cohen for useful discussions during the work on these topics. Alexandros Saplaouras gratefully acknowledges the financial support from the DFG Research Training Group 1845 ``Stochastic Analysis with Applications in Biology, Finance and Physics''. Dylan Possama\"i gratefully acknowledges the financial support from the ANR project {\sc Pacman} (ANR-16-CE05-0027). Moreover, all authors gratefully acknowledge the financial support from the {\sc Procope} project ``Financial markets in transition: mathematical models and challenges''.}


\subjclass[2010]{}

\date{}

\begin{abstract}
In this paper, we obtain stability results for martingale representations in a very general framework.
More specifically, we consider a sequence of martingales each adapted to its own filtration, and a sequence of random variables measurable with respect to those filtrations. 
We assume that the terminal values of the martingales and the associated filtrations converge in the extended sense, and that the limiting martingale is quasi--left--continuous and admits the predictable representation property.
Then, we prove that each component in the martingale representation of the sequence converges to the corresponding component of the martingale representation of the limiting random variable relative to the limiting filtration, under the Skorokhod topology. 
This extends in several directions earlier contributions in the literature, and has applications to stability results for backward SDEs with jumps and to discretisation schemes for stochastic systems.
\end{abstract}

\maketitle \frenchspacing

\iftoggle{full}{\tableofcontents}{}

\vspace{-2em}
\section{Introduction}

Consider a sequence $(X^n)_{n\in\mathbb N}$ of square--integrable martingales, which is assumed to converge to another square--integrable martingale $X^\infty$, the convergence holding either in the \textit{strong} sense, meaning in particular that all the martingales $X^n$, as well as $X^\infty$, are defined on a common probability space, or in the \textit{weak} sense, that is, the convergence is in distribution, and each $X^n$ is then defined on its own probability space. 
For every $n\in\mathbb N$, let us denote by $\mathbb G^n:=(\mathcal G_t^n)_{t\geq 0}$ the filtrations with respect to which $X^n$ are martingales and by $\mathbb G^\infty:=(\mathcal G^\infty_t)_{t\geq 0}$ the one associated to $X^\infty$. 
Given now a sequence of random variables $(\xi^n)_{n\in\mathbb N}$, where $\xi^n$ is respectively $\mathcal G^n_\infty-$measurable, based on a well--known result (see, for instance, \citet[Lemma III.4.24]{jacod2003limit}) the martingales $Y^n_\cdot:=\mathbb E[\xi^n|\mathcal G^n_\cdot]$ admit a so--called orthogonal decomposition with respect to $X^n$. 
In other words, for every $n\in\mathbb N$, let $X^{n,c}$ be the continuous part of $X^n$ and $\widetilde \mu^{X^{n,d}}$ be the (compensated) random measure of jumps associated to $X^{n,d}$, \emph{i.e.} the purely discontinuous part of $X^n$, then 
\begin{equation}\label{eq:rep-n}
Y^n_\cdot = Y^n_0 + \int_0^\cdot Z^n_s \ud X^{n,c}_s + \int_0^\cdot\int_{\mathbb R^\ell}U^n_s(x)\widetilde \mu^{X^{n,d}}(\ds,\dx) + N^n_\cdot,	
\end{equation}
where $Z^n$ and $U^n$ are respectively a predictable process and a predictable function, while $N^n$ is another martingale, appropriately orthogonal to both the continuous and the discontinuous martingale parts of $X^n$. 

\medskip
Assume now that the sequence of pairs $(X^n,\mathbb G^n)_{n\in\mathbb N}$ converges (in the extended sense) to $(X^\infty,\mathbb G^\infty)$, and that the sequence $(\xi^n)_{n\in\mathbb N}$ converges, in an appropriate sense, to some $\mathcal G^\infty_\infty-$measurable random variable $\xi^\infty$, such that the following orthogonal decomposition for $Y_\cdot:=\mathbb E[\xi^\infty|\mathcal G^\infty_\cdot]$ with respect to $X^\infty$ holds
\begin{equation}\label{eq:rep}
Y^\infty_\cdot = Y^\infty_0 + \int_0^\cdot Z^\infty_s \ud X^{\infty,c}_s + \int_0^\cdot\int_{\mathbb R^\ell}U^\infty_s(x)\widetilde{\mu}^{X^{\infty,d}}(\ds,\dx) + N^\infty_\cdot.
\end{equation}
A natural question is then whether the convergence of $Y^n$ to $Y^\infty$ implies also the convergence of the martingale parts on the right--hand side of \eqref{eq:rep-n} to their respective counterparts on the right--hand side of \eqref{eq:rep}.
A weaker version of the posed question is whether the sequence consisting of the sum of the stochastic integrals in \eqref{eq:rep-n} converges to the sum of the stochastic integrals in \eqref{eq:rep} and therefore also the sequence of the orthogonal martingales $(N^n)_{n\in\mathbb N}$ converges to $N^\infty$. 
 
\medskip
This problem of approximations of certain martingale representations has a long history, which was mainly motivated by applications in mathematical finance.
There the random variables $\xi^n$ can be understood as contingent claims to be hedged using financial assets whose prices are given by $X^n$, and where $Z^n$ are then appropriate hedging strategies (usually risk minimising). 
In this context, $U^n$ would typically appear when the price processes $X^n$ can have jumps, and $N^n$ would sum up all the information in the filtration $\mathbb G^n$ which cannot be generated by $X^{n,c}$ or $\widetilde \mu^{X^{n,d}}$. 
This is the typical situation encountered in so--called incomplete financial markets. 
Furthermore, the approximation of $X$ by $X^n$ usually stems from computational considerations, typically using discretisation schemes for practical and efficient implementations. 
The question of whether the associated hedging strategies converge or not, and in which sense, is then of paramount importance. 
This was notably the subject of \citet{jacod2000explicit}, which considers a setting where the $U^n$ do not appear, since the stochastic integral is an integral with respect to $X^n$ (and not only its continuous martingale part; this is the celebrated Galtchouk--Kunita--Watanabe decomposition from martingale theory), and where the $\xi^n$ are Markovian functionals of $X^n$. 
Earlier contributions by \citet{jakubowski1989convergence} and then \citet{kurtz1991weak,kurtz1996weak} had already studied, from a theoretical point of view, the simpler question of the weak convergence of stochastic integrals of the form $\int_0^\cdot Z^n_s \ud X^{n}_s,$ while Duffie and Protter had investigated the aforementioned financial applications in \cite{duffie1992discrete}.
 
\medskip
The problem posed above is also intimately linked to the study of weak convergence of discretisation schemes for stochastic systems, which has been a topic of continued interest in stochastic numerical analysis and its applications. 
As illustrated in several articles, there are discretisation schemes for such systems which do not lead to satisfactory stability properties, especially for the simplest and elementary processes, such as stochastic integrals and stochastic differential equations. 
These questions, in a context similar to ours, have been investigated for instance by \citet{barlow1990convergence}, for the stability of special semimartingale decompositions, by \citet{coquet1999convergence} for Dirichlet processes, and by \'Emery \cite{emery1978stabilite,emery1979equations}, \citet{protter1978H,protter1985approximations,mackevicius1986sp,mackevicius1987sp} and \citet{slominski1989stability} for strong solutions of stochastic differential equations. 
More recently, this also was the direction followed by \citet{leao2013weak}, where the authors aimed at describing readable structural conditions on a given optional process adapted to a Brownian filtration, in order to construct explicit, robust and feasible approximating skeletons for smooth semimartingales.
Closedness results for stochastic integrals with respect to (local) martingales are also part of the folklore of the general theory of processes. 
This, roughly speaking, corresponds to the case where one simply studies integrals of the formv$\int_0^\cdot Z_s {\rm d}X^{n}_s.$
Hence, the case of integrals in $\mathbb L^2$ (or more generally in $\mathbb L^p$, $p>1$) is straightforward, coming almost directly from the Hilbert space isometry of stochastic integrands and integrals, see for instance \citet[p. 153]{protter2005stochastic} or \citet{jacod1979calcul}. 
The much more subtle case of martingales in $\mathbb L^1$ was settled by \citet{yor1978sous}, see also \citet{delbaen1999compactness} for a survey of these results, as well as additional compactness criteria. 
The case where $X$ is allowed to be a semimartingale is naturally quite more involved, and comprehensive results in this direction were obtained by \citet{memin1980espaces}, \citet{schweizer1995variance}, \citet{monat1994fermeture,monat1995follmer}, and \citet{delbaen1994inegalites,delbaen1997weighted}.
 
\medskip
Once the semimartingales considered have more structural properties, other interesting results can be obtained. 
\citet{barrieu2008closedness}, for instance, and later \citet{barrieu2013monotone} were interested in what they coined ``continuous quadratic semimartingales'' (see also related articles by \citet{mocha2012quadratic}, still in the continuous case, and recent extensions to jump processes by \citet{ngoupeyou2010optimisation} and \citet{karoui2016quadratic}), for which, roughly speaking, the bounded variation process part in the semimartingale $X$ is absolutely continuous with respect to the quadratic variation of the martingale part of $X$.
\cite{barrieu2008closedness,barrieu2013monotone} obtained associated stability results for these processes.
 
\medskip
An important common feature of the articles mentioned so far, is that they actually only consider the \textit{strong} framework we described at the beginning of this introduction, in the sense that there is a always a fixed probability space and all processes (meaning here mainly $X^n$ and $X^\infty$) are adapted to the same fixed filtration $\mathbb G$.
An important exception is \citet{slominski1989stability}, where the probability space is fixed, but not the filtration. 
However, for practical purposes, and especially for the analysis of numerical schemes, it is well--known that the \textit{weak} framework is also of paramount importance, as illustrated for instance by the famous Donsker theorem. 
There has thus been a certain number of studies of stability properties for semimartingale or martingale decompositions when the underlying filtration itself is also allowed to change. 
In that direction, \citet{antonelli2000filtration}, followed by \citet{coquet1998stability,coquet1999corrigendum}, \citet{ma2002numerical}, \citet{briand2001donsker,briand2002robustness}, and then \citet{cheridito2013bs}, studied such stability properties for continuous backward stochastic differential equations (BSDEs for short), a type of non--linear martingale representation.
\citet{memin2003stability} looked into the stability of the canonical decomposition for semimartingales, \citet{kchia2011semimartingales} (see also \citet{kchia2015progressive}) extended Barlow and Protter's result \cite{barlow1990convergence} for stability of special semimartingale decompositions to a framework allowing changing filtrations, while \citet{possamai2015weak} extended results of \cite{briand2001donsker,briand2002robustness} to the case of so--called second order BSDEs. 
Let us also mention the recent paper by \citet{madan2015convergence}, which considers stability results for BSDEs with jumps (that is to say that both the processes $Z$ and $U$ are present in the solution), when the driving c\`adl\`ag martingale is approximated by random walks. 
Several of these works make a strong use of the notions of extended convergence, introduced by \citet{aldous1981weak}, as well as that of convergence of filtrations, introduced by \citet{hoover1991convergence} and further developed by \citet{coquet2000some} and \citet{coquet2001weak}, which also plays a major role in the present paper. Let us also mention the recent contributions by Le\~ao, Ohashi and Simas \cite{leao2018weak2,leao2018weak}, who study stability of Wiener functionals under weak convergence of filtration beyond semimartingales, in the context of functional It\=o calculus.

\medskip
Our work follows this latter strand of literature and studies the problem of stability for the martingale representation (or the orthogonal decomposition) of the martingales $X^n$, when their filtration is also allowed to change. 
A very important difference compared to the existing literature is that we basically make no assumption on the filtration $\mathbb G^n$, besides the minimal ones, \textit{i.e.} that they satisfy the usual assumptions of right--continuity and completeness under a fixed reference probability measure $\mathbb P$, and that they converge in an appropriate sense to the filtration $\mathbb G^\infty$ associated to $X^\infty$. 
This means, in particular, that the filtrations $\mathbb G^n$ are not constrained to be quasi--left--continuous, a property often considered in the existing literature, and whose relaxation highly complicates the problem. 
We believe that this is somehow the highest degree of generality one can consider while still remaining in the martingale framework. 
However, this level of generality comes at the price that we have to assume more properties for the limiting filtration $\mathbb G^\infty$. 
More precisely, we have to assume that $\mathbb G^\infty$ is quasi--left--continuous (hence also the martingale $X^\infty$) and that the predictable representation property holds for $\mathbb G^\infty$ and $X^\infty$ (meaning that $N^\infty$ in \eqref{eq:rep} above must vanish). 
Although the first assumption is somehow unavoidable in such a setting, as illustrated by \citet{memin2003stability}, the second one is slightly more restrictive. 
A proper discussion of the reasons why our approach cannot work without it requires lengthy preliminaries, hence we postpone it to Subsection \ref{sec:comp} below. 
Our results stipulate that under these assumptions, the extended convergence of $(\xi^n,\mathbb G^n)$ implies the joint convergence in the Skorokhod topology of $(Y^n,Z^n\cdot X^n+U^n\star\widetilde{\mu}^{X^{n,d}},N^n)$, but also of the angle brackets $(\langle Y^n\rangle,\langle Y^n,X^n\rangle,\langle N^n\rangle)$\footnote{The convergence of the respective square bracket processes is also obtained, but this is well-known in the existing literature.}. 
In case the processes $X^n$ have in addition independent increments, we can obtain that the above convergences also hold in law when we work under the natural filtration associated to $X^n$ (see \cref{corr:PII}). 
Besides, if we assume that there exist two sequences that converge to the continuous and the purely discontinuous part of the limiting martingale respectively, then the angle brackets of $Y^n$ with respect to these sequences converge to the angle brackets of $Y^\infty$ with respect to the continuous and the purely discontinuous part of the limiting martingale, see \cref{StrongRobMartRep}.

\medskip
On the way to proving our results, we needed to apply the Burkholder--Davis--Gundy inequality, as well as the Doob inequality in their general form, namely for a (suitable) moderate Young function $\Phi$.\footnote{For the definition see \cref{App:Young}.} This finally allowed us to have a ``sharp'' $\mathbb L^2-$convergence in our results, and not simply $\mathbb L^{2+\epsilon}$, for some $\varepsilon>0$, as it is usually imposed in the literature in order to have sufficient integrability. 
As the reader may suspect, this result was possible due to the special role that $p=2$ plays in the general $\mathbb L^p-$theory.

\medskip
The rest of the paper is organised as follows. 
Section \ref{sec:prel} introduces all the relevant notions from stochastic analysis and stochastic integration, as well as from the study of the Skorokhod space and the extended convergence of Aldous. 
Section \ref{sec:RobMartRepSection} is then devoted to the statement of our main results, a comprehensive comparison with the existing literature, as well as a very detailed explanation of our strategy of proof. 
The proof itself follows, while the appendices collect important technical results.

\medskip
{\bf Notation.} 
Let $\mathbb R_+$ denote the set of non--negative real numbers, $\overline{\mathbb R}_+:=[0,\infty]$ and $\overline{\N}:=\N\cup \{\infty\}$. 
For any positive integer $\ell$, any $x\in \R \ell$ will be identified as a \emph{column} vector of length $\ell$, $x^i$ will denote the $i-$th element of $x$ and $\pi^i$ will denote the canonical $i-$projection 
$\mathbb R^\ell\ni x\longmapsto x^i\in \mathbb R,$ for $1\leq i \leq \ell.$ 
The identity function $\mathbb R^\ell\ni x\longmapsto x\in\mathbb R^\ell$ will be denoted by ${\rm Id_\ell},$ where we will suppress the index when the dimension is clear. 
By $\abs x$ we will denote the usual Euclidean norm of $x$, while the metric compatible with the topology imposed by the Euclidean norm will be denoted by $d_{|\cdot|}.$
For any additional positive integer $q$, a $q\times \ell-$matrix with real entries will be considered as an element of $\R {q\times \ell}$. 
For any $z\in\R {q\times \ell}$, its transpose will be denoted by $z^\top\in\R {\ell\times q}$.
The element at the $i-$th row and $j-$th column of $z\in\R {q\times \ell}$ will be denoted by $z^{ij}$, for $1\leq i\leq q$ and $1\leq j\leq \ell$. 
The trace of a square matrix $z\in\R{\ell\times\ell}$ is ${\rm Tr}[z]:=\sum_{i=1}^{\ell}z^{ii}.$
We endow $\R {q\times \ell}$ with the $\norm{\cdot}_2-$norm defined for any $z\in\R {q\times \ell}$ by $\norm z^2_2:={\rm Tr}[z^\transp z]$ and remind the reader that this norm is derived from the inner product defined for any 
$(z,u)\in\R {q\times \ell}\times\R {q\times \ell}$ by ${\rm Tr}[z^{\transp}u]$. 
Moreover, we will also make use of the $\norm{\cdot}_1-$norm, which is defined as $\Vert z\Vert_1:=\sum_{i=1}^q\sum_{j=1}^{\ell} |z^{ij}|$, for $z\in\R{q\times \ell}$.
We abuse notation and denote by $0$ the neutral element in the groups $(\mathbb R^{\ell},+)$ and $(\mathbb R^{q\times\ell},+)$. 
Throughout the rest of the paper $p,q, \ell$ will always denote natural integers and, in particular, $\ell$ will be fixed.

\vspace{0.5em}
Let $E$ denote a finite dimensional topological space, then $\mathcal B(E)$ will denote the associated Borel $\sigma-$algebra. 
Furthermore, for any other finite dimensional topological space $G$ and for any non--negative measure $\rho$ defined on $(\mathbb R_+,\mathcal B(\mathbb R_+))$, we will denote the Lebesgue--Stieltjes integral, with respect to some measure $\rho$ on $(\mathbb R_+,\mathcal B(\mathbb R_+))$, of any measurable map $f:(\mathbb R_+,\mathcal B(\mathbb R_+)) \longrightarrow (G,\mathcal B(G))$  by 
$$ 	\int_{(u,t]}f(s)\rho(\ds) \ \text{ and } \ \int_{(u,\infty)}f(s)\rho(\ds),  \text{ for any } u,t\in\mathbb R_+.$$
In case $\rho$ is a finite measure with associated distribution function $F^{\rho}(\cdot):=\rho([0,\cdot])$, we will  indifferently denote the above integrals by
$$ 	\int_{(u,t]}f(s)\ud F^{\rho}_s \ \text{ and } \ \int_{(u,\infty)}f(s)\ud F^{\rho}_s, \text{ for any } u,t\in\mathbb R_+.$$
When there is no confusion as to which measure the distribution function $F^{\rho}$ is associated to, we will omit the upper index and simply write $F$.
More generally, for any measure $\varrho$ on $(\mathbb R_+\times E,\mathcal B(\mathbb R_+)\otimes\mathcal B(E))$ and for any measurable map
$g:(\mathbb R_+\times E,\mathcal B(\mathbb R_+)\otimes\mathcal B(E))\longrightarrow (G,\mathcal B(G))$ we will denote the Lebesgue--Stieltjes integral by
\begin{align*}
\int_{(u,t]\times A}g(s,x)\varrho(\ds,\dx) \  \text{ and }\  \int_{(u,\infty)\times A}g(s,x)\varrho(\ds,\dx), \text{ for any } t,u\in\mathbb R_+, A\in\mathcal B(E).
 \end{align*}
The integrals above are to be understood in a component--wise sense. 

\vspace{0.5em}
Finally, we recall, for the convenience of the reader, some classical terminology.
Let $(E,d_E)$ be a Polish space.
We denote by $\mathcal P(E)$ the set of all probability measures on $(E,\mathcal B(E)).$
We endow $\mathcal P(E)$ with the \emph{weak topology}\footnote{\label{weaklimit} In functional analysis, this is called the \emph{weak$^\star-$topology}.}, \emph{i.e.} the coarsest topology for which the mappings 
$\mathcal P(E)\ni\varrho \longmapsto \int_{E}f\ud \varrho\in\mathbb R$, are continuous for all bounded continuous functions $f$ on $E.$ 
It is well--known that $E$ is Polish if and only if $\mathcal P(E)$ is Polish;  see \citet[Theorem 15.15]{aliprantis2006infinite} or \citet[Thoerem 6.5]{parthasarathy1972probability}.
Moreover, a subset $\Gamma$ of $\mathcal P(E)$ is \emph{relatively compact} for the weak topology if and only if it is \emph{tight}, see \cite[Theorem 15.22]{aliprantis2006infinite} or \cite[Theorem 6.7]{parthasarathy1972probability}.
For a random variable $\Xi:(\Omega,\mathcal G)\longrightarrow (E,\mathcal B(E))$, its \emph{law} $\mathcal L(\Xi)$ is defined as $\mathcal L(\Xi)(A):=\mathbb P(\{\omega\in\Omega, \Xi(\omega)\in A\}),$ for every $A\in \mathcal B(E).$
We will say that the sequence of random variables $(\Xi^k)_{k\in\N}$ \emph{converges in law} to the random variable $\Xi^\infty,$ and we will write $\Xi^k\xrightarrow{\hspace{0.2cm}\mathcal L\hspace{0.2cm}}\Xi^\infty,$ if $\mathcal L(\Xi^k)$ converges weakly to $\mathcal L (\Xi^\infty),$ which will be denoted as $\mathcal L(\Xi^k)\xrightarrow{\hspace{0.2cm}{\rm w}\hspace{0.2cm}}\mathcal L(\Xi^\infty).$
Moreover, we will say that the sequence of random variables $(\Xi^k)_{k\in\overline\N}$ is tight, if the associated sequence of laws $\big(\mathcal L(\Xi^k)\big)_{k\in\N}$ is tight.
Finally, for a tight sequence $(\Xi^k)_{k\in\N}$ of random variables, we will say that $\overline \Xi$ is a weak--limit point if there exists a subsequence $(\Xi^{k_l})_{l\in\overline\N}$ such that $\mathcal L(\Xi^{k_l})\xrightarrow{\hspace{0.2cm}{\rm w}\hspace{0.2cm}} \mathcal L(\overline \Xi)$  holds.

\section{Preliminaries}\label{sec:prel}

\subsection{The stochastic basis}

Let $(\Omega,\mathcal G,\mathbb P)$ be a probability space, which is fixed for the remainder of this paper. 
Expectations under $\mathbb P$ will be denoted by $\mathbb E[\cdot]$. 
For any filtration\footnote{\label{UCFil}We assume that all filtrations considered satisfy the usual conditions of right--continuity and $\mathbb P-$completeness.} $\mathbb F$ on $(\Omega, \mathcal G, \Pm)$ and for any $\mathbb F-$stopping time $\tau$, we will denote the set of $\R{q}-$valued and square--integrable $\mathbb F-$martingales stopped\footnote{For a process $M$, the corresponding process stopped at $\tau$, denoted by $M^\tau$, is defined by $M^\tau_t:=M_{t\wedge\tau}$, $t\geq 0$.} at $\tau$ by $\mathcal{H}^2(\mathbb F, \tau;\R{q})$.
A process $(M_t)_{t\in\mathbb R_+}$ will be denoted also as $M$, and the usual augmentation of its natural filtration will be denoted by $\mathbb F^M$.
Let $M\in\mathcal{H}^2(\mathbb{F}, \tau;\R{q})$, then its norm is defined by $\norm{M}^2_{\mathcal{H}^2(\mathbb F,\tau; \R{q})}:=\Expect{{\rm Tr}[\, \pqc{M}_\tau]}.$ In the sequel, we will say that the real--valued $\mathbb F-$martingales $L, M$ are \emph{$($mutually$)$ orthogonal}, denoted by $L \mutorthog M$, if their product $LM$ is an $\mathbb F-$martingale; see \citet[Definition I.4.11.a, Lemma I.4.13.c]{jacod2003limit}.
An $\R{q}-$valued, $\mathbb F-$martingale $L$ will be called a \emph{continuous} martingale if $L_0=0$ and $L^i$ is a continuous $\mathbb F-$martingale, for each $i=1,\ldots,q$.
Moreover, an $\R{q}-$valued, $\mathbb F-$martingale $M$ will be called a \emph{purely discontinuous} martingale if $M_0=0$ and $M^i$ is orthogonal to all continuous real--valued $\mathbb F-$martingales, for each $i=1,\ldots,q$. 
Using \cite[Corollary I.4.16]{jacod2003limit} we can decompose the space of square integrable $\mathbb F-$martingales as follows
\begin{align*}
\mathcal{H}^{2}(\mathbb F,\tau;\R{q})
	&=\big(\mathcal{H}^{2}(\mathbb F,\tau;\R{}) \big)^{q}=\big(\mathcal{H}^{2,c}(\mathbb F,\tau;\R{})\oplus\mathcal{H}^{2,d}(\mathbb F,\tau;\R{})\big)^{q}=\mathcal{H}^{2,c}(\mathbb F,\tau;\R{q})\oplus\mathcal{H}^{2,d}(\mathbb F,\tau;\R{q}),
\end{align*}
where we have defined for any $m\in\N$
\begin{align*}
\mathcal{H}^{2,c}(\mathbb F,\tau;\R{m})&:=\big\{ M\in \mathcal{H}^{2}(\mathbb F,\tau;\R{m}), \, M \text{ is continuous}\big\}, \\
\mathcal{H}^{2,d}(\mathbb F,\tau;\R{m})&:=\big\{ M\in \mathcal{H}^{2}(\mathbb F,\tau;\R{m}), \, M \text{ is purely discontinuous} \big\}.
\end{align*}
Then, it follows from \cite[Theorem I.4.18]{jacod2003limit}, that any $M\in\mathcal{H}^{2}(\mathbb F,\tau;\R{q})$ admits a unique decomposition, up to $\mathbb P-$indistinguishability $M_\cdot=M_0+M^c_\cdot+M^d_\cdot,$ where $M^c_0=M^d_0=0$.
The process $M^c=(M^{c,1}, \dots,M^{c,{q}})\in\mathcal{H}^{2,c}(\mathbb F,\tau;\R{q})$ will be called the \emph{continuous $($martingale$)$ part of $M$} and the process $M^d=(M^{d,1}, \ldots,M^{d,{q}})\in\mathcal{H}^{2,d}(\mathbb F,\tau;\R{q})$ will be called the \emph{purely discontinuous $($martingale$)$ part of $M$}.

\subsection{Stochastic integrals}\label{StochIntegrals}

Let us fix an arbitrary filtration $\mathbb F$ on $(\Omega,\mathcal G, \Pm)$ and an arbitrary $\mathbb F-$stopping time $\tau$.
The predictable $\sigma-$field generated by $\mathbb F-$adapted and left--continuous processes on $\Omega\times\mathbb R_+$ is denoted by $\mathcal P^\mathbb F$.

\subsubsection{\texorpdfstring{It\=o stochastic integral}{Ito stochastic integral}}

We will follow \cite[Section~III.6a]{jacod2003limit} throughout this sub--sub--section.
Let $X\in\mathcal H^{2}(\mathbb F,\tau;\mathbb R^\ell).$
There exists an $\mathbb F-$predictable, \cadlag
and increasing process $C^{X}$ such that 
$$\langle {X} \rangle^{\tau}_\cdot= \int_{(0,\cdot\wedge\tau]}^\cdot \frac{\ud \langle {X} \rangle_s}{\ud C^{X}_s} \,\ud C^{X}_s,$$ 
where $\ud \langle {X} \rangle/\ud C^{X}$ is a positive definite, symmetric and $\mathbb F-$predictable ${\ell\times\ell}-$matrix whose elements are defined by 
$$\Big(\frac{\ud \langle {X} \rangle_\cdot}{\ud C^{X}_\cdot}\Big)^{ij} := \frac{\ud \langle {X} \rangle^{ij}_\cdot}{\ud C^{X}_\cdot},\;\; \text{for $i,j=1,\ldots,\ell.$}$$
Let us now proceed by defining
\begin{align*}
\mathbb H^2(X,\mathbb F,\tau;\R{\ell}) 
	:=\bigg\{& 
		Z:(\Omega\times\Rp,\Pred^\mathbb F) \longrightarrow (\R{\ell},\mathcal B(\R{\ell})),\;
		\E\bigg[
		\int_{(0,\tau]} Z_t^{\transp} \frac{\ud \langle {X} \rangle_t}{\ud C^{X}_t} Z_t\ud C^{X}_t
		\bigg]<\infty
		\bigg\}.
\end{align*}
Notice that this space does not depend on the choice of $C^X$.
For any $Z\in\mathbb H^2(X,\mathbb F,\tau;\R{\ell})$, the It\=o stochastic integral of $Z$ with respect to $X$ is well defined and is an element of $\mathcal{H}^{2}(\mathbb F,\tau;\R{})$. 
It will be denoted interchangeably by $\int_0^{\cdot}Z_s \ud {X}_s$ or $Z\cdot {X}$. 
Moreover, we have the following equality
\begin{align*}
\norm{Z}_{\mathbb H^2(X,\mathbb F,\tau;\R{\ell})}^2:=
		\E\bigg[
		\int_{(0,\tau]} Z_t^{\transp}\frac{\ud \langle X \rangle_t}{\ud C^{X}_t} Z_t\ud C^{X}_t\bigg]
	=	\E\big[\rm{Tr}[\langle Z\cdot X\rangle_\tau]\big].
\end{align*}
We denote the space of It\=o stochastic integrals of processes in the space $\mathbb H^2(X,\mathbb F,\tau;\R{\ell})$, with respect to $X$, by $\mathcal{L}^2(X,\mathbb F,\tau;\R{})$, and remind the reader that $\mathcal{L}^2(X,\mathbb F,\tau;\R{})\subset\mathcal H^2(\mathbb F,\tau;\R{}).$

\iftoggle{full}{\begin{remark}
In case $X$ is a continuous martingale, \textit{i.e.} $X\in\mathcal H^{2,c}(\mathbb F,\tau;\mathbb R^\ell)$, then for any $Z\in\mathbb H^2(X,\mathbb F,\tau;\R{\ell})$ it holds that $Z\cdot X$ is an element of $\mathcal{H}^{2,c}(\mathbb F,\tau;\R{})$, from which it follows from \cite[Section~III.4a]{jacod2003limit} that
$$\mathcal{L}^2(X,\mathbb F,\tau;\R{})\subset\mathcal H^{2,c}(\mathbb F,\tau;\R{})\subset\mathcal H^2(\mathbb F,\tau;\R{}).$$
\end{remark}}{}

\subsubsection{Stochastic integral with respect to an integer--valued random measure}\label{sec:StochIntIVM}

Let us now define the space $\widetilde{\Omega}:=\Omega\times\Rp\times\R {\ell}$ as well as the $\sigma-$algebra $\widetilde{\Pred}^\mathbb F:=\Pred^\mathbb F\otimes \borel{\R {\ell}}.$ 
A measurable function $U:\big(\widetilde{\Omega},\widetilde{\Pred}^\mathbb F\big)\longrightarrow \left(\R{} ,\borel{\R{} } \right)$ is called an \emph{$\widetilde{\Pred}^\mathbb F-$measurable function} or simply \emph{$\mathbb F-$predictable function}.

\vspace{0.5em}
Let $\mu := \set{\mu\pair{\omega;\dt,\dx}}_{\omega\in\Omega}$ be a random measure on $\Rp\times\R{\ell}$, \textit{i.e.} a family of non--negative measures defined on $\pair{\Rp\times\R{\ell},\borel{\Rp}\otimes\borel{\R{\ell}}}$ satisfying $\mu\pair{\omega;\set{0}\times\R{\ell}}=0$, identically.
For an $\mathbb F-$predictable function $U$, we define the process
$$
U*\mu_\cdot(\omega) :=
			\begin{cases}
			\displaystyle \int_{(0,\cdot]\times\R{\ell}} U\pair{\omega,s,x} \mu\odsdx,\textrm{ if } \displaystyle\int_{(0,\cdot]\times\R{\ell}} \abs{U\pair{\omega,s,x}} \mu\odsdx<\infty,\\
		\displaystyle	\infty,\textrm{ otherwise}.
			\end{cases}
$$
Let us fix an arbitrary \cadlag $\mathbb F-$adapted process $X$ until the end of the present sub--sub--section, from which we can define the processes 
$X_-=(X_{t-})_{t\in\mathbb R+}$ and $\Delta X = (\Delta X_t)_{t\in\mathbb R_+}$\iftoggle{full}{, where
\begin{align*}
X_{0-}:=X_0, \hspace{0.3em} X_{t-}:=\lim_{\substack{s\to t \\ s<t}} X_s \hspace{0.3em}\text{ and } \hspace{0.3em} \Delta X_t:=X_t-X_{t-},\; t>0.
\end{align*}}{.}
Observe that $\Delta X_0=0.$
We assume that $X$ satisfies $\mathbb E\big[\sum_{0<s\le\tau}|\Delta X_s|^2\big]<\infty$.
\iftoggle{full}{We can associate to $X$ the $\mathbb F-$optional integer--valued random measure $\mu^{X}$ on $\pair{\Rp\times\R{\ell},\borel{\Rp}\otimes\borel{\R{\ell}}}$ defined by
\[
\mu^{X}\pair{\omega;\dt,\dx} 
	:= \sum_{s>0} \mathds{1}_{\set{\Delta {X}_s(\omega) \neq 0}} \delta_{(s,\Delta {X}_s(\omega))}\pair{\dt,\dx},
\]
see \cite[Proposition II.1.16]{jacod2003limit}.
Here $\delta_z$ denotes the Dirac measure at the point $z$, for any $z\in\Rp\times\R{\ell}$. }
{We can associate to $X$ an $\mathbb F-$optional integer--valued random measure denoted by $\mu^{X}$, see \cite[Proposition II.1.16]{jacod2003limit}.}
Notice that $\mu^{X}(\omega;\mathbb R_+\times\{0\})=0$, and that 
$$\mathbb E\bigg[\int_{(0,\tau]}|x|^2 \mu^X(\ds,\dx)\bigg ]=\mathbb E\Bigg[\sum_{0<s\le\tau}|\Delta X_s|^2\Bigg]<\infty.$$
In view of the latter condition, the compensator of $\mu^{X}$ under $\Pm$ exists.
This is the unique, up to a $\Pm-$null set, $\mathbb F-$predictable random measure $\nu^{({X},\mathbb F)}$ on $\pair{\Rp\times\R{\ell},\borel{\Rp}\otimes\borel{\R{\ell}}}$, for which the equality
\begin{align*}
\E \Big[ W*\mu^{X}_\infty \Big]=\Expect{W*\nu^{({X},\mathbb F)}_\infty},
\end{align*}
holds for every non--negative $\mathbb F-$predictable function $W$; see \cite[Theorem II.1.8]{jacod2003limit}.
Moreover, we define the \emph{compensated integer--valued random measure \mbox{$\mutilde^{({X},\mathbb F)}\!:=\mu^{X}-\nu^{({X},\mathbb F)}$}}.

\vspace{.5em}
In order to define the stochastic integral of an $\mathbb F-$predictable function $U$ with respect to $\mutilde^{({X},\mathbb F)}$, we will consider the following class
\begin{align*}
G_2(\mu^{X},\mathbb F,\tau;\R{}):=\Bigg\{
	U:\big(\widetilde{\Omega},\widetilde{\Pred}^\mathbb F\big)\longrightarrow \big(\R{},\mathcal B(\R{})\big),\, 
	\E\Bigg[\sum_{t\leq \tau} \bigg|\int_{\R{\ell}}U(t,x)\mutilde^{({X},\mathbb F)}(\{t\}\times\dx)\bigg|^2\Bigg]
	<\infty
\Bigg\},
\end{align*}
see also \citet[Section 2.2]{papapantoleon2016existence} for more details.
It is a well--known result that any element of $G_2(\mu^{X},\mathbb F,\tau;\R{})$ is associated to an element of 
$\mathcal{H}^{2,d}(\mathbb F,\tau;\R{}),$ which is unique up to $\mathbb P-$indistinguishability, see \cite[Defintion II.1.27, Proposition II.1.33.a]{jacod2003limit}
\begin{align*}
G_2(\mu^{X},\mathbb F,\tau;\R{})\ni U\longmapsto U\star\mutilde^{({X},\mathbb F)}\in \mathcal{H}^{2,d}(\mathbb F,\tau;\R{}).
\end{align*}
We call $U\star\mutilde^{({X},\mathbb F)}$ the \emph{stochastic integral of $U$ with respect to $\mutilde^{({X},\mathbb F)}$}, and point out that for $U\in G_2(\mu^{X},\mathbb F,\tau;\R{})$ it holds 
$\Delta(U\star \mutilde^{({X},\mathbb F)})_t=\int_{\R{\ell}}U(t,x)\mutilde^{({X},\mathbb F)}(\{t\}\times\dx)$ by definition.
Let us also introduce the following convenient notation
\begin{align*}
\int_{\tau_1}^{\tau_2}\int_{\R{\ell}} U_s(x)\,\mutilde^{({X}, \mathbb F)}(\ds,\dx)
	:= 	U\star\mutilde^{({X},\mathbb F)}_{\tau_2}
	-	U\star\mutilde^{({X},\mathbb F)}_{\tau_1},
\end{align*}
where $\tau_1,\tau_2$ are $\mathbb F-$stopping times such that $0\leq \tau_1\leq \tau_2\leq \infty,$ $\Pm-a.s.$

\begin{remark}\label{rem:PurDisJum}
Observe that the canonical projections satisfy $\pi^i\in G_2(\mu^X,\Fil,\tau;\mathbb R)$, for every $i=1,\dots,\ell$.
Therefore, we can associate to the process $X$ the $\mathbb F-$martingale $(\pi^1\star \widetilde{\mu}^{(X,\mathbb F)},\dots ,\pi^\ell\star \widetilde{\mu}^{(X,\mathbb F)})\in\mathcal H^{2,d}(\mathbb F,\tau;\mathbb R^\ell).$ In case $X\in\mathcal H^{2}(\mathbb F,\tau;\mathbb R^\ell)$, it is clear that $X^d=(\pi^1\star \widetilde{\mu}^{(X,\mathbb F)},\dots ,\pi^\ell\star \widetilde{\mu}^{(X,\mathbb F)}),$
\emph{i.e.} the purely discontinuous part of the martingale $X$ is indistinguishable from  $(\pi^1\star \widetilde{\mu}^{(X,\mathbb F)},\dots ,\pi^\ell\star \widetilde{\mu}^{(X,\mathbb F)}).$
Henceforth, when $X\in\mathcal H^{2}(\mathbb F,\tau;\mathbb R^\ell)$, we will make no distinction between these two purely discontinuous martingales.
Moreover, assuming that $X$ is an $\mathbb F-$martingale, when we refer to the jump process $\Delta X^d$ we will mean the $\mathbb R^\ell-$valued process 
\begin{align}\label{def:JumpMartin}
\Delta X^d_{\cdot}:=\bigg(\int_{\mathbb R^\ell} \pi^1(x)\widetilde{\mu}^{(X,\mathbb F)}(\{\cdot\}\times\dx), \dots, \int_{\mathbb R^\ell} \pi^\ell(x)\widetilde{\mu}^{(X,\mathbb F)}(\{\cdot\}\times\dx) \bigg),
\end{align}
while, assuming that $X$ is simply an $\mathbb F-$adapted process, when we refer to the jump process $\Delta X$ we will mean the $\mathbb R^\ell-$valued process
\begin{align}\label{def:JumpProcess}
\Delta X_{\cdot}=\bigg(\int_{\mathbb R^\ell} \pi^1(x)\mu^X(\{\cdot\}\times\dx), \dots, \int_{\mathbb R^\ell} \pi^\ell(x)\mu^X(\{\cdot\}\times\dx) \bigg).
\end{align}
This subtle difference arises because the process $X$ is not quasi--left--continuous, hence the compensator $\nu^{(X,\mathbb F)}$ can have jumps.
In other words, the jumps of $\mu^X$ and $\widetilde{\mu}^{(X,\mathbb F)}$ are not identical, in general.
\end{remark}
\iftoggle{full}{\vspace{0.5em}
\begin{remark}
Let $X \in \mathcal H^{2}(\mathbb F,\tau;\mathbb R^\ell)$.
We can also associate an integer--valued random measure to the jumps of the martingale $X^d$, denoted by $\mu^{X^d}$, and then \eqref{def:JumpMartin} can be written in a similar form to \eqref{def:JumpProcess}, \textit{i.e.}
\begin{align}\label{def:JumpMartin-2}
\Delta X^d_{\cdot}=\bigg(\int_{\mathbb R^\ell} \pi^1(x)\mu^{X^d}(\{\cdot\}\times\dx), \dots, \int_{\mathbb R^\ell} \pi^\ell(x)\mu^{X^d}(\{\cdot\}\times\dx) \bigg).
\end{align}
\end{remark}

\begin{remark}\label{rem:NotationIVM}
We would like to clarify a subtle detail in our notation at this point, namely the difference between the mappings $*$ and $\star$.
Apart from the fact that they have different domains, let us comment on the way they are defined.
Let $U\in G_2(\mu^{X},\mathbb F,\tau;\R{})$, then
$$U*\mutilde^{({X},\mathbb F)}_{\cdot}=\int_{(0,\cdot]\times \R{\ell}}U(\omega,s,x)\mutilde^{({X},\mathbb F)}(\omega;\ds,\dx),$$
\textit{i.e.} this is the Lebesgue--Stieljes integral of $U$ with respect to the signed measure $\mutilde^{(X,\mathbb F)}$, for which the only information we have regarding its integrability is the square summability of its jumps.
Clearly this does not imply the finiteness of the process in any time interval.
On the contrary, by $U\star\mutilde^{({X},\mathbb F)}$ we denote the square--integrable purely discontinuous $\mathbb F-$martingale whose jump at each time $t$ is given by $\int_{\R{\ell}}U(t,x)\mutilde^{({X},\mathbb F)}(\{t\}\times\dx).$ 
A specific case where the two processes $U*\mutilde^{({X},\mathbb F)}_{\cdot}$ and $U\star\mutilde^{({X},\mathbb F)}_{\cdot}$ coincide is given by {\rm \cite[Proposition II.1.28]{jacod2003limit}} and corresponds to the finite variation case.
\end{remark}}{}

The space of real--valued square--integrable stochastic integrals with respect to $\mutilde^{({X},\mathbb F)}$ will be denoted by
\begin{align*}
\mathcal{K}^2\pair{\mu^{X},\mathbb F,\tau;\R{}}
	:=\set{U\star\mutilde^{({X},\mathbb F)},\; U\in G_2(\mu^{X},\mathbb F,\tau;\R{})}. 
\end{align*}
By \cite[Theorem II.1.33]{jacod2003limit}, or \citet[Theorem 11.21]{he1992semimartingale}, and the Kunita--Watanabe inequality, see {\it e.g.} \cite[Corollary 6.34]{he1992semimartingale}, we have
\begin{align*}
\E \big[ \langle U\star\mutilde^{({X},\mathbb F)}\rangle_{\tau}\big] <\infty, \ \textrm{ if and only if } \
U\in G_2\pair{\mu^{X},\mathbb F,\tau;\R{}},
\end{align*}
which enables us to define the following more convenient space
\begin{align*}
\Htwo{}{\mu^X,\mathbb F,\tau;\R{}}&
	:=\left\{ 
			U:\big(\widetilde{\Omega},\widetilde{\Pred}^\mathbb F\big) \longrightarrow \big(\R{},\mathcal B(\R{})\big),\ 
				\E \big[\langle U\star\mutilde^{({X},\mathbb F)}\rangle_{\tau}\big]<\infty
		\right\},
\end{align*}
and we emphasise that we have the direct identification $\Htwo{}{\mu^X,\mathbb F,\tau;\R{}}=G_2(\mu^{X},\mathbb F,\tau;\R{}).$

\subsubsection{Orthogonal decompositions}\label{sec:OrthDec}

We close this subsection with a reminder on orthogonal decompositions of square integrable martingales.

\begin{definition}
Let ${X}\in\mathcal{H}^2(\mathbb F,\tau;\R{\ell})$. 
${X}$ is said to possess the $\mathbb F-$\emph{predictable representation property} if
\begin{align*}
\mathcal{H}^2_0(\mathbb F,\tau;\R{}) = \mathcal{L}^2({X}^c,\mathbb F,\tau;\R{})\oplus \mathcal{K}^2(\mu^{X},\mathbb F,\tau;\R{}),
\end{align*}
where $\mathcal{H}^2_0\pair{\mathbb F,\tau;\R{}}: = \{M\in\mathcal{H}^2(\mathbb F,\tau;\R{}), M_0=0\}.$ 
In other words, for any $M\in\mathcal{H}^2(\mathbb F,\tau;\R{})$, there exists a pair $(Z,U)\in\Htwo{}{{X}^c,\mathbb F,\tau;\R{\ell}}\times\Htwo{}{\mu^X,\mathbb F,\tau;\R{}}$ such that
$$M_\cdot=M_0+Z\cdot{X}^c_\cdot+U\star\mutilde^{({X},\mathbb F)}_\cdot.$$
\end{definition}

In the sequel, we adapt the notation of \citet[Sections 13.2--3]{cohen2015stochastic}.
We associate the measure $M_{\mu}:(\widetilde{\Omega}, \G\otimes \borel{\Rp}\otimes\borel{\R{\ell}})\longrightarrow \Rp$ to a random measure $\mu$, which is defined as $M_{\mu}(B)=\E[ \mathds{1}_B*\mu_{\infty}]$. 
We will refer to $M_{\mu}$ as the \emph{Dol\'eans measure associated to $\mu.$}
If there exists an $\mathbb F-$predictable partition $(A_n)_{n\in\N}$ of $\widetilde{\Omega}$ such that $M_{\mu}(A_n)<\infty,$ for every $n\in\N,$ then we will say that $\mu$ is $\mathbb F-$predictably $\sigma-$integrable and we will denote it by $\mu\in\widetilde{\mathcal{A}}_{\sigma}(\mathbb F).$ 
For a sub--$\sigma$--algebra $\mathcal{A}$ of $ \G\otimes \borel{\Rp}\otimes\borel{\R{\ell}}$, the restriction of the measure $M_{\mu}$ to $(\widetilde{\Omega}, \mathcal{A})$ will be denoted by $M_{\mu}|_{\mathcal{A}}.$
Moreover, for $W:(\widetilde{\Omega},\G\otimes \borel{\Rp}\otimes\borel{\R{\ell}})\longrightarrow (\R{},\borel{\R{}})$, we define the random measure 
$W\mu$ as follows
$$(W\mu)(\omega;\ds,\dx):= W(\omega,s,x)\mu(\omega;\ds,dx).$$ 

\begin{definition}\label{CondFPredProj}
Let $\mu\in\widetilde{\mathcal{A}}_{\sigma}(\mathbb F)$ and 
$W:(\widetilde{\Omega},\G\otimes \borel{\Rp}\otimes\borel{\R{\ell}})\longrightarrow (\R{},\borel{\R{}})$ be such that $|W|\mu\in\widetilde{\mathcal{A}}_{\sigma}(\mathbb F).$
Then we define the \emph{conditional $\mathbb F-$predictable projection of $W$ on $\mu$}, denoted by 
$M_{\mu}\big[W|\widetilde{\mathcal{P}}^{\mathbb F}\big]$ as follows
$$M_{\mu}\big[W|\widetilde{\mathcal{P}}^{\mathbb F}\big]:=
	\frac{\ud M_{W\mu}|_{\widetilde{\mathcal{P}}^{\mathbb F}}}{\ud M_{\mu}|_{\widetilde{\mathcal{P}}^{\mathbb F}}}.$$ 
\end{definition}

The following definition is justified by \cite[Lemma III.4.24]{jacod2003limit}.

\begin{definition}
Let ${Y}\in\mathcal{H}^2(\mathbb F,\tau;\R{})$ and consider a triple $(Z,U,N)\in\mathbb{H}^2({X}^c,\mathbb F,\tau;\R{\ell})\times\mathbb{H}^2(\mu^X,\mathbb F,\tau;\R{})\times\mathcal{H}^2(\mathbb F,\tau;\R{})$ such that
\begin{align}\label{eq:orth-deco}
Y= Y_0 + Z\cdot {X}^c + U\star\mutilde^{({X},\mathbb F)} + {N},
\end{align}
with $\langle N,{X}^{c,i}\rangle = 0,$ for $i=1,\ldots,\ell$, and $M_{\mu^{X}}\big[ \Delta {N} |\widetilde{\mathcal{P}}^{\mathbb F} \big]=0.$
Then, \eqref{eq:orth-deco} is called the \emph{orthogonal decomposition of $\ {Y}$ with respect to $(X^c,\mu^X,\mathbb F).$}
\end{definition}

We conclude this subsection with a useful corollary, which must be preceded by the definition of the following space
\begin{align*}
\mathcal{H}^2(X^\perp,\mathbb F,\tau;\R{}) 		:=\Big\{N\in\mathcal{H}^2(\mathbb F,\tau;\R{}),\; N\mutorthog L,\; \forall\ L\in\mathcal{L}^2(X^c,\mathbb F,\tau;\R{})\oplus\mathcal{K}^2(\mu^{X},\mathbb F,\tau;\R{})\Big\}.
\end{align*}  

\begin{corollary}\label{cor:orthogSpace}
Let $Y\in\mathcal{H}^2(\mathbb F,\tau;\R{})$, $X\in\mathcal{H}^2(\mathbb F,\tau;\R{\ell})$ and $Y= Y_0 + Z\cdot {X}^c + U\star\mutilde^{({X},\mathbb F)} + {N}$ be the orthogonal decomposition of $\ Y$ with respect to $(X^c,\mu^X,\mathbb F)$.
Then $N\in \mathcal{H}^2(X^\perp,\mathbb F,\tau;\R{}).$ 
In particular $N\mutorthog X^i,$ for every $i=1,\ldots,\ell.$ 
\end{corollary}

\begin{proof}
This is immediate by \cite[Proposition I.4.15, Theorem III.4.5]{jacod2003limit} and \cite[Theorem 13.3.16]{cohen2015stochastic}.
\end{proof}

\subsection{\texorpdfstring{The Skorokhod space $(\D,\J_1)$ and convergence in the extended sense} 
			{The Skorokhod space (D, J 1) and convergence in the extended sense}}\label{sec:skor}

The natural path--space for an $\R{p\times q}-$valued process, which is adapted to some filtration, is the Skorokhod space 
\begin{align*}
\D([0,\infty);\R{p\times q}):=\set{\alpha:\Rp\longrightarrow \R{p\times q},\  \alpha \textrm{ is \cadlag }},
\end{align*}
which we equip with the Skorokhod $\J_1(\R{p\times q})-$topology, see \cite[Section VI.1b]{jacod2003limit}. 
We denote by $d_{\J_1(\R{p\times q})}$ the metric which is compatible with the $\J_1(\R{p\times q})-$topology.
\iftoggle{full}{Due to \cite[Comments VI.1.21-22]{jacod2003limit}, and since in the remainder of the paper we will have to distinguish between joint and separate convergence on products of Skorokhod spaces, we will always indicate the state space in our notation, for the sake of clarity. We remind the reader that $(\D([0,\infty);\R{p\times q}), d_{\J_1(\R {p\times q})})$ is a Polish space, see \cite[Theorem VI.1.14]{jacod2003limit}.}{}
We will denote by $d_{\lu}$ the metric on $\D([0,\infty);\R{p\times q})$ which is compatible with the topology of \emph{locally uniform} convergence, see \cite[Section VI.1a]{jacod2003limit}, and by $d_{\norm{\cdot}_\infty}$ the metric on $\D([0,\infty);\R{p\times q})$ which is compatible with the topology of \emph{uniform} convergence. 
Clearly $d_{\norm{\cdot}_\infty}$ is stronger than $d_{\lu}$. 
Moreover, it is well--known that $d_{\lu}$ is stronger than $d_{\J_1(\R{p\times q})}$, see \cite[Proposition VI.1.17]{jacod2003limit}.
In order to simplify notations, we set $\D^{p\times q}\!:=\D\pair{[0,\infty);\R{p\times q}}$ and, in case $p=q=1,$ $\D:=\D^1$. 
\iftoggle{full}{Moreover, to avoid any misunderstanding, we will not introduce any shorthand notation for the space $(\D([0,\infty);\R{}))^{p\times q}.$ }{}
We will postpone all proofs of the present section, except for the very short ones, to Appendix \ref{subsec:JointConv} for the sake of readability.

\iftoggle{full}{
\begin{definition} 
Let $(M^k)_{k\in\overline{\N}}$ be an arbitrary sequence such that $M^k$ is an $\R{p\times q}-$valued \cadlag process, for every $k\in\overline{\N}$.

\begin{enumerate}
\item The sequence $(M^k)_{k\in\N}$ converges \emph{in probability under the $\J_1(\R{p\times q})-$topology} to $M^\infty$ if 
$$\mathbb P\Big(d_{\J_1(\R{p\times q})}(M^k,M^\infty)>\varepsilon\Big)\xrightarrow[\hspace{0.2cm}k\to\infty\hspace{0.2cm}]{} 0,\ \text{ for every $\varepsilon>0$,}$$
and we denote\footnote{Notice that we omit the index associated to the convergence (\emph{i.e.} $k\longrightarrow\infty$). We will do the same in the remainder of the paper, when it is clear to which index the convergence refers.} it by $M^k\sipmulti{p\times q} M^\infty$.

\vspace{.5em}
\item Let $\vartheta\in[1,\infty)$. The sequence $(M^k)_{k\in\N}$ converges \emph{in $\mathbb{L}^{\vartheta}-$mean under the $\J_1(\R{p\times q})-$topology} to $M^\infty$ if
$$\E \Big[\big(d_{\J_1(\R{p\times q})}(M^k,M^\infty)\big)^{\vartheta}\Big]\xrightarrow[\hspace{0.2cm}k\to\infty\hspace{0.2cm}]{} 0,$$
and we denote it by $M^k\silmulti{p\times q}{\vartheta} M^\infty.$

\vspace{0.5em}
\item Analogously, we denote by $M^k\luip M^\infty$, resp. $M^k\luil{\vartheta} M^\infty$, the convergence \emph{in probability}, resp. \emph{in $\mathbb{L}^{\vartheta}-$mean}, \emph{under the locally uniform topology}.

\vspace{0.5em}
\item Let $(p_1,p_2,q_1,q_2)\in\N^4$.
	Moreover, let $(M^{k})_{k\in\overline{\N}}$ be a sequence of $\R{p_1\times q_1}-$valued and \cadlag processes and $(N^{k})_{k\in\overline{\N}}$ be a sequence of $\R{p_2\times q_2}-$valued and \cadlag processes. 
For $\vartheta_1,\vartheta_2\in[1,\infty)$ we will write
\begin{align*}
(M^{k}, N^{k})
\xrightarrow{\hspace{0.2cm} \big(\J_1(\R{p_1\times q_1}\times\R{p_2\times q_2}),\mathbb{L}^{\vartheta_1}(\D^{p_1\times q_1})\times\mathbb{L}^{\vartheta_2}(\D^{p_2\times q_2})\big) \hspace{0.2cm}} 
(M^\infty, N^\infty),
\end{align*}
if the following convergences hold
\begin{align*}
&(M^{k}, N^{k})
	\xrightarrow{\hspace{0.2cm} \big(\J_1(\R{p_1\times q_1}\times\R{p_2\times q_2}),\Pm)\big) \hspace{0.2cm}} 
	(M^\infty, N^\infty),\ M^{k}\silmulti{q_1}{\vartheta_1}M^\infty,	
	 \; \text{and}\;		
	N^k\silmulti{q_2}{\vartheta_2}N^\infty.
 \end{align*}
\end{enumerate}
\end{definition}}{For an arbitrary sequence $(M^k)_{k\in\overline{\N}}$ such that $M^k$ is an $\R{p\times q}-$valued \cadlag process, for every $k\in\overline{\N}$, we denote by $M^k\sipmulti{p\times q} M^\infty$ the convergence in probability under the $\J_1(\R{p\times q})-$topology to $M^\infty$, by $M^k\silmulti{p\times q}{\vartheta} M^\infty$ the convergence in $\mathbb{L}^{\vartheta}-$mean under the $\J_1(\R{p\times q})-$topology to $M^\infty$, for $\vartheta\in[1,\infty)$, and we use similar notations for convergence in probability, and in $\mathbb{L}^{\vartheta}-$mean, under the locally uniform topology. We also extend naturally these notations to convergence of pairs of processes $(M^{k}, N^{k})$.}

\vspace{0.5em}
Let us now introduce notions related to the convergence of $\sigma-$fields and filtrations. 
We will need the filtrations to be indexed by $[0,\infty]$, hence, given a filtration $\mathbb F:=(\mathcal F_t)_{t\ge 0}$, we define the $\sigma-$algebra $\mathcal F_\infty$ by using the convention 
$$\mathcal F_\infty:=\mathcal F_{\infty-}=\bigvee_{t\geq 0} \mathcal F_t.$$

We recall also the following notation for every sub$-\sigma-$field $\mathcal F$ of $\mathcal G$ and $\vartheta\in[1,\infty)$
$$\mathbb L^{\vartheta}(\Omega,\mathcal F,\mathbb P;\mathbb R^q):=\left\{\xi,\ \text{$\mathbb R^q-$valued and $\mathcal F-$measurable such that } \mathbb E[|\xi|^{\vartheta}]<\infty\right\}.$$

\begin{definition}\label{FiltrationConvergence}
${\rm (i)}$ A sequence of $\sigma-$\emph{algebrae} $(\mathcal F^k)_{k\in\mathbb{N}}$ converges \emph{weakly} to the $\sigma-$algebra $\mathcal F^\infty$ if, for every $\xi\in\mathbb{L}^1(\Omega,\mathcal F^\infty,\Pm;\mathbb R)$, we have
		$$
		\mathbb E[ \xi | \mathcal F^k]\inp\mathbb E[\xi|\mathcal F^\infty].
		$$
		We denote the weak convergence of $\sigma-$algebrae by  $\mathcal F^k\weakFil \mathcal F^\infty$.

\vspace{0.5em}
${\rm (ii)}$ A sequence of filtrations $\big(\mathbb F^k:=(\mathcal F^k_t)_{t\geq 0}\big)_{k\in\mathbb{N}}$ converges \emph{weakly} to $\mathbb F^\infty:=(\mathcal F^\infty_t)_{t\geq 0}$,
		if, for every $\xi\in\mathbb{L}^1(\Omega,\mathcal F^\infty_\infty,\Pm;\mathbb R)$, we have
		\[
		\mathbb E[ \xi  | \mathcal F^k_\cdot]\underset{}{\sipmulti{}}\mathbb E[\xi|\mathcal F^\infty_\cdot].
		\]
		We denote the weak convergence of the filtrations by $\mathbb F^k\weakFil \mathbb F^\infty$.

\vspace{0.5em}
${\rm (iii)}$ Consider the sequence $\big( (M^k, \mathbb F^k)\big)_{k\in\overline{\N}}$, where $M^k$ is an $\R{q}-$valued \cadlag process and $\mathbb F^k$ is a filtration, for any $k\in\overline{\N}$.
The sequence $\big((M^k, \mathbb F^k)\big)_{k\in\mathbb{N}}$ converges \emph{in the extended sense} to $(M^\infty, \mathbb F^\infty)$ if for every $\xi\in\mathbb{L}^1(\Omega,\mathcal F^\infty_\infty,\Pm;\mathbb R)$,
		\begin{equation}\label{extended}
		\bigg(
			\begin{array}{c}
				M^{k}\\ \mathbb E[ \xi | \mathcal F^k_\cdot]
			\end{array}
		\bigg)
		\sipmulti{q+1}
		\bigg(
			\begin{array}{c}
				M^\infty\\ \mathbb E[ \xi | \mathcal F^\infty_\cdot]
			\end{array}
		\bigg).
		\end{equation}
		We denote the convergence in the extended sense by
		$\left(M^k, \mathbb F^k\right)\extsense \left(M^\infty, \mathbb F^\infty\right)$.
\end{definition}

\iftoggle{full}{\begin{remark}
For the definition of weak convergence of filtrations, we could have used only random variables $\xi$ of the form $\mathds{1}_A$, for $A\in\mathcal F^\infty_\infty$. 
Indeed, the two definitions are equivalent, see {\rm \citet[Remark 1.1)]{coquet2001weak}}.
\end{remark}}{}

\iftoggle{full}{The following result, which is due to \citet[Theorem 7.4]{hoover1991convergence}, provides a sufficient condition for weak convergence of $\sigma-$\emph{algebrae} which are generated by random variables.
\begin{example}
Let $(\xi^k)_{k\in\overline{\N}}$ be a sequence of random variables such that  $\xi^k\inp\xi^\infty.$
Then the convergence $\sigma(\xi^k)\weakFil\sigma(\xi^\infty)$ holds, where $\sigma(\psi)$ denotes the $\sigma-$\emph{algebra} generated by the random variable $\psi.$
\end{example} }{}

In the \iftoggle{full}{next}{following} example, which is \cite[Proposition 2]{coquet2001weak}, a sufficient condition for the weak convergence of the natural filtrations of stochastic processes is provided.
\begin{example}\label{example:PII}
Let $M^k$ be a process with independent increments, for every $k\in\overline{\N}$. 
If $M^k\sipmulti{q} M^\infty$, then $\Fil^{M^k}\!\weakFil\Fil^{M^\infty}$.
\end{example} 

In the remainder of this section, we fix an arbitrary sequence of filtrations $(\mathbb F^k)_{k\in\overline{\N}}$ on $(\Omega,\mathcal G,\Pm)$ (recall Footnote \ref{UCFil}), with $\mathbb F^k:=(\mathcal F^k_t)_{t\geq 0},$ and an arbitrary sequence $(M^k)_{k\in\overline{\N}}$, where $M^k$ is an $\R{q}-$valued, uniformly integrable, $\mathbb F^k-$martingale, for every $k\in\overline{\N}.$
Then, it is well known that the random variables $M^k_\infty:=\lim_{t\to\infty} M^k_t$ are well--defined $\Pm-a.s.$, and 
$M^k_\infty\in\mathbb{L}^1(\Omega,\mathcal F^k_\infty,\Pm;\mathbb R^q)$, for $k\in\overline{\N}$; see \cite[Theorem I.1.42]{jacod2003limit}.

\vspace{0.5em}
Next, we would like to discuss how to deduce the extended convergence of martingales and filtrations from individual convergence results. 
Such properties have already been obtained by \citet[Proposition 1.(iii)]{memin2003stability}, where he refers to \citet[Proposition 7]{coquet2001weak} for the proof.
However, the authors in \cite{coquet2001weak} proved the result under the additional assumption that the processes are adapted to their natural filtrations. 
Moreover, they consider a finite time horizon $T$, which gives the time point $T$ a special role for the $\J_1(\R{})-$topology on 
$\D([0,T];\R{})$, see also \cite[Remark VI.1.10]{jacod2003limit}. 
In addition, in \cite[Remark 1.2)]{coquet2001weak}, the convergence $M^k_\infty\xrightarrow{\hspace{0.2cm}\mathbb{L}^{1}(\Omega,\mathcal F^\infty_\infty,\Pm;\mathbb R^q)\hspace{0.2cm}} M^\infty_\infty$ is assumed, although it is not necessary (note that we have translated their results into our notation). 
This is restrictive, in the sense that they have to assume in addition the $\mathcal F^\infty_\infty-$measurability of $M^k_\infty,$ for each $k\in\N.$

\vspace{0.5em}
We present below, for the sake of completeness, the statement \iftoggle{full}{and proof}{} of the aforementioned results for the infinite time horizon case, under the condition 
$M^k_\infty\xrightarrow{\hspace{0.2cm}\mathbb{L}^{1}(\Omega,\mathcal G,\Pm;\mathbb R^q)\hspace{0.2cm}} M^\infty_\infty$. 
\iftoggle{full}{}{The proofs, aside from a very short one, can be found in Appendix A.1.2 in the extended version of the present article \cite{papapantoleon2018stability}.}
 
\begin{proposition}\label{MeminProp1}
Assume the convergence $M_\infty^k\xrightarrow{\hspace{0.2cm}\mathbb{L}^1(\Omega,\mathcal G,\Pm;\mathbb R^q)\hspace{0.2cm}} M^\infty_\infty$ holds.
Then, the convergence $\mathbb F^k\xrightarrow{\hspace{0.2cm}\w\hspace{0.2cm}}\mathbb F^\infty$ is equivalent to the convergence $(M^k, \mathbb F^k)\extsense (M^\infty, \mathbb F^\infty)$.
\end{proposition}

The following two results, which are essentially \cite[Theorem 11, Corollary 12]{memin2003stability}, constitute the cornerstone for the convergence in the extended sense.
Here we state and prove them in the multi--dimensional case.
Before we proceed, let us recall some further definitions.
An $\mathbb F-$adapted process $M$ is called $\mathbb F-$\emph{quasi--left--continuous} if $\Delta M_{\sigma}=0$, $\mathbb P-a.s.$, for every $\mathbb F-$predictable time $\sigma.$  
An $\mathbb F-$adapted process $S$ is called an $\mathbb F-$\emph{special semimartingale} if $S=S_0 + M + A,$ where $S_0$ is finite--valued and $\mathcal F_0-$measurable, $M$ is a local $\mathbb F-$martingale with $M_0=0$ and $A$ is an $\mathbb F-$predictable, finite variation process with $A_0=0$; see \cite[Definition I.4.21]{jacod2003limit}. 
This decomposition of an $\mathbb F-$special semimartingale is unique, and for this reason we will call it the $\mathbb F-$canonical decomposition of $S$.  
For a process $A$ of finite variation, we denote by ${\rm Var} (A)$ the \emph{$($total$)$ variation process} of $A$, \emph{i.e.} ${\rm Var}(A)_t(\omega)$ is the total variation of the function $\Rp \ni s\longmapsto A_s(\omega)\in\Rp$ in the interval $[0,t]$.
For $A\in\D^{p\times q},$ we denote by ${\rm Var}(A)\in\D^{p\times q}$ the process for which ${\rm Var}(A)^{ij}:={\rm Var}(A^{ij}),$ for $i=1,\dots,p$ and $j=1,\ldots,q.$ 

\begin{theorem}\label{MeminTheorem}
Let $(S^k)_{k\in\overline{\N}}$ be a sequence of $\R{q}-$valued $\mathbb F^k-$special semimartingales with $\mathbb F^k-$canonical decomposition $S^k=S^k_0+M^k+A^k,$ for every $k\in\overline{\N}.$
Assume that $S^\infty$ is $\Fil^\infty-$quasi--left--continuous and the following properties hold
\begin{enumerate}[label={\rm (\roman*)}]
\item\label{MeminTh-i} the sequence $\big([S^{k,i}]_\infty^{1/2}\big)_{k\in\overline{\N}}$ is uniformly integrable, for every $i=1,\ldots,q$,

\vspace{0.3em}
\item\label{MeminTh-ii} the sequence $(\Vert{\rm Var}(A^k)_\infty\Vert_1)_{k\in\N}$ is tight,

\vspace{0.3em}
\item\label{MeminTh-iii} the extended convergence $(S^k,\mathbb F^k)\extsense (S^\infty,\mathbb F^\infty)$ holds.
\end{enumerate}
Then 
\[
	(S^k, M^k, A^k)\sipmulti{q\times3} (S^\infty,M^\infty,A^\infty).
\]
\end{theorem}

\begin{proof}
By \cite[Theorem 11]{memin2003stability}, we obtain for every $i=1,\ldots,q$ the following convergence 
\[
	(S^{k,i}, M^{k,i}, A^{k,i})\sipmulti{3} (S^{\infty,i},M^{\infty,i},A^{\infty,i}). 	
\]
Then, by assumption $S^k\sipmulti{q} S^\infty,$ and using Corollary \ref{cormultiskorokhod} and Remark \ref{rem:cormultiskorokhod} we obtain the required result.
\end{proof}

\begin{theorem}\label{MeminCorollary}
Let $M^k\in\mathcal{H}^2(\mathbb F^k,\infty;\R{q})$ for any $k\in\overline{\N}$ and $M^\infty$ be $\mathbb F^\infty-$quasi--left--continuous.
If the following convergences hold
\begin{align*}
	(M^k,\mathbb F^k)\extsense (M^\infty,\mathbb F^\infty)
		\ \text{ and }\
	M^k_\infty\xrightarrow{\hspace{0.2cm}\mathbb{L}^2(\Omega,\mathcal G,\Pm;\mathbb R^q)\hspace{0.2cm}}M^\infty_\infty,
\end{align*}
then
\begin{enumerate}[label={\rm(\roman*)}]
\item\label{MeminCorollaryJ1P}
	$(M^k, [M^k], \langle M^k\rangle)	
		\xrightarrow{\hspace{0.2cm} \left(\J_1(\R{q}\times\R{q\times q}\times\R{q\times q}),\Pm\right)\hspace{0.2cm}} 
	(M^\infty, [M^\infty], \langle M^\infty\rangle)$,
	
	\vspace{0.3em}
\item\label{BracketConvinL} 
	for every $i=1,\ldots,q$, we have
	\begin{align*}
		[M^{k,i}]_\infty 
			\underset{k\to\infty}{\xrightarrow{\hspace{0.2cm}\mathbb{L}^1(\Omega,\mathcal G,\Pm;\mathbb R^{})\hspace{0.2cm}}}
		[M^{\infty,i}]_\infty\
	\text{ and }\ 
		\langle M^{k,i}\rangle_\infty 
			\underset{k\to\infty}{\xrightarrow{\hspace{0.2cm}\mathbb{L}^1(\Omega,\mathcal G,\Pm;\mathbb R^{})\hspace{0.2cm}}}
		\langle M^{\infty,i}\rangle_\infty.
	\end{align*}
\end{enumerate}
\end{theorem}

We conclude this subsection with the following technical lemma which will be of utmost importance for us in the proof of our robustness result for martingale representations.

\begin{lemma}\label{UIplusL2Bounded}
Let $(L^k)_{k\in\N}$ be a sequence of $\R{p}-$valued processes such that $\big({\rm Tr}\big[[L^k]_\infty\big] \big)_{k\in\N}$ is uniformly integrable and 
$(N^k)_{k\in\N}$ be a sequence of $\R{q}-$valued processes such that $\big({\rm Tr}\big[[N^k]_\infty\big] \big)_{k\in\N}$ is bounded in $\mathbb L^1(\Omega,\mathcal G,\Pm)$, 
\emph{i.e.} $\sup\limits_{k\in\N}\E\big[{\rm Tr}\big[[N^k]_\infty\big]\big]<\infty$.
Then $\big(\big\Vert{\rm \Var}\big([L^k,N^k]\big)_\infty \big\Vert_1\big)_{k\in\N}$ is uniformly integrable.  
\end{lemma}

\section{Stability of martingale representations}\label{sec:RobMartRepSection}	

\subsection{Framework and statement of the main theorem}\label{subsec:RoMRFrame}

We start by presenting and discussing the main assumptions that will be used throughout this section.
Let us fix an arbitrary sequence of \cadlag $\R{\ell}-$valued processes $(X^k)_{k\in\overline{\N}}$ for which we assume that
\begin{align}\label{SqIntofIVM}
\sup_{k\in\overline\N}\mathbb E\bigg[\int_{(0,\infty)\times\mathbb R^\ell} |x|^2\,\mu^{X^k}(\ds,\dx)\bigg]<\infty.
\end{align}
Then we fix an arbitrary sequence of filtrations $(\FilG^k)_{k\in\overline{\N}}$ (recall Footnote \ref{UCFil}) with $\FilG^k:=(\G^k_t)_{t\geq 0},$ for every $k\in\overline{\N}$, on the probability space $(\Omega,\mathcal G,\mathbb P)$, and an arbitrary sequence of real--valued random variables $(\xi^k)_{k\in\overline{\N}}$.

\vspace{0.5em}
The following assumptions will be in force throughout this section.
\begin{enumerate}[label=\textup{\textbf{(M\arabic*)}},itemindent=1.15cm]
\item\label{MFilqlc} 	
	The filtration $\FilG^{\infty}$ is quasi--left--continuous and the process $X^\infty$ is $\mathbb G^\infty-$quasi--left--continuous.
\item\label{MSI}
	The process $X^k\in\mathcal H^2(\FilG^k,\infty;\R{\ell})$, for every $k\in\overline{\N}.$
	Moreover
	$X^k_{\infty} {\xrightarrow{\hspace*{0.2cm}\mathbb{L}^2(\Omega,\mathcal{G},\mathbb P;\mathbb R^\ell)\hspace*{0.2cm}}} X^{\infty}_{\infty}.$
\item\label{MXWPRP}
	{The martingale $X^{\infty}$ possesses the ${\FilG^{\infty}}-$predictable representation property.} \phantom{$\xrightarrow{L^2}$}
\item\label{MFilweak}
	The filtrations converge weakly, \textit{i.e.}
	$\FilG^k\weakFil \FilG^{\infty}$. \phantom{$\xrightarrow{L^2}$}
\item\label{Mfinalrv} 	
	The random variable $\xi^k\in \mathbb{L}^2(\Omega,\mathcal{G}^k_{\infty},\mathbb P;\mathbb R)$, for every $k\in\overline{\N}$, and 
			$\xi^k \xrightarrow{\hspace*{0.2cm}\mathbb{L}^2(\Omega,\mathcal{G},\mathbb P;\mathbb R)\hspace*{0.2cm}} \xi^{\infty}.$
\end{enumerate}

\begin{remark}\label{EquivalentExtended}
In view of {\rm Proposition \ref{MeminProp1}}, conditions {\rm \ref{MSI}} and {\rm \ref{MFilweak}} imply that 
$$(X^k,\FilG^k)\extsense (X^\infty,\FilG^\infty).$$
\end{remark}

\begin{remark}\label{rem:AlternCond}
In {\rm\ref{Mfinalrv}}, we have imposed an additional measurability assumption for the sequence of random variables $(\xi^k)_{k\in\overline{\N}}$, since we require that $\xi^k$ is $\mathcal G^k_\infty-$measurable for any $k\in\overline{\N}$, instead of just being $\mathcal G-$measurable. 
We could spare that additional assumption at the cost of a stronger hypothesis in {\rm\ref{MFilweak}}, namely that the weak convergence of the $\sigma-$algebrae
$$\G^k_{\infty}\weakFil \G^{\infty}_{\infty},$$
holds in addition. 
To sum up, the pair {\rm\ref{MFilweak}} and {\rm\ref{Mfinalrv}} can be substituted by the following
\begin{enumerate}[itemindent=1.25cm]
\item[{\bf(M4$^{\prime}$)}] The filtrations converge weakly as well as the final $\sigma-$algebrae, that is $\FilG^k\weakFil \FilG^{\infty}$ and $\G^k_{\infty}\weakFil \G^{\infty}_{\infty}.$
\item[{\bf(M5$^{\prime}$)}] The sequence $(\xi^k)_{k\in\overline \N}\subset \mathbb{L}^2(\Omega,\mathcal{G},\mathbb P;\mathbb R)$ and satisfies 
			$\xi^k \xrightarrow{\hspace*{0.2cm}\mathbb{L}^2(\Omega,\mathcal{G},\mathbb P;\mathbb R)\hspace*{0.2cm}} \xi^{\infty}.$
\end{enumerate}
\end{remark}

In the sequel we are going to update our notation as follows: for $k\in\overline{\N}$, $X^{k,c}$, respectively $X^{k,d}$, will denote the continuous part, respectively the purely discontinuous part, of the martingale $X^k$.
Moreover, for $i=1,\dots,\ell$, $X^{k,c,i}$, respectively $X^{k,d,i}$, will denote the $i-$th element of the continuous part, respectively the $i-$th element of the purely discontinuous part, of the martingale $X^k$.
An $\Fil-$predictable process $A$ will be denoted by $A^{\Fil}$, whenever the coexistence of several filtrations may create confusion. 

\begin{theorem}\label{RobMartRep}
Let conditions {\rm \ref{MFilqlc}--\ref{Mfinalrv}} hold and define the $\mathbb G^k-$martingales	$Y^{k} 	:= \E[\,\xi^{k}|\, \G^{k}_{\cdot}],$ for $k\in\overline{\mathbb{N}}. $
The orthogonal decomposition of $Y^{\infty}$ with respect to $(X^{\infty,c},\mu^{X^{\infty,d}},\FilG^{\infty})$ 
\begin{align*}
	Y^{\infty} 	&= Y^{\infty}_0 + Z^{\infty}\cdot X^{\infty,c}\ + U^{\infty}\star\mutilde^{(X^{\infty,d},\mathbb G^\infty)},
\end{align*}
and the orthogonal decomposition of $Y^k$ with respect to $(X^{k,c},\mu^{X^{k,d}}, \FilG^k)$
\begin{align*}
	Y^k &= Y^k_0 + Z^k\cdot X^{k,c} + U^k\star \mutilde^{(X^{k,d},\mathbb G^k)} + N^k, \text{ for } k\in\N,
\end{align*}
satisfy the following convergences 
\begin{align}\label{RobMartRepi}
\big(Y^k, Z^k\cdot X^{k,c} + U^k \star\mutilde^{(X^{k,d},\mathbb G^k)}, N^k \big)
&\silmulti{3}{2}
	\big(Y^{\infty}, Z^{\infty}\cdot X^{\infty,c} + U^{\infty}\star\mutilde^{(X^{\infty,d},\mathbb G^\infty)}, 0 \big),\\
\big(\langle Y^k\rangle, \langle Y^k,X^k\rangle,  \langle N^k\rangle \big)&
\xrightarrow{ \hspace{.3cm} \left(\J_1(\R{}\times\R\ell\times\R{}),\mathbb L^1\right) \hspace{.3cm} }
\big(\langle Y^{\infty}\rangle, \langle Y^{\infty},X^{\infty}\rangle,  0\big), \label{RobMartRepii}
\end{align}
where $\langle Y^k,X^k\rangle := (\langle Y^k,X^{k,1}\rangle,\dots,\langle Y^k,X^{k,\ell}\rangle)^\top$ for all $k\in\overline\N$.
\end{theorem}

\subsection{Examples and applications}\label{subsec:Ex_and_Appl}

In order to apply the above result in a concrete scenario, we need to check that Assumptions \ref{MFilqlc}--\ref{Mfinalrv} are satisfied.
The input data would then be a random variable $\xi^\infty$, a martingale $X^\infty$ and a filtration $\mathbb G^\infty$, and we assume they satisfy \ref{MFilqlc} and \ref{MXWPRP}.
Moreover, we can construct sequences $(\xi^k)_{k\in\N}$ and $(X^k)_{k\in\N}$ such that \ref{MSI} and \ref{Mfinalrv} are also satisfied.
Therefore, what remains to be shown and is not trivial, is the weak convergence of the filtrations, \textit{i.e.} \ref{MFilweak}.

\smallskip
The following two cases describe situations where we can easily check that this condition is satisfied.
\begin{itemize}[leftmargin=*]
\item According to \citet[Proposition 2]{coquet2001weak}, which we have stated as \cref{example:PII}, if $X^k$ is a martingale with independent increments that converges to $X^\infty$, and the filtrations $\mathbb G^k$ and $\mathbb G^\infty$ are the natural filtrations generated by the respective martingales, then \ref{MFilweak} is automatically satisfied.
\item According to \citet[Theorem 1]{coquet1998stability}, if $X^\infty$ is a \cadlag Markov process, $X^k$ is a discretization of $X^\infty$, and $\mathbb G^k$ and $\mathbb G^\infty$ are the natural filtrations generated by the respective martingales, then \ref{MFilweak} is again automatically satisfied.
\end{itemize}

\vspace{0.5em}
The following two sub--sub--sections provide two corollaries of \cref{RobMartRep} that are relevant for applications, in particular for numerical schemes.
 
\subsubsection{The case of processes with independent increments}

In this sub--sub--section, we focus on processes with independent increments, and we are interested in convergence in law, which is the relevant convergence for numerical schemes, such as the Euler--Monte Carlo method.
Convergence results for numerical schemes typically involve stochastic processes that are defined on distinct probability spaces, while the very definition of weak convergence of filtrations in \ref{MFilweak} requires that processes are defined on the same space.
In order to reconcile these opposing facts, we will work with the natural filtration of stochastic processes with independent increments and, in this case, the weak convergence of the filtrations follows from the results of \citet[Proposition~2]{coquet2001weak}.
Then, as a corollary of the main theorem, we can show that the convergence results in \eqref{RobMartRepi} and \eqref{RobMartRepii} hold also in law.

\vspace{.5em}
Let us set the framework for the results that follow.
Let $(X^k)_{k\in\overline{\mathbb N}}$ be a sequence of $\mathbb R^\ell-$valued \cadlag processes, where $X^k$ is defined on the space $(\Omega^k, \mathcal G^k, \mathbb P^k)$ for each $k\in\overline{\mathbb N}$, and $(\xi^k)_{k\in\overline{\mathbb N}}$ be a sequence of real--valued random variables, where each $\xi^k$ is defined on $(\Omega^k, \mathcal G^k, \mathbb P^k)$ for each $k\in\overline{\mathbb N}$.
Moreover, we assume that 
\begin{align}\label{WeakSchem:SqIntofIVM}
\sup_{k\in\overline{\mathbb N}} \mathbb E^k\bigg[\int_{(0,\infty)\times\mathbb R^\ell} |x|^2\,\mu^{X^k}(\ds,\dx)\bigg]<\infty,
\end{align}
where by $\mathbb E^k[\cdot]$ we have denoted the expectation under $\mathbb P^k$, for every $k\in\overline{\mathbb N}$.
We will denote analogously by $\mathbb E^k[\cdot|\mathcal F_\cdot]$ the conditional expectation with respect to (an element of) a filtration $\mathbb F$ under the measure $\mathbb P^k$, for every $k\in\overline{\mathbb N}$. 
Moreover, we will denote the set of $\mathbb R^q-$valued and square--integrable $\mathbb F^{X^k}-$martingales by
$\mathcal H^2(\mathbb F^{X^k},\infty;\mathbb R^q)$, for every $k\in\overline{\mathbb N}$, where we have notationally suppressed the dependence on $\Omega^k$.

\vspace{0.5em}
The following assumptions will be in force throughout this sub--sub--section.
\begin{enumerate}[label={\rm\bf(W\arabic*)},itemindent=.0cm,leftmargin=!]
\item\label{Weak-MFilqlc} The filtration $\mathbb F^{X^\infty}$ is quasi--left--continuous.
\item\label{Weak-MSI} $X^k\in\mathcal H^2(\mathbb F^{X^k},\infty;\R{\ell})$, for every $k\in\overline{\N}$.
			Moreover, $X^k\xrightarrow{\hspace{0.2em}\mathcal L\hspace{0.2em}} X^\infty \text{ in }\mathbb D([0,\infty);\mathbb R^\ell)$
			as well as 
			$X^k_{\infty} \xrightarrow{\hspace*{0.2cm}\mathcal L\hspace*{0.2cm}} X^{\infty}_{\infty}$ in $\mathbb R^\ell$, where  
			$(|X^k_\infty|^2)_{k\in\overline{\mathbb N}}$ is in addition uniformly integrable.
\item\label{Weak-MXWPRP} The martingale $X^{\infty}$ possesses the ${\mathbb F^{X^\infty}}-$predictable representation property.
\item\label{Weak-MFilweak} The process $X^k$ has independent increments relative to the filtration $\mathbb F^{X^k}$, for every $k\in\overline{\mathbb N}$.
\item\label{Weak-Mfinalrv} The random variable $\xi^k\in \mathbb{L}^2(\Omega^k,\mathcal{F}^{X^k}_{\infty},\mathbb P^k;\mathbb R)$, for every $k\in\overline{\N}$, and is such that $(|\xi^k|^2)_{k\in\overline{\mathbb N}}$ is uniformly integrable and $\xi^k \xrightarrow{\hspace*{0.2cm}\mathcal L\hspace*{0.2cm}} \xi^{\infty}.$
\end{enumerate}

\begin{corollary}\label{corr:PII}
Let conditions \ref{Weak-MFilqlc}--\ref{Weak-Mfinalrv} hold and define the martingales	$Y^{k} 	:= \mathbb E^k[\,\xi^{k}|\, \mathcal F^{X^k}_{\cdot}],$ for $k\in\overline{\mathbb{N}}. $
The orthogonal decomposition of $Y^{\infty}$ with respect to $(X^{\infty,c},\mu^{X^{\infty,d}},\mathbb F^{X^\infty})$ 
\begin{align*}
	Y^{\infty} 	&= Y^{\infty}_0 + Z^{\infty}\cdot X^{\infty,c}\ + U^{\infty}\star\mutilde^{(X^{\infty,d},\mathbb F^{X^\infty})},
\end{align*}
and the orthogonal decomposition of $Y^k$ with respect to $(X^{k,c},\mu^{X^{k,d}},\mathbb F^{X^k})$
\begin{align*}
	Y^k &= Y^k_0 + Z^k\cdot X^{k,c} + U^k\star \mutilde^{(X^{k,d},\mathbb F^{X^k})} + N^k, \text{ for } k\in\N,
\end{align*}
satisfy the following convergences 
\begin{align}\label{Weak-RobMartRepi}
	\big(Y^k, Z^k\cdot X^{k,c} + U^k \star\mutilde^{(X^{k,d},\mathbb F^{X^k})}, N^k \big)
		\xrightarrow{\hspace{0.2cm}\mathcal L \hspace{0.2cm}}
	\big(Y^{\infty}, Z^{\infty}\cdot X^{\infty,c} + U^{\infty}\star\mutilde^{(X^{\infty,d},\mathbb F^{X^\infty})}, 0 \big),
\end{align}
\begin{align}\label{Weak-RobMartRepii}
	\big(\langle Y^k\rangle, \langle Y^k,X^k\rangle,  \langle N^k\rangle \big)
		\xrightarrow{\hspace{0.2cm}\mathcal L \hspace{0.2cm}}
	\big(\langle Y^{\infty}\rangle, \langle Y^{\infty},X^{\infty}\rangle,  0\big).
\end{align}
\end{corollary}

\iftoggle{full}{

The proof is deferred to \cref{Proof_Corr_main_thm}.

}{

\begin{proof}
See Appendix A.4 in the extended version of the present article \cite{papapantoleon2018stability}.
\end{proof}

}

\subsubsection{A stronger version of the main theorem}

In this sub--sub--section, we strengthen the main theorem in the following sense: if we assume that there exist two sequences that converge to the continuous and the purely discontinuous part of the limiting martingale, then the angle brackets of $Y^k$ with respect to these sequences converge to the angle brackets of $Y^\infty$ with respect to the continuous and the purely discontinuous part of the limiting martingale.
This framework is very useful when considering discrete--time approximations of continuous--time processes.
To this end, we need to allow the It\=o integrator to exhibit jumps, \emph{i.e.} it should not necessarily be a continuous martingale.
Therefore, we also need to generalize the notion of orthogonal decompositions which was described in Sub--sub--subsection \ref{sec:OrthDec}.

\begin{definition}\label{def:OrthogDecomp}
Let $\mathbb F$ be a filtration, $(X^\circ,X^\natural)\in\mathcal H^2(\mathbb F,\infty;\mathbb R^m)\times\mathcal H^{2,d}(\mathbb F,\infty;\mathbb R^n)$ and $Y\in\mathcal H^2(\mathbb F,\infty;\mathbb R)$.
The decomposition
\begin{align*}
Y= Y_0 + Z\cdot X^\circ + U\star \widetilde{\mu}^{X^\natural} + N,
\end{align*}
where the equality is understood componentwise, will be called the \emph{orthogonal decomposition of $Y$ with respect to $(X^\circ,X^\natural)$} if
\begin{enumerate}[label=$(\roman*)$]
  \item $Z\in\mathbb H^2(X^\circ,\mathbb F,\infty;\mathbb R^m)$ and $U\in\mathbb H^2(\mu^{X^\natural},\mathbb F,\infty;\mathbb R),$
  \item $Z\cdot X^\circ \mutorthog U\star \widetilde{\mu}^{X^\natural}$, 
  \item $N\in\mathcal H^2(\mathbb F,\infty;\mathbb R)$ with $\langle N,X^{\circ}\rangle^{\mathbb F}=0$ and $M_{\mu^{X^\natural}}[\Delta N|\widetilde{\mathcal P}^\mathbb F]=0.$
\end{enumerate}
\end{definition}

The following results will allow us to obtain the orthogonal decomposition as understood in Definition \ref{def:OrthogDecomp}.
Their proofs can be found in \cite[Appendix A]{papapantoleon2016existence}.
For the statement of the following results we will fix an arbitrary filtration $\mathbb F$.
\begin{lemma}\label{lem:StochIntOrthog}
Let $(X^\circ,X^\natural)\in\mathcal H^2(\mathbb F,\infty;\mathbb R^m)\times\mathcal H^{2,d}(\mathbb F,\infty;\mathbb R^n)$ with $M_{\mu^{X^\natural}}[\Delta X^\circ|\widetilde{\mathcal P}^{\mathbb F}]=0$, where the equality is understood componentwise.
Then, for every $Y^\circ\in\mathcal L^2(X^\circ,\mathbb F,\infty;\mathbb R)$, $Y^\natural\in\mathcal K^2(\mu^{X^\natural},\mathbb F,\infty;\mathbb R)$, we have $\langle Y^\circ,Y^\natural\rangle=0.$
In particular, $\langle X^\circ, X^\natural\rangle=0$.
\end{lemma}

In view of Lemma \ref{lem:StochIntOrthog}, we can provide in the next proposition the desired orthogonal decomposition of a martingale $Y$ with respect to a pair $(X^\circ,X^\natural)\in\mathcal H^2(\mathbb F,\infty;\mathbb R^m)\times\mathcal H^{2,d}(\mathbb F,\infty;\mathbb R^n)$, \emph{i.e.} we do not necessarily use the pair $(X^c,X^d)$ which is naturally associated to the martingale $X$. 
Observe that in this case we do allow the first component to have jumps.
This is particularly useful when one needs to decompose a discrete--time martingale as a sum of an It\=o integral, a stochastic integral with respect to an integer--valued random measure and a martingale orthogonal to the space of stochastic integrals. 
\begin{proposition}\label{prop:OrthogDecomp}
Let $(Y,X^\circ,X^\natural)\in\mathcal H^2(\mathbb F,\infty;\mathbb R)\times\mathcal H^2(\mathbb F,\infty;\mathbb R^m)\times\mathcal H^{2,d}(\mathbb F,\infty;\mathbb R^n)$ with $M_{\mu^{X^\natural}}[\Delta X^\circ|\widetilde{\mathcal P}^{\mathbb F}]=0$, where the equality is understood componentwise.
Then, there exists a pair $(Z,U)\in\mathbb{H}^2({X}^\circ,\mathbb F,\infty;\mathbb R^m)\times\mathbb{H}^2(\mu^{X^\natural}\mathbb F,\infty;\mathbb R)$ and $N\in\mathcal H^2(\mathbb F,\infty;\mathbb R)$ such that
\begin{align}\label{eq:orth-deco2}
Y= Y_0 + Z\cdot {X}^\circ + U\star\mutilde^{X^\natural} + {N},
\end{align}
with $\langle {X}^{\circ},N\rangle = 0$ and $M_{\mu^{X^\natural}}\big[ \Delta {N} |\widetilde{\mathcal{P}}^{\mathbb F} \big]=0$. 
Moreover, this decomposition is unique, up to indistinguishability. 
\end{proposition}

In other words, the orthogonal decomposition of $Y$ with respect to the pair $(X^\circ,X^\natural)$ is well--defined under the above additional assumption on the jump parts of the martingales $X^\circ$ and $X^\natural$.

\vspace{0.5em}
We conclude this subsection with some useful results. 
Let $\oX:=(X^\circ,X^\natural)\in\mathcal H^2(\mathbb R,\infty;\mathbb R^m)\times\mathcal H^{2,d}(\mathbb F,\infty;\mathbb R^n)$ with $M_{\mu^{X^\natural}}[\Delta X^\circ|\widetilde{\mathcal P}^{\mathbb F}]=0$.
Then we define
\begin{align*}
\mathcal{H}^2(\oX^\perp,\mathbb F,\infty;\mathbb R)  :=\big(\mathcal{L}^2(X^\circ,\mathbb F,\infty;\mathbb R)\oplus\mathcal{K}^2(\mu^{X^\natural},\mathbb F,\infty;\mathbb R)\big)^\perp.
\end{align*}  

\begin{proposition}\label{prop:CharacterOrthogSpace}
Let $\oX:=(X^\circ,X^\natural)\in\mathcal H^2(\mathbb F,\infty;\mathbb R^m)\times\mathcal H^{2,d}(\mathbb F,\infty;\mathbb R^n)$ with $M_{\mu^{X^\natural}}[\Delta X^\circ|\widetilde{\mathcal P}^\mathbb F]=0$.
Then,
\begin{align*}
\mathcal{H}^2(\oX^\perp,\mathbb F,\infty;\mathbb R) = \big\{ L\in\mathcal H^2(\mathbb F,\infty;\mathbb R), \langle X^{\circ}, L\rangle^\mathbb F=0 \text{ and } M_{\mu^{X^\natural}}[\Delta L|\widetilde{\mathcal P}^{\mathbb F}]=0\big\}.
\end{align*}
Moreover, the space $\big(\mathcal{H}^2(\oX^\perp,\mathbb F,\infty;\mathbb R), \Vert \cdot\Vert_{\mathcal H^2(\mathbb F,\infty;\mathbb R)}\big)$ is closed.
\end{proposition}

\begin{corollary}\label{cor:SpaceDecomposition}
Let $\oX:=(X^\circ,X^\natural)\in\mathcal H^2(\mathbb F,\infty;\mathbb R^m)\times\mathcal H^{2,d}(\mathbb F,\infty;\mathbb R^n)$ with $M_{\mu^{X^\natural}}[\Delta X^\circ|\widetilde{\mathcal P}^{\mathbb F}]=0$.
Then, 
\begin{align*}
\mathcal H^2(\mathbb R) = \mathcal L^2(X^\circ,\mathbb F,\infty;\mathbb R) \oplus \mathcal K^2(\mu^{X^\natural},\mathbb F,\infty;\mathbb R)\oplus \mathcal H^2(\oX^\perp,\mathbb F,\infty;\mathbb R),
\end{align*}
where each of the spaces appearing in the above identity is closed.
\end{corollary}

\vspace{.5em}
In view of the above results, we are going to strengthen Condition \ref{MSI} to the following one
\begin{enumerate}[label=\textbf{(M2$^{\prime}$)},itemindent=-.cm,leftmargin=!]
\item\label{MSIprime}
There is a pair $(X^{k,\circ}, X^{k,\natural})\in\mathcal H^2(\mathbb G^k,\infty;\mathbb R^\ell)\times\mathcal H^{2,d}(\mathbb G^k,\infty;\mathbb R^\ell)$ with
$M_{\mu^{X^{k,\natural}}}\big[\Delta X^{k,\circ}\big| \widetilde{\mathcal P}^{\mathbb G^k}\big] = 0$
 for every $k\in\overline{\mathbb N}$, such that
in addition $X^k = X^{k,\circ} + X^{k,\natural}$,
and
\begin{gather}
\big(X^{k,\circ}_{\infty}, X^{k,\natural}_{\infty}\big) \LGconvmulti{2}{\ell\times 2} \big(X^{\infty, c}_{\infty},X^{\infty,d}_{\infty}\big).
\label{MSISeparateFinalRV}
\end{gather}
\end{enumerate}

\begin{corollary}\label{StrongRobMartRep}
Let conditions {\rm\ref{MFilqlc}}, {\rm\ref{MSIprime}} and {\rm\ref{MXWPRP}--\ref{Mfinalrv}} hold. 
Then Theorem \ref{RobMartRep} is valid and the convergence \eqref{RobMartRepii} can be improved into the following one
\begin{align}\label{RobStrongMartRepii}
	\big(\langle Y^k,X^{k,\circ}\rangle,\langle Y^k,X^{k,\natural}\rangle\big)
	\xrightarrow{\hspace{0.2em}(\textup{J}_1(\mathbb R^\ell\times\mathbb R^\ell), \mathbb L^1)\hspace{0.2em}}
	\big(\langle Y^{\infty,c},X^{\infty,c}\rangle,\langle Y^{\infty,d},X^{\infty,d}\rangle\big),
\end{align} 
where we have defined $\langle Y^k, X^{k,\circ}\rangle^i:=\langle Y^k, X^{k,\circ,i}\rangle,$ for $i=1,\dots,\ell$, and $k\in\mathbb N$, and analogously for the processes $\langle Y^k, X^{k,\natural}\rangle$, $k\in\mathbb N$, $\langle Y^\infty, X^{\infty,c}\rangle$ and $\langle Y^\infty, X^{\infty,d}\rangle$.
\end{corollary}

\begin{proof} 
Theorem \ref{RobMartRep} obviously applies in the current framework, since the sequence $(X^{k,\circ})_{k\in\mathbb N}$ approximates $X^{\infty,c}$, \emph{i.e.} the associated sequence of jump processes will finally vanish.
Due to the bilinearity of the dual predictable projection, we obtain the following convergence
\begin{align}\label{eq:311-help}
\big\langle Y^k, X^{k,\circ}\big\rangle+\big\langle Y^k,X^{k,\natural}\big\rangle=\big\langle Y^k,X^k\big\rangle 	
\silmulti{\ell}{1} 	
\big\langle Y^\infty,X^\infty\big\rangle=\big\langle Y^\infty,X^{\infty,c}\big\rangle+\big\langle Y^\infty,X^{\infty,d}\big\rangle.
\end{align}
By convergence  \eqref{MSISeparateFinalRV} and \ref{Mfinalrv},  we obtain the following convergence 
\begin{align}
\big(Y^k_{\infty}+X^{k,\circ,i}_{\infty}, Y^k_{\infty}-X^{k,\circ,i}_{\infty}\big) 
\LGconvmulti{2}{2} 
\big(Y^\infty_{\infty}+X^{\infty,c,i}_{\infty},Y^\infty_{\infty}-X^{\infty,c,i}_{\infty}\big),
\end{align}
for every $i=1,\dots,\ell$.
Theorem \ref{MeminCorollary}.(i) yields that the predictable quadratic covariation of the processes above also converge, and then the polarisation identity allows us to deduce that 
$$\langle Y^k, X^{k,\circ} \rangle \silmulti{\ell}{1}  \langle Y^\infty,X^{\infty,c}\big\rangle = \langle Y^{\infty,c},X^{\infty,c}\big\rangle.$$ 
The statement now follows from the continuity of the angle bracket of the limiting processes, due to the quasi--left--continuity of the limiting filtrations, the convergence in \eqref{eq:311-help} and Lemma \ref{JointSkorokhodConv}.
\end{proof}

\iftoggle{full}{\begin{remark}
The last corollary generalizes \citet[Corollary 2.6]{madan2015convergence}.
Indeed, they consider a discrete--time process approximating a L\'evy process, and work with the natural filtrations, while we can deal both with discrete-- and continuous--time approximations of general martingales, with arbitrary filtrations.
\end{remark}}{}

\subsection{Comparison with the literature}\label{sec:comp}

In this section, we will compare the results presented in Subsections \ref{subsec:RoMRFrame} and \ref{subsec:Ex_and_Appl} with analogous results in the existing literature, namely, with \citet[Theorem~5]{briand2002robustness} and \citet[Theorem~3.3]{jacod2000explicit}.
This discussion serves also as an introduction to the next section, where we will try to clarify some technical points of the proof. 
For the convenience of the reader, we have adapted the notation of the aforementioned articles to our notation.

\vspace{0.5em}
We will follow the chronological order for our discussion, \emph{i.e.} we will start with the comparison of \cref{RobMartRep} with \cite[Theorem 3.3]{jacod2000explicit}.
There, the authors consider a single filtration, \emph{i.e.} $\FilG^k=\FilG^\infty$ for every $k\in\overline{\N}$, where $\mathbb G^\infty$ is an arbitrary filtration. 
The reader should observe that under this framework, $\ref{MFilweak}$ reduces to a triviality. 
Additionally, since the filtration is chosen to be arbitrary, conditions \ref{MFilqlc} and \ref{MXWPRP} are not necessarily satisfied.
Regarding the stochastic integrators, $X^{k,\circ}$ is a locally square--integrable real--valued $\mathbb G^\infty-$martingale and $X^{k,\natural}=0$, for every $k\in\overline{\mathbb N}$.
Therefore, the authors deal with Kunita--Watanabe decompositions, \emph{i.e.} each $Y^k$ consists of an It\=o process and its orthogonal martingale, for every $k\in\overline{\mathbb N}$.
Schematically, the convergence indicated by solid arrows in the following scheme holds, and it is proved in \cite[Theorem~3.3]{jacod2000explicit} that the convergence of the respective parts of the Kunita--Watanabe decompositions (indicated with the dashed arrows) also holds
\begin{center}
\begin{tabular}{l|@{\hspace{1em}}lcccc}
$X^{k}$									&$Y^k$ 								&$=$& 	$Z^k\cdot X^{k} $ 						&$+$& 	$N^k$		\\
\hspace{0.15em}$\big\downarrow$ 		&\hspace{-0.05em}$\big\downarrow$ 	& 	& 	\hspace{-0.2em}$\raisebox{2.0ex}{\rotatebox{-90}{$\dashrightarrow$}}$				& 	& \hspace{-0.2em}\raisebox{2.0ex}{\rotatebox{-90}{$\dashrightarrow$}}\\
$X^{\infty}$ 								&$Y^\infty$							&$=$& 	$Z^\infty\cdot X^{\infty}$ 	&$+$& 	$N^\infty$
\end{tabular}.
\end{center}
In other words, the result of \cite{jacod2000explicit} is more general than \cref{RobMartRep}, in the sense that the orthogonal martingale part is non--zero, and more restrictive in the sense that it considers a single filtration and Kunita--Watanabe decompositions.\medskip

Let us rewrite the above scheme by means of an It\=o stochastic integral with respect to the continuous martingale part $X^{k,c}$, and of a stochastic integral with respect to the integer--valued measure $\mu^{X^{k,d}}$, for every $k\in\overline{\mathbb N}$
\footnote{For simplicity, we will also assume in the following discussion that every local martingale is a square--integrable martingale. This will of  course not be the case in the subsequent sections.\label{ComparisonSimplification}},
\emph{i.e.}
\begin{center}
\begin{tabular}{l|@{\hspace{1em}}lcccc}
$X^{k}$									&$Y^k$ 								&$=$& 	$Z^k\cdot X^{k,c} + (Z^k\text{Id})\star\widetilde{\mu}^{(X^{k,d},\mathbb G^\infty)} $ 						&$+$& 	$N^k$		\\
\hspace{0.15em}$\big\downarrow$ 		&\hspace{-0.05em}$\big\downarrow$ 	& 	& 	\hspace{-4.4em}$\raisebox{2.0ex}{\rotatebox{-90}{$\dashrightarrow$}}$				& 	& \hspace{-0.2em}\raisebox{2.0ex}{\rotatebox{-90}{$\dashrightarrow$}}\\
$X^{\infty}$ 								&$Y^\infty$							&$=$& 	$Z^\infty\cdot X^{\infty,c}+ (Z^\infty\text{Id})\star\widetilde{\mu}^{(X^{\infty,d},\mathbb G^\infty)}$ 	&$+$& 	$N^\infty$
\end{tabular}.
\end{center}
Observe that we have written the purely discontinuous part of the martingale $Z^k\cdot X^{k}$ as $(Z^k\text{Id})\star\widetilde{\mu}^{(X^{k,d},\mathbb G^\infty)}$. This follows from \cite[Proposition II.1.30]{jacod2003limit}, in conjunction with the fact that $X^{k,d}=\text{Id}\star\widetilde{\mu}^{(X^{k,d},\mathbb G^k)}$, for every $k\in\overline{\mathbb N}$.
Assume moreover, that we are interested in proving a result analogous to \cite[Theorem 3.3]{jacod2000explicit} for orthogonal decompositions, and not only for Kunita--Watanabe decompositions.
Then, in view of the above scheme, we should not restrict ourselves to integrands of the form $W\text{Id}$ for the stochastic integral with respect to an integer--valued measure, where $W$ is a predictable process which is determined by the It\=o integrand.
In other words, given that the convergence indicated by solid arrows in the following scheme holds, we would like to prove that the convergence indicated by the dashed arrows also holds
\begin{align*}
\begin{tabular}{l|@{\hspace{1em}}lcccc}
$X^{k}$				&$Y^k$ 								&$=$& 	$Z^k\cdot X^{k,c}+U^k\star\widetilde{\mu}^{(X^{k,d},\mathbb G^\infty )} $ 						&$+$& 	$N^k$		\\
\hspace{0.15em}$\big\downarrow$ 		&\hspace{-0.05em}$\big\downarrow$ 	& 	& 	\hspace{-2.6em}$\raisebox{2.0ex}{\rotatebox{-90}{$\dashrightarrow$}}$				& 	& \hspace{-0.5em}\raisebox{2.0ex}{\rotatebox{-90}{$\dashrightarrow$}}\\
$X^{\infty}$ 				&$Y^\infty$							&$=$& 	$Z^\infty\cdot X^{\infty,c}+U^\infty\star\widetilde{\mu}^{(X^{\infty,d},\mathbb G^\infty)}$ 	&$+$& 	$N^\infty$
\end{tabular}.
\end{align*}
However, if we were to try to follow arguments analogous to \cite[Theorem 3.3]{jacod2000explicit} in order to prove the result described by the last scheme, we would not be able to conclude.
In order to intuitively explain why, let us introduce some further notations.
To this end, recall that we are under the framework of \cite[Theorem 3.3]{jacod2000explicit}, \emph{i.e.} $\mathbb G^k=\mathbb G^\infty$ for every $k\in\mathbb N$, and fix a $k\in\mathbb N$.
Moreover, recall by the last scheme that $N^k$ is the martingale obtained by the orthogonal decomposition of $Y^k$ with respect to $X^k$, and assume the following orthogonal decompositions hold
\begin{align*}
X^{k} &= X^k_0 + \gamma^k\cdot X^{\infty,c} + \beta^k\star\widetilde{\mu}^{(X^{\infty,d},\mathbb G^\infty)} + L^{k}, 
	\text{ where }\langle X^{\infty,c},L^{k,c}\rangle=0\text{ and }M_{\mu^{X^{\infty,d}}}[\Delta L^k|\mathcal P^{\mathbb G^\infty}]=0,
\\
N^{k} &=  \lambda^k\cdot X^{\infty,c} + \zeta^k\star\widetilde{\mu}^{(X^{\infty,d},\mathbb G^\infty)} + R^k, 
	\text{ where }\langle X^{\infty,c},R^{k,c}\rangle=0\text{ and }M_{\mu^{X^{\infty,d}}}[\Delta R^k|\mathcal P^{\mathbb G^\infty}]=0.
\end{align*}
Using arguments analogous to \cite[Theorem 3.3]{jacod2000explicit}, we would be required at some point to write the stochastic integrals $U^k\star\widetilde{\mu}^{(X^{k,d},\mathbb G^k)}$ and $\zeta^k\star\widetilde{\mu}^{(X^{\infty,d},\mathbb G^k)}$ as the sum of stochastic integrals with respect to $\widetilde{\mu}^{(X^{\infty,d},\mathbb G^\infty)}$ and $\widetilde{\mu}^{(L^{k,d},\mathbb G^\infty)}$.
This is indeed possible when the predictable function $\beta^k$, respectively $U^k,\zeta^k$, is of the form $W_1\text{Id}$, respectively $W_2\text{Id},W_3\text{Id}$, for some predictable process $W_1$, respectively $W_2,W_3$.
However, nothing guarantees that this is possible in the general case, which means that we cannot follow through with the proof.

\vspace{0.5em}
Our work intends to provide a general limit theorem. 
Therefore, a condition which imposes a specific relationship between the approximating filtrations $\mathbb{G}^k$ and the limiting filtration $\mathbb{G}^\infty$ would seem restrictive, \emph{e.g.}, $\mathbb{G}^k\subset \mathbb{G}^\infty$ for every $k\in\mathbb{N}$.
On the other hand, it is well--known that any right-continuous martingale can be embedded into a Brownian motion, where the filtration  associated to the Brownian motion is in general larger than its natural filtration and also depends on the family of stopping times that are used to embed the right--continuous martingale; see \cite[Theorem 11]{monroe1972}.
This may lead one to think that we can, without loss of generality, embed every element of the convergent sequence of martingales $(X^k)_{k\in\mathbb{N}\cup\{\infty\}}$ into a Brownian motion, though, in this case, the associated filtration $\mathbb{G}^\infty$ would automatically depend on the sequence of families of stopping times used.
In special cases, \emph{e.g.}, in the embedding of a sequence of random walks into a Brownian motion, the stopping times used for the embedding are stopping times with respect to the natural filtration of said Brownian motion. Therefore, for Donsker type approximations, it seems that one may use this assumption without loss of generality.  
However, for the general case, it is not clear how one can verify the weak predictable representation property of $X^\infty$ with respect to $\mathbb{G}^\infty$. Indeed, observe that if we consider only the natural filtrations then there is no reason why the assumption 
$\mathbb{G}^k\subset \mathbb{G}^\infty$ should hold and in the best case we fall back in the framework presented in our work. 
As one can see, the two questions posed in this comment are closely related. 
Returning to our initial formulation and in view of the Jacod--Yor theorem, see \cite[Theorem 18.3.6]{cohen2015stochastic} in conjunction with \cite[Theorem 14.5.7]{cohen2015stochastic}, the case of a sequence $(X^k)_{k\in\mathbb{N}\cup\{\infty\}}$ of additive martingales with respect to their natural filtrations provides a concrete example.

\vspace{0.5em}
We proceed now to the comparison of \cref{StrongRobMartRep} and \cite[Theorem 5]{briand2002robustness}. 
The authors of \cite{briand2002robustness} consider the Brownian motion case on a finite and deterministic time interval $[0,T]$.
Therefore, in order to translate the finite time horizon into the positive real half--line, we will assume in the following that the processes are indexed on $\mathbb R_+$ but are constant on the time interval $[T,\infty)$. 
More precisely and by means of the notation we have introduced, $X^{\infty}$ is a real--valued $\mathbb F^{X^\infty}-$Brownian motion which is approximated by the sequence $(X^{k,\circ})_{k\in\mathbb N}$ under the locally uniform topology in $\mathbb L^2-$mean, where $X^{k,\circ}$ is a square--integrable $\mathbb F^{X^k}-$martingale for every $k\in\mathbb N$. 
Obviously, $X^{k,\natural}=0$ for every $k\in\overline{\mathbb N}$.
The reader may recall that the convergence $X^{k,\circ}\xrightarrow{\hspace{0.3em} (\text{lu},\mathbb L^2)\hspace{0.3em}} X^\infty$
is equivalent to $X^{k,\circ}\xrightarrow{\hspace{0.3em} (\text{J}_1(\mathbb R),\mathbb L^2)\hspace{0.3em}} X^\infty$, due to the continuity of the limit $X^\infty$.
In view of the aforementioned convergence and due to the special role of the time $T$, we have that 
$$X_\infty^{k,\circ}=X^{k,\circ}_T \xrightarrow{\hspace{0.3em} \mathbb L^2(\Omega, \mathcal G,\mathbb P;\mathbb R)\hspace{0.3em}} X^\infty_T=X^\infty_\infty.$$
Until this point, we can verify that conditions \ref{MFilqlc}, \ref{MSIprime} and \ref{MXWPRP} are satisfied. 
Condition \ref{MFilweak} is guaranteed by \cite[Proposition~3]{briand2002robustness}, while condition \ref{Mfinalrv} is \cite[Assumption (H2)]{briand2002robustness}.
It is immediate now that the framework of \cite[Theorem 5]{briand2002robustness} is stronger than the one we have described for \cref{StrongRobMartRep}.
Moreover, Briand, Delyon and M\'emin consider the Kunita--Watanabe decomposition of the martingale $Y^k$ with respect to $X^{k,\circ}$, for every $k\in\overline{\mathbb N}$, \emph{i.e.} the stochastic integral is an It\=o integral.
Schematically, the convergences indicated by solid arrows in the following scheme hold, and it is proved in \cite[Theorem~5]{briand2002robustness} that the convergence of the respective parts of the Kunita--Watanabe decompositions (indicated with dashed arrows) also holds
\begin{center}
\begin{tabular}{l|@{\hspace{1em}}lcccc}
$X^{k,\circ}$									&$Y^k$ 								&$=$& 	$Z^k\cdot X^{k,\circ} $ 						&$+$& 	$N^k$		\\
\hspace{0.15em}$\big\downarrow$ 		&\hspace{-0.05em}$\big\downarrow$ 	& 	& 	\hspace{-0.2em}$\raisebox{2.0ex}{\rotatebox{-90}{$\dashrightarrow$}}$				& 	& \hspace{-0.2em}\raisebox{2.0ex}{\rotatebox{-90}{$\dashrightarrow$}}\\
$X^\infty$ 								&$Y^\infty$							&$=$& 	$Z^\infty\cdot X^{\infty}$ 	&$+$& 	$0$
\end{tabular}.
\end{center}
Therefore, we can identify \cite[Theorem 5]{briand2002robustness} as a special case of \cref{StrongRobMartRep}. 

\vspace{0.5em}
Let us also briefly describe the technique for the proof that the authors have followed in \cite{briand2002robustness}.
We will do so, because we are going to follow the same technique, \textit{mutatis mutandis}, in order to prove Theorem \ref{RobMartRep}.
A sufficient condition to conclude the required result is to prove that an arbitrary weak--limit point of the sequence $(N^k)_{k\in\N}$, say $\overline{N}$, is orthogonal to $X^\infty$, 
\emph{i.e.} $[X^\infty,\overline{N}]$ is an $\Fil^{(X^\infty,\overline{N})^\top}-$martingale or equivalently 
\begin{equation}\label{eq:suffBM}
\langle X^{\infty},\overline{N}\rangle^{\Fil^{(X^\infty,\overline N)^\top}}=0.
\end{equation}
The reader should keep in mind that the predictable quadratic variation of $[X^\infty,\overline{N}]$ is determined with respect to the filtration $\Fil^{(X^\infty,\overline{N})^\top}$ and not merely the filtration $\Fil^{X^\infty}$. 
The reason is that the $\mathbb F^{X^\infty}-$measurability of $\overline{N}$ cannot be \emph{a priori} guaranteed.
Let us now briefly argue why \eqref{eq:suffBM} is a sufficient condition.
We start by observing that by definition we have 
$$\Fil^{X^\infty}\subset\Fil^{(X^\infty,\overline N)^\top},\text{ and therefore }\mathcal{P}^{\Fil^{X^\infty}}\subset\mathcal{P}^{\Fil^{(X^\infty,\overline N)^\top}}.$$
Hence, using the well--known property (see \cite[Theorem I.4.40.d)]{jacod2003limit})
\begin{equation}\label{eq:suffBMb}
H\cdot\langle X^{\infty},\overline{N}\rangle^{\Fil^{(X^\infty,\overline N)^\top}}=\langle H\cdot X^{\infty}, \overline{N}\rangle^{\Fil^{(X^\infty,\overline N)^\top}},
\text{ for $H\in\mathbb{H}^2\big(X^{\infty},\Fil^{(X^{\infty},\overline N)^\top},T;\R{}\big)$},
\end{equation} 
the authors can conclude in particular that $\overline{N}\mutorthog L,$ for every $L\in\mathcal L^2(X^{\infty},\Fil^{X^{\infty}},T;\R{})$.
For the equivalent forms of orthogonality see \cite[Proposition I.4.15]{jacod2003limit}.
Assume now that $(N^{k_l})_{l\in\mathbb{N}}$ is the subsequence that approximates $\overline N$, \emph{i.e.} $N^{k_l}\xrightarrow[k\to\infty]{\hspace{0.1cm}\mathcal{L}\hspace{0.1cm}}\overline N.$ 
Using condition \eqref{eq:suffBMb}, the fact that $Y^\infty$ can be written in the form $Y^\infty=Z^\infty\cdot X^\infty$ for some $Z\in\mathbb H^2(X^\infty,\mathbb F^{X^\infty}, T;\mathbb R)$, and recalling that  $X^{\infty}$ possesses the $\Fil^{X^{\infty}}-$predictable representation property,
Briand, Delyon and M\'emin can conclude that $\langle N^{k_l}\rangle\xrightarrow{\hspace{0.1cm}(\textrm{lu},\mathbb{L}^1)\hspace{0.1cm}} 0$ and consequently $\overline N=0.$
Since the last convergence holds for every weak--limit point of the sequence $(N^k)_{k\in\mathbb N}$ and due to the sequential compactness of $(N^k)_{k\in\N}$ (recall that $\big(C([0,T];\mathbb R),d_{\rm lu}\big)$ is Polish), the authors can conclude that $N^k\xrightarrow{\hspace{0.1cm}(\textrm{lu},\mathbb{L}^2)\hspace{0.1cm}} 0.$
The convergence $Z^k\cdot X^{k,\circ} \xrightarrow{\hspace{0.1cm}(\textrm{lu},\mathbb{L}^2)\hspace{0.1cm}} Z^\infty\cdot X^\infty$ follows then automatically.

\vspace{0.5em}
After this discussion, the reader may have already wondered what difficulties arise if we do not impose \ref{MXWPRP}. 
In order to explain briefly the issue, let us omit \ref{MXWPRP} from the set of assumptions, \emph{i.e.} $N^\infty\neq 0$ in the Convergence \eqref{RobMartRepi}.
Then, by recalling the arguments above, we have that $(N^k)_{k\in\mathbb N}$ is a tight sequence, therefore we can assume in general that the set $\{\overline{N}: \overline{N} \text{ is a weak limit point of} (N^k)_{k\in\mathbb N}\}$ is not a singleton and we can approximate every element up to a subsequence.
Recall also that the arbitrary weak limit of $(N^k)_{k\in\mathbb N}$ is not necessarily $\mathbb G^\infty-$adapted
\footnote{If such a property could be proved, then using the fact that $\overline N$ is orthogonal to every element of $\mathcal L^2(X^{\infty,c},\mathbb G^\infty,\infty;\mathbb R)\oplus\mathcal K^2(\mu^{X^\infty,d},\mathbb G^\infty,\infty;\mathbb R)$, we could then easily conclude that $\overline{N} =N^\infty$, by the uniqueness of the orthogonal representation.}.
Then, the best we can say about the relationship between $N^\infty$ and the arbitrary weak--limit $\overline{N}$ is that the $\mathbb G^\infty-$optional projection of $\overline{N}$ is indistinguishable from $N^\infty$.
Despite our best efforts, we did not manage to improve this result, leading us to assume that \ref{MXWPRP} holds.   

\medskip
We close this subsection with a short discussion on the set of conditions \ref{MFilqlc}--\ref{Mfinalrv}.
As we have already stated in Section \ref{sec:skor}, Theorems \ref{MeminTheorem} and \ref{MeminCorollary} (which are respectively \cite[Theorem 11, Corollary 12]{memin2003stability} restated in the multidimensional case) constitute the cornerstones for our work.
The former will be used in order to construct convergent sequences of martingales, while the latter will guarantee that the respective sequences of dual predictable quadratic variations will be also convergent.
Therefore, it is natural for our results to be built on the framework of the aforementioned results. 
Indeed, conditions \ref{MFilqlc}, \ref{MSI}, \ref{MFilweak} and \ref{Mfinalrv} are those which guarantee that $(X^k,Y^k,\mathbb G^k)\xrightarrow{\hspace{0.3em}\text{ext}\hspace{0.3em}}(X^\infty,Y^\infty,\mathbb G^\infty)$, 
recall \cref{EquivalentExtended}.
Condition \ref{MXWPRP} will be needed in order to characterize the arbitrary weak--limit of the sequence $(N^k)_{k\in\mathbb N}$; for more details see \cref{subsec:Proof}.  

\subsection{Outline of the proof}\label{sec:outline}

In this subsection, we present the strategy and an overview of the main arguments used to prove Theorem \ref{RobMartRep}, in order to ease the understanding of the technical parts that follow.

\vspace{0.5em}
The first part of the statement amounts to showing the following convergences
\begin{align}\label{eq:outline1}
\begin{matrix}
Y^k 		&=& Y^k_0 		 &+& Z^k\cdot X^{k,c} 			   + U^k\star \mutilde^{(X^{k,d},\mathbb G^k)} 				&+& N^k\\
\hspace{-0.6em}\raisebox{2.0ex}{\rotatebox{-90}{$\dashrightarrow$}}	&& \hspace{-0.4em}\raisebox{2.0ex}{\rotatebox{-90}{$\dashrightarrow$}}   && \hspace{-2.5em}\raisebox{2.0ex}{\rotatebox{-90}{$\dashrightarrow$}}                                     && \hspace{-0.4em}\raisebox{2.0ex}{\rotatebox{-90}{$\dashrightarrow$}}\\  
Y^{\infty} 	&=& Y^{\infty}_0 &+& Z^{\infty}\cdot X^{\infty,c}\ + U^{\infty}\star\mutilde^{(X^{\infty,d},\mathbb G^\infty)} 	&+& 0.
\end{matrix} 
\end{align}
The definition of $Y^k$ as the optional projection of $\xi^k$ with respect to $\mathbb G^k$, for $k\in\overline{\mathbb N}$, together with Assumptions \ref{MFilweak} and \ref{Mfinalrv} and Proposition \ref{MeminProp1}, yield directly that $Y^k\sipmulti{} Y^\infty$.
In particular, $Y^k_0\inp Y^\infty_0$.
Thus, the sum on the right-hand of \eqref{eq:outline1} converges, and we can conclude if we show that $N^k\sipmulti{}0$.

\vspace{0.5em}
\begin{enumerate}[leftmargin=0cm, itemindent=1.02cm, labelsep=0cm, label={\rm \bf Step $\arabic*$}]
\item\label{StepA}:\textit{\hfil$N^k\sipmulti{}0$.}	

\vspace{0.5em}
\noindent We will show that $\langle N^k \rangle \xrightarrow{\hspace{0.2cm}\mathcal L \hspace{0.2cm}} 0$, which is equivalent to $\langle N^k \rangle \sipmulti{} 0$.
The integrability of the sequence $(\langle N^k\rangle)_{k\in\mathbb N}$, in conjunction with Doob's maximal inequality, allow us then to deduce that $N^k\sipmulti{} 0$.

\vspace{0.5em}
\item\label{StepB}:\hfil \textit{$\langle N^k\rangle\xrightarrow{\hspace{0.2cm}\mathcal L \hspace{0.2cm}}0$.}

\vspace{0.5em}
\noindent We cannot show the statement directly, since we do not even know if $(N^k)_{k\in\N}$ has a well--defined limit.
Now, notice that $\langle N^k \rangle = \langle Y^k, N^k \rangle$, for every $k\in\overline{\mathbb N}$, by the orthogonality in the martingale representation.
We will thus show instead that $\langle Y^\infty, \overline N \rangle^\mathbb F = 0$, where $\overline N$ is a weak--limit point of $(N^k)_{k\in\N}$, $\mathbb F$ is a filtration such that  $\mathcal G_t^\infty\subset \mathcal F_t$, for every $t\in[0,\infty)$, and 
$\langle Y^\infty, \overline N \rangle^\mathbb F$ is the dual $\mathbb F-$predictable projection of $[Y^\infty,\overline N]$. 
This is equivalent to proving that $[Y^\infty, \overline N]$ is an $\mathbb F-$martingale, and a sufficient condition for the latter is the following
\begin{align}\label{eq:outline2}
\langle X^{\infty,c}, \overline N^c \rangle^\Fil = 0, \ 
M_{\mu^X}[\Delta \overline N | \widetilde{\mathcal P}^\Fil]= 0, \ 
\text{and}\  \overline N \text{ is an } \Fil-\text{martingale}.
\end{align}
Hence, having showed that $\overline N=0$ for every weak--limit point $\overline N$ of $(N^k)_{k\in\mathbb N}$, we can show \textit{a posteriori} the required convergence.

\vspace{0.5em}
\item\label{StepC}:\hfil\textit{A sufficient condition for \eqref{eq:outline2}.}

\vspace{0.5em}
\noindent This amounts to showing that 
\begin{align}\label{eq:outline3}
[X^\infty, \overline N] \text{ and }
\left[ \int_0^\cdot \int_{\R \ell} h(x) \mathds{1}_I(x) \widetilde\mu^{(X^{\infty,d},\mathbb G^\infty)}(\ds,\dx), \overline N \right]
\text{ are uniformly integrable } \mathbb F-\text{martingales},
\end{align}
for a suitable positive and deterministic function $h$ and for a suitable family of sets $I$.

\vspace{0.5em}
\item\label{StepD}:\hfil \textit{A sufficient condition for \eqref{eq:outline3}.}

\vspace{0.5em}
\noindent This amounts to proving convergence \eqref{conv:Step4} below.  
We have that $N^{k_l} \xrightarrow{\hspace{0.2cm}\mathcal L \hspace{0.2cm}} \overline N$, $X^k \xrightarrow{\hspace{0.2cm}\mathcal L \hspace{0.2cm}} X^\infty$, while both sequences possess the $\PUT$ property; see \cite[Definition VI.6.1]{jacod2003limit}.
Hence we can conclude that 
$[X^{k_l},N^{k_l}]\xrightarrow{\hspace{0.2cm}\mathcal L \hspace{0.2cm}}[X^\infty,\overline N]$, 
and $[X^\infty,\overline N]$ is a uniformly integrable martingale as the limit of uniformly integrable martingales. 
Then, we also need to show that
\begin{align}\label{conv:Step4}
\int_0^\cdot \int_{\R \ell} h(x) \mathds{1}_I(x) \widetilde\mu^{(X^{k,d},\mathbb G^k)}(\ds,\dx)
\silmulti{}{2}
\int_0^\cdot \int_{\R \ell} h(x) \mathds{1}_I(x) \widetilde\mu^{(X^{\infty,d},\mathbb G^\infty)}(\ds,\dx),
\end{align}
again for a suitable deterministic and positive function $h$, and a suitable family of sets $I$, in order to obtain the convergence
\begin{align*}
\Big[\int_0^\cdot \int_{\R \ell} h(x) \mathds{1}_I(x) \widetilde\mu^{(X^{k,d},\mathbb G^k)}(\ds,\dx), N^k \Big]
\xrightarrow{\hspace{0.2cm}\mathcal L \hspace{0.2cm}}
\Big[\int_0^\cdot \int_{\R \ell} h(x) \mathds{1}_I(x) \widetilde\mu^{(X^{\infty,d},\mathbb G^\infty)}(\ds,\dx), \overline N \Big].
\end{align*}
\end{enumerate}


\subsection{\texorpdfstring{Step $4$ is valid for $\mathcal J(X^\infty)$}
							{Step 4 is valid for J(X-infty)}}\label{subsec:StepD}

This subsection is devoted to proving that \eqref{conv:Step4} is true for a family $\mathcal J$ of open subsets of $\mathbb R^\ell$.
Throughout this section, we will consider $X^k$ to be a $\mathbb G^k-$martingale, for every $k\in\overline\N$. 
In particular, its jump process is given by \eqref{def:JumpMartin}.
Before we proceed, let us introduce some notations that will be used throughout the rest of Section \ref{sec:RobMartRepSection}.
For a fixed $i\in\{1,\dots,\ell\}$, following the notations used in Appendix \ref{subsec:JointConv}, we introduce the sets
\begin{align*}
W(X^{\infty,i}(\omega))&:= \{u\in\mathbb{R}\setminus\{0\}, \exists t>0 \text{ with }\Delta X^{\infty,d,i}_t(\omega)=u \},\\
V(X^{\infty,i})&:= \{u\in\mathbb{R}\setminus\{0\}, \mathbb P\big(\Delta X^{\infty,d,i}_t=u, \text{ for some } t>0\big)>0 \},\\
\mathcal{I}(X^{\infty,i})&:=\{(v,w)\subset \R{}\setminus \{0\}, vw>0 \text{ and } v,w\notin V(X^{\infty,i})\}.
\end{align*}
By \cite[Lemma VI.3.12]{jacod2003limit}, we have that the set $V(X^{\infty,i})$ is at most countable, for every $i=1,\dots,\ell.$
%
%
%
%
Moreover, we define
$$\mathcal{J}(X^\infty):=\bigg\{\prod_{i=1}^\ell I_i, \text{ where } I_i\in\mathcal I(X^{\infty,i})\cup\big\{ \mathbb R\big\} \text{ for every } i=1,\dots,\ell\bigg\}\setminus\big\{\mathbb R^\ell\big\},$$
and, for every $I:=I_1\times\dots\times I_\ell\in\mathcal J(X^\infty)$, we set $
J_I:=\big\{i\in\{1,\dots,\ell\}, I_i\neq \mathbb R\}\neq\emptyset.$
Let $k\in\overline\N$, $I:=I_1\times\dots\times I_\ell\in\mathcal J(X^\infty)$ and $g:\Omega\times\mathbb R^\ell\longrightarrow \mathbb R$.
Then we define the $\mathbb R^{\ell + 1}-$valued process 
\begin{align*}
\widehat{X}^{k}[g,I]&:=\big((X^k)^\top, 
X^{k,g,I}\big)^\top,
\end{align*}
where
\begin{align*}
X^{k,g,I}_\cdot(\omega)&:= \sum_{0< t\le \cdot}g(\omega,\Delta X^{k,d}_t(\omega)) \mathds{1}_I(\Delta X^{k,d}_t(\omega)) =\int_{(0,\cdot]\times\mathbb R^\ell}g(\omega,x)\mathds{1}_I(x){\mu}^{X^{k,d}}(\omega;\ud s,\ud x).
\end{align*}
Observe that, due to \ref{MSI}, the random variable $X^{k}_\infty$ exists $\mathbb P-a.s.$
Consequently, the process $X^{k,g,I}$, as well as the process $\widehat{X}^{k}[g,I]$, are $\mathbb P-a.s.$ well--defined.
%
%
%
%
%
\begin{proposition}
Let condition {\rm\ref{MSI}} hold.
Fix an $I\in\mathcal J(X^\infty)$ 
and a function $g:(\Omega\times\mathbb R^\ell,\mathcal G\otimes\mathcal B(\mathbb R^\ell))\longrightarrow (\mathbb R,\mathcal B(\mathbb R))$ such that there exists $\Omega^C\subset \Omega$ with $\mathbb P(\Omega^C)=1$, and 
\begin{align*}
&\mathbb R^\ell\ni x\longmapsto g(\omega,x) \in\mathbb R \hspace{0.3em}\text{ is continuous on }\hspace{0.2em}C(\omega) \hspace{0.2em}\text{ for every }\hspace{0.2em}\omega\in\Omega^C,
\shortintertext{where }
C(\omega)&:=\prod_{i=1}^\ell A_i(\omega), \text{ with }A_i(\omega):=
\begin{cases}
W(X^{\infty,i}(\omega)), 		\ \text{if }i\in J_I,\\
W(X^{\infty,i}(\omega))\cup\{0\}, \ \text{if }i\in\{1,\dots,\ell\}\setminus J_I,
\end{cases}
\end{align*}
Then, it holds
$$\widehat{X}^{k}[g,I]
	   	\xrightarrow[k\to\infty]{\hspace*{0.2cm}(\J_1(\mathbb R^{\ell + 1}),\mathbb P)\hspace*{0.2cm}}
 \widehat{X}^{\infty}[g,I].$$
\end{proposition}
%
%
\begin{proof}
Let us fix an $I:=I_1\times \dots \times I_\ell\in\mathcal J(X^\infty)$.
Since the space $(\D(\mathbb R^{\ell + 1}),d_{\J_1(\mathbb R^{\ell + 1})})$ is Polish,
by \citet[Theorem 9.2.1]{dudley2002real}, it is therefore sufficient to prove that for every subsequence 
$(\widehat{X}^{k_l}[g,I])_{l\in\N}$, there exists a further subsequence $(\widehat{X}^{k_{l_m}}[g,I])_{m\in\N}$ for which
\begin{align}\label{eq:SkorContOperatorPas}
\hspace{0.3cm}\widehat{X}^{k_{l_m}}[g,I]\xrightarrow[m\to\infty]{\hspace*{0.2cm}\J_1(\mathbb R^{\ell + 1})\hspace*{0.2cm}}\widehat{X}^{\infty}[g,I],\; \Pm-a.s.
\end{align}
Let $(\widehat{X}^{{k_l}}[g,I])_{l\in\N}$ be fixed hereinafter. 
For every $i\in J_I$, since $I_i\in \mathcal{I}(X^{\infty,i}),$ there exists $\Omega^{I_i}\subset\Omega$ such that 
\begin{enumerate}[label={\rm (\roman*)}]
\item $\Pm\big( \Omega^{I_i}\big)=1$,
\item\label{Prop:hatXgI-ii} 
 $\Delta X^{\infty,d,i}_t(\omega)\notin \partial I_i$, for every $t\in\Rp$ and $\omega \in\Omega^{I_i},$
\end{enumerate}
where $\partial A$ denotes the $|\cdot|-$boundary of the set $A.$

Condition \ref{MSI} implies the convergence $X^k\sipmulti{\ell}X^\infty.$
Hence, for the subsequence $(X^{k_l})_{l\in\N}$, there is a further subsequence $(X^{k_{l_m}})_{m\in\N}$ for which holds
\begin{align}\label{conv:SubPas}
X^{k_{l_m}}\xrightarrow{\hspace*{0.2cm}\J_1(\R{\ell})\hspace*{0.2cm}} X^\infty,\ \Pm-a.s.
\end{align}
Let $\Omega_{\text{sub}}\subset \Omega $ be such that $\Pm(\Omega_{\text{sub}})=1$ and such that the convergence \eqref{conv:SubPas} holds for every $\omega\in\Omega_{\text{sub}}.$ 
Define $\Omega^{J_I}_{\text{sub}}:=\Omega_{\text{sub}}\cap \big(\cap_{i\in J_I}\Omega^{I_i}\big)$. 
Then 
\begin{enumerate}[label={\rm (\roman*$'$)},]
\item  $\Pm(\Omega^{J_I}_{\text{sub}})=1$, 
\item $\Delta X^{\infty,d,i}_t(\omega)\notin \partial I_i$, for every $t\in\Rp$, for every $i\in J_I$ and every $\omega \in\Omega^{J_I}_{\text{sub}}.$
\end{enumerate}
By the last property, we can conclude that $\partial I_i\cap W(X^{\infty,i}(\omega))=\emptyset$, for every $\omega \in \Omega^{J_I}_{\text{sub}}, i\in J_I.$
Therefore, the function $\mathbb R^\ell \ni x\overset{g_{_I}(\omega)}{\longmapsto} g(\omega,x)\mathds{1}_I(x)$ is continuous on the set $C(\omega)$ for $\mathbb P-$almost every $\omega \in \Omega^{J_I}_{\text{sub}}$. We can now conclude once we apply, for each $\omega \in  \Omega^{J_I}_{\text{sub}}$, Proposition \ref{prop:ContFuncInSkor} for the sequence $(X^{k_{l_m}}(\omega))_{m\in\mathbb N}$ and for the function $\mathbb R^\ell\ni x\longmapsto g(\omega,x)\in\mathbb R$. This gives us \eqref{eq:SkorContOperatorPas}.
\end{proof}
%
%
%
%
In order to abridge the notation in the following results, we introduce, for any $p\in\mathbb R_+$, the continuous function
$$\mathbb R^\ell \ni x\overset{{\rm R}_p}{\longmapsto}\sum_{i=1}^\ell (|x^i|\wedge 1)^p\in\mathbb R,$$
where $a\wedge b = \min\{a,b\}$ for $a,b\in\mathbb R$.

\begin{corollary}\label{cor:R2-J1Cont}
Let condition {\rm\ref{MSI}} hold.
Then, for every $I\in \mathcal J(X^\infty)$ it holds
$$\widehat X^k[{\rm R}_p,I]\xrightarrow{\hspace{0.3em}(\textup{J}_1(\mathbb R^{\ell+1}),\mathbb P) \hspace{0.3em}} \widehat X^\infty[{\rm R}_p,I], \ \text{ for every } p\in\mathbb R_+.$$
\end{corollary}
\begin{proof}
Let $p\in\mathbb R_+.$ It suffices to apply the above proposition to the function ${\rm R}_p$, which is continuous.
\end{proof}
%
%
%
%
%
%
%
%
%
Let us now provide more details on the strategy of the proof for this step. 
Proposition \ref{prop:AppMemin} below is the most important result of this subsection, since it provides us with a rich enough family of converging martingale sequences.
It is this family that we are going to use in order to show convergence \eqref{conv:Step4}. 
To this end, we are going to apply Theorem \ref{MeminTheorem} to the sequence $(X^{k,{\rm R}_2,I})_{k\in\overline\N}$, where $I\in\mathcal J(X^{\infty})$. 
However, all of this requires to make sure that this sequence indeed verifies the requirements of the aforementioned theorem, for every $I\in\mathcal J(X^\infty)$. 
This is the subject of the remainder of this subsection. \label{important-comment}



\vspace{.5em}
Let us fix $I=\prod_{i=1}^\ell I_i\in\mathcal J(X^{\infty})$, hence the set $J_I$ is a fixed non--empty subset of $\{1,\dots,\ell\}.$
Moreover, we define $\mathbb R^\ell\ni x \xmapsto{{\rm R}_{2,A}} {\rm R}_2(x) \mathds{1}_A(x),$ for every $A\subset \mathbb R^\ell.$ 


%
%
%
\begin{lemma}\label{lem:XIisSpecial}
The process $X^{k,{\rm R}_2,I}$ is a $\mathbb G^k-$special semimartingale, for every $k\in\overline\N$.
In particular, its $\mathbb G^k$--canonical decomposition is given by
\begin{align*}
X^{k,{\rm R}_2,I}_\cdot
	=&\	\int_0^\cdot\int_{\mathbb R^\ell} {\rm R}_{2,I}(x)\mutilde^{(X^{k,d},\FilG^k)}(\ds,\dx)+ 	\int_{(0,\cdot]\times I} {\rm R}_2(x)\nu^{(X^{k,d},\FilG^k)}(\ds,\dx),
\end{align*}
or equivalently
\begin{align*}
{\rm R}_{2,I}*{\mu}^{X^{k,d}}	=&\ {\rm R}_{2,I}\star \widetilde{\mu}^{(X^{k,d},\mathbb G^k)} +  {\rm R}_{2,I}*{\nu}^{(X^{k,d},\mathbb G^k)}. 
\end{align*}
\end{lemma}
\begin{proof}
Let $k\in\overline\N$.
Observe that by construction the process $X^{k,{\rm R}_2,I}$ is $\mathbb G^k-$adapted and \cadlag.
The function ${\rm R}_{2,I}$ is positive, hence the process $X^{k,{\rm R}_2,I}$ is a $\mathbb G^k-$submartingale of finite variation, as its paths are $\mathbb P-a.s.$ non--decreasing.
Before we proceed, we need to show the integrability of 
${\rm R}_{2,I}* \mu^{X^{k,d}}_{\infty}$ in order to make use of \cite[Proposition II.1.28]{jacod2003limit}. 
But, we have
\begin{align*}
\mathbb E\bigg[ \int_{(0,\infty)\times \mathbb R^\ell}{\rm R}_{2,I}(x)\, \mu^{X^{k,d}}(\ds,\dx) \bigg]
	&\le 	\mathbb E\bigg[ \int_{(0,\infty)\times \mathbb R^\ell}  \sum_{i=1}^\ell |x^i|^2 \mu^{X^{k,d}}(\ds,\dx) \bigg] \\
	&\le 	4 \mathbb E\bigg[ \int_{(0,\infty)\times \mathbb R^\ell}  |x|^2 \mu^{X^k}(\ds,\dx) \bigg]
	\overset{\eqref{SqIntofIVM}}{<}\infty,
	\numberthis \label{cond:IntegMuXk}
\end{align*}
which yields also that $X^{k,{\rm R}_2,I}$ is special, by \cite[Proposition~I.4.23.(iv)]{jacod2003limit}.
Moreover, by \cite[Theorem II.1.8]{jacod2003limit} and condition \eqref{SqIntofIVM} we obtain
\begin{align}\label{cond:IntegNuXk}
\mathbb E\bigg[ \int_{(0,\infty)\times \mathbb R^\ell}{\rm R}_{2,I}(x)\, \nu^{(X^{k,d},\mathbb G^k)}(\ds,\dx) \bigg]<\infty,\text{ for every }k\in\overline{\mathbb N}.
\end{align}
Therefore, we have
\begin{align*}
X^{k,{\rm R}_2,I}_\cdot 
		&= 	\sum_{0< s\le \cdot} {\rm R}_2(\Delta X^{k,d}_s)\mathds{1}_{I}(\Delta X^{k,d}_s) 
		 = \int_{(0,\cdot]\times\mathbb R^\ell} {\rm R}_2(x)\mathds{1}_{I}(x) \mu^{X^{k,d}}(\ud s,\ud x)\\
		&\hspace{-2.7em}
		\overset{\text{\cite[Theorem II.1.28]{jacod2003limit}}}{\underset{\eqref{cond:IntegMuXk}}{=}} 
			\int_0^\cdot\int_{\mathbb R^\ell} {\rm R}_2(x) \mathds{1}_{I}(x)\widetilde{\mu}^{(X^{k,d},\mathbb G^k)}(\ds,\dx) + \int_{(0,\cdot]\times\mathbb R^\ell} {\rm R}_2(x)\mathds{1}_{I}(x) \nu^{(X^{k,d},\mathbb G^k)}(\ud s,\ud x)\\
		&= {\rm R}_{2,I}\star \widetilde{\mu}^{(X^{k,d},\mathbb G^k)} +  {\rm R}_{2,I}*{\nu}^{(X^{k,d},\mathbb G^k)}.
\end{align*}
The finite variation part ${\rm R}_{2,I}*{\nu}^{(X^{k,d},\mathbb G^k)}$ is predictable, since ${\rm R}_{2,I}$ is deterministic and the random measure $\nu^{(X^{k,d},\mathbb G^k)}$ is predictable (see \cite[Definition II.1.6, Theorem II.1.8]{jacod2003limit}).
Hence, we can conclude also via this route that $X^{k,{\rm R}_2,I}$ is a special semimartingale, since it admits a representation as the sum of a martingale and a predictable part of finite variation.
\end{proof}

%
%
%
%
%
\begin{lemma}\label{lem:Optional-UI}
\begin{enumerate}[label={\rm (\roman*)},itemindent=0cm, leftmargin=1cm]
\item \label{lem:Optional-UI-i} The sequence $\big({\rm Tr}\big[[X^{k}]_\infty\big]\big)_{k\in\overline{\mathbb N}}$ is uniformly integrable.
\item \label{lem:Optional-UI-ii} The sequence $\big(X^{k,{\rm R}_2,I}_\infty\big)_{k\in\overline{\mathbb N}}$ is uniformly integrable.
\item \label{lem:Optional-UI-iii} The sequence $\big([X^{k,{\rm R}_2,I}]_{\infty} \big)_{k\in\overline \N}$ is uniformly integrable.
\end{enumerate}
\end{lemma}
\begin{proof}
\ref{lem:Optional-UI-i} In view of conditions \ref{MFilqlc}, \ref{MSI} and \ref{MFilweak}, we have from Theorem \ref{MeminCorollary}.\ref{BracketConvinL} that the sequence $([X^{k,i}]_{\infty})_{k\in\overline\N}$ is uniformly integrable.
By \cite[Corollary 1.10]{he1992semimartingale} we can conclude.

\vspace{0.5em}
\ref{lem:Optional-UI-ii}
Using the definitions of ${\rm R}_2$ and $X^{k,{\rm R}_{2},I}$, we get
\begin{align*}
X^{k,{\rm R}_{2},I}_\infty 
	& = \int_{(0,\infty)\times\mathbb R^\ell}{\rm R}_{2,I}(x)\,\mu^{X^{k,d}}(\ud s,\ud x)
	  \le \int_{(0,\infty)\times \mathbb R^\ell} |x|^2 \, \mu^{X^{k,d}}(\ud s,\ud x) 
	  = {\rm Tr}\big[ [X^{k,d}]_\infty \big]
	  \le {\rm Tr}\big[ [X^{k}]_\infty \big].
\end{align*}
Hence, from \ref{lem:Optional-UI-i} and \cite[Theorem 1.7]{he1992semimartingale} we can conclude.

\vspace{0.5em}
\ref{lem:Optional-UI-iii}
By Lemma \ref{lem:XIisSpecial}, the process $X^{k,{\rm R}_2,I}$ is a $\mathbb G^k-$special semimartingale for every $k\in\overline\N$, whose martingale part is purely discontinuous.
Therefore, we have by \cite[Theorem I.4.52]{jacod2003limit} that
\begin{align*}
[X^{k,{\rm R}_2,I} ]_{\infty} 
	&= 	\sum_{s>0}\big|{\rm R}_2\big(\Delta X^{k,d}_s\big)\big|^2\mathds{1}_{I}(\Delta X^{k,d}_s)
	\le \sum_{s>0}\big|{\rm R}_2\big(\Delta X^{k,d}_s\big)\big|^2
	= 	\sum_{s>0}\bigg(\sum_{i=1}^\ell\big(\big|\Delta X^{k,d,i}_s\big|\wedge 1\big)^2\bigg)^2\\
	&=	\sum_{s>0}\bigg[\sum_{i=1}^\ell\Big(\big|\Delta X^{k,d,i}_s\big|^2\mathds{1}_{(0,1)}\big(\big|\Delta X^{k,d,i}_s\big|\big) + \mathds{1}_{[1,\infty)}\big(\big|\Delta X^{k,d,i}_s\big|\big)\Big)\bigg]^2\\
	&\le 2\ell	\sum_{s>0}\bigg[\sum_{i=1}^\ell\Big(\big|\Delta X^{k,d,i}_s\big|^4\mathds{1}_{(0,1)}\big(\big|\Delta X^{k,d,i}_s\big|\big) + \mathds{1}_{[1,\infty)}\big(\big|\Delta X^{k,d,i}_s\big|\big)\Big)\bigg]\\	
	&\le 2\ell \sum_{s>0}\bigg[\sum_{i=1}^\ell\big(\Delta X^{k,d,i}_s\big)^2\bigg] 
		= 2\ell {\rm Tr}\big[[X^{k,d}]_{\infty}\big]
		\leq 2\ell {\rm Tr}\big[[X^{k}]_\infty\big].
\end{align*}
Thus, using \ref{lem:Optional-UI-i} and \cite[Theorem 1.7]{he1992semimartingale} again, we have the required result.
\end{proof}
%
%
%
%
\begin{lemma}\label{lem:FV}
\begin{enumerate}[label={\rm (\roman*)}, itemindent=0cm, leftmargin=1cm]
\item\label{lem:FV-i} The sequence $\big({\rm Var}({\rm R}_{2,I}*\nu^{(X^{k,d},\mathbb G^k)})_\infty\big)_{k\in\overline{\mathbb N}}$ is tight in $(\mathbb R,|\cdot|)$.
\item\label{lem:FV-ii} The sequence $\big(\sum_{s>0}\big(\Delta ({\rm R}_{2,I}*\nu^{(X^{k,d},\mathbb G^k)})_s\big)^2 \big)_{k\in\overline{\mathbb N}}$ is uniformly integrable.
\end{enumerate}
\end{lemma}
\begin{proof}
\ref{lem:FV-i}
We have already observed that $X^{k,{\rm R}_2,I}$ is a $\mathbb G^k-$submartingale for every $k\in\overline{\mathbb N}$ and consequently ${\rm R}_{2,I}*\nu^{(X^{k,d},\mathbb G^k)}$ is non--decreasing for every $k\in\overline{\mathbb N}$; a property which is also immediate since ${\rm R}_{2,I}$ is a positive function and $\nu^{(X^{k,d},\mathbb G^k)}$ is a (positive) measure. 
Therefore, it holds
\begin{align*}
{\rm Var} \big({\rm R}_{2,I}*\nu^{(X^{k,d},\mathbb G^k)}\big) = {\rm R}_{2,I}*\nu^{(X^{k,d},\mathbb G^k)}, \text{ for every }k\in\overline{\mathbb N}.
\end{align*}
In view of the above and due to Markov's inequality, it suffices to prove that $\sup_{k\in\overline{\mathbb N}}\mathbb E [ {\rm R}_{2,I}*\nu^{(X^{k,d},\mathbb G^k)}_\infty ]<\infty$.
Indeed, for every $\varepsilon >0$ it holds for $K:=\frac{1}{\varepsilon}\sup_{k\in\overline{\mathbb N}}\mathbb E [ {\rm R}_{2,I}*\nu^{(X^{k,d},\mathbb G^k)}_\infty ]>0$ that
\begin{align*}
\sup_{k\in\overline{\mathbb N}}\mathbb P\big[ {\rm R}_{2,I}*\nu^{(X^{k,d},\mathbb G^k)}_{\infty}>K\big]\le \frac{1}{K}\sup_{k\in\overline{\mathbb N}}\mathbb E [ {\rm R}_{2,I}*\nu^{(X^{k,d},\mathbb G^k)}_\infty ]<\varepsilon,
\end{align*}
which yields the required tightness.
Now, observe that we have
\begin{align*}
\mathbb E[{\rm R}_{2,I}*\nu^{(X^{k,d},\mathbb G^k)}_\infty] 
		&\overset{\text{\cite[Proposition II.1.28]{jacod2003limit}}}{=} 
				\mathbb E[ {\rm R}_{2,I}* \mu^{X^{k,d}}_\infty]\overset{\eqref{cond:IntegMuXk}}{<}\infty.
		\numberthis \label{con:Fkinf-integr}
\end{align*}
We have concluded using inequality \eqref{cond:IntegMuXk}, which in turn makes use of Assumption \eqref{SqIntofIVM}. 
Therefore \eqref{con:Fkinf-integr} yields that $\sup_{k\in\overline\N}\mathbb E[{\rm R}_{2,I}*\nu^{(X^{k,d},\mathbb G^k)}_\infty]<\infty.$

\vspace{0.5em}
\ref{lem:FV-ii}
We have, for every $k\in\overline{\mathbb N}$, that the following holds:
\begin{align*}
&\hspace{-1em}\sum_{s>0}\big(\Delta ({\rm R}_{2,I}*\nu^{(X^{k,d},\mathbb G^k)})_s\big)^2 
 	  = 	\sum_{s>0}\Big(\int_{\mathbb R^\ell}{\rm R}_{2,I}(x)\,\nu^{(X^{k,d},\mathbb G^k)}(\{s\}\times\ud x)\Big)^2\\
	& \hspace{-1.4em}\overset{\text{Jensen Ineq.}}{\le} 	
			\sum_{s>0}\int_{\mathbb R^\ell}{\rm R}_{2,I}^2(x)\,\nu^{(X^{k,d},\mathbb G^k)}(\{s\}\times\ud x)
	   \le 	\int_{(0,\infty)\times\mathbb R^\ell}{\rm R}_{2,I}^2(x)\,\nu^{(X^{k,d},\mathbb G^k)}(\ud s,\ud x)\\
	&  = 	\int_{(0,\infty)\times I} \Big[\sum_{i=1}^\ell (|x^i|\wedge1)^2\Big]^2 \nu^{(X^{k,d},\mathbb G^k)}(\ud s,\ud x)
	   \le 	2\ell\int_{(0,\infty)\times I} \sum_{i=1}^\ell (|x^i|^2\wedge1)\, \nu^{(X^{k,d},\mathbb G^k)}(\ud s,\ud x)\\
	&  = 2\ell \, {\rm R}_{2,I}*\nu^{(X^{k,d},\mathbb G^k)}_\infty.  
	\numberthis\label{ineq:SquareFVJumps-UI}
\end{align*}

Using Lemma \ref{lem:Optional-UI}.\ref{lem:Optional-UI-ii} and the first lemma in \citet[p.770]{meyer1978sur}, there exists a moderate\footnote{Cf. Definition \ref{def:Young-constants}.}, Young function $\Phi$ such that 
\begin{align*}
\sup_{k\in\overline{\mathbb N}}\mathbb E\big[ \Phi\big(X^{k,{\rm R}_{2},I}_\infty\big)\big]<\infty.
\end{align*}
Then, using that $X^{k,{\rm R}_{2},I}$ is an increasing process, and is thus equal to its supremum process, the decomposition of Lemma \ref{lem:XIisSpecial} and applying \citet[Th\'eor\`eme 3.2.1]{lenglart1980presentation}, we can conclude that
\begin{align*}
\sup_{k\in\overline{\mathbb N}}\mathbb E\big[ \Phi\big({\rm R}_{2,I}*\nu^{(X^{k,d},\mathbb G^k)}_\infty\big)\big]<\infty.
\end{align*}
By de La Vall\'ee--Poussin's criterion, the latter condition is equivalent to the uniform integrability of the sequence $\big({\rm R}_{2,I}*\nu^{(X^{k,d},\mathbb G^k)}_\infty\big)_{k\in\overline{\mathbb N}}$.
Then, by \eqref{ineq:SquareFVJumps-UI} and \cite[Theorem 1.7]{he1992semimartingale} we can conclude the uniform integrability of the required sequence.
\end{proof}

\begin{corollary}\label{cor:R2I-Integr}
\begin{enumerate}[label={\rm (\roman*)}, itemindent=0cm, leftmargin=1cm]
\item\label{cor:R2I-Integr-ii} The sequence $\big(\big[ {\rm R}_{2,I}\star \widetilde{\mu}^{(X^{k,d},\mathbb G^k)}\big]_\infty\big)_{k\in\overline{\mathbb N}}$ is uniformly integrable.
\item\label{cor:R2I-Integr-i} Consequently, it holds that ${\rm R}_{2,I}\star \widetilde{\mu}^{(X^{k,d},\mathbb G^k)}\in\mathcal H^{2,d}(\mathbb G^k,\infty;\mathbb R)$, for every $k\in\overline{\mathbb N}$.
\end{enumerate}

\end{corollary}
\begin{proof}
Using that ${\rm R}_{2,I}\star \widetilde{\mu}^{(X^{k,d},\mathbb G^k)}$ is a martingale of finite varation, we have
\begin{align*}
\big[ {\rm R}_{2,I}\star \widetilde{\mu}^{(X^{k,d},\mathbb G^k)}\big]_{\infty}
	& = \sum_{s>0} \Big(\Delta({\rm R}_{2,I}\star \widetilde{\mu}^{(X^{k,d},\mathbb G^k)})_s \Big)^2 
	  = \sum_{s>0} \Big(\Delta(X^{k,{\rm R}_2,I})_s - \Delta({\rm R}_{2,I}*\nu^{(X^{k,d},\mathbb G^k)})_s\Big)^2\\
	& \le 	2\sum_{s>0} \Big(\Delta(X^{k,{\rm R}_2,I})_s\Big)^2 + 2\sum_{s>0}\Big(\Delta({\rm R}_{2,I}*\nu^{(X^{k,d},\mathbb G^k)})_s\Big)^2\\
	& = 	2 \, [ X^{k,{\rm R}_{2,I}}]_\infty + 2 \sum_{s>0}\Big(\Delta({\rm R}_{2,I}*\nu^{(X^{k,d},\mathbb G^k)})_s\Big)^2,
\end{align*}
where in the last equality we have used that $X^{k,{\rm R}_{2,I}}$ is a semimartingale whose paths have finite variation and \cite[Theorem I.4.52]{jacod2003limit}.
In view now of the above inequality, Lemma \ref{lem:Optional-UI}.\ref{lem:Optional-UI-iii}, Lemma \ref{lem:FV}.\ref{lem:FV-ii} and \cite[Theorem 1.7, Corollary 1.10]{he1992semimartingale}, we can conclude the required property.
This shows \ref{cor:R2I-Integr-ii}.

\vspace{0.5em}
In addition, \ref{cor:R2I-Integr-ii} implies the integrability of $\big[ {\rm R}_{2,I}\star \widetilde{\mu}^{(X^{k,d},\mathbb G^k)}\big]$, hence from \cite[Proposition~I.4.50.c)]{jacod2003limit} we get that ${\rm R}_{2,I}\star \widetilde{\mu}^{(X^{k,d},\mathbb G^k)}\in\mathcal H^{2,d}(\mathbb G^k,\infty;\mathbb R)$.
\end{proof}
%
%
%
%
%
%
\begin{proposition}\label{prop:AppMemin}
Let conditions {\rm\ref{MFilqlc}}, {\rm\ref{MSI}} and {\rm\ref{MFilweak}} hold. 
Then the following convergence holds
\begin{align}\label{conv:AppMemin}
\big((X^k)^\top, {\rm R}_{2,I}\star\widetilde{\mu}^{(X^{k,d},\mathbb G^k)}\big)^\top  
\xrightarrow{\hspace{0.2cm}(\J_1(\mathbb R^{\ell+1}),\mathbb{L}^2)\hspace{0.2cm}}
\big((X^\infty)^\top,{\rm R}_{2,I}\star\widetilde{\mu}^{(X^{\infty,d},\mathbb G^\infty)}\big)^\top.  
\end{align}
\end{proposition}
%
%
%
%
%
%
\begin{proof}
As we have already pointed out on page \pageref{important-comment}, we are going to apply Theorem \ref{MeminTheorem} to the sequence $(X^{k,{\rm R}_2,I})_{k\in\overline\N}.$ 
By Lemma \ref{lem:XIisSpecial}, this is a sequence of $\mathbb G^k-$special semimartingales, for every $k\in\overline \N$. 

\vspace{0.5em}
In view of \ref{MFilqlc}, which states that $X^\infty$ is $\mathbb G^\infty-$quasi--left--continuous, and using \cite[Corollary II.1.19]{jacod2003limit}, we get that the compensator $\nu^{(X^{\infty,d},\mathbb G^\infty)}$ associated to $\mu^{X^{\infty,d}}$ is an atomless random measure.
Therefore, the finite variation part of the $\mathbb G^\infty-$canonical decomposition of $X^{\infty,{\rm R}_2,I}$ is a continuous process.
Moreover, by \cite[Theorem 5.36]{he1992semimartingale} and \ref{MFilqlc}, which states that the filtration $\mathbb G^\infty$ is quasi--left--continuous, it suffices to show that the martingale part of $X^{\infty,{\rm R}_2,I}$ is uniformly integrable. 
The latter holds by \cref{cor:R2I-Integr}.\ref{cor:R2I-Integr-i}.

\vspace{0.5em}
Lemma \ref{lem:Optional-UI}.\ref{lem:Optional-UI-iii} yields that condition \ref{MeminTh-i} of Theorem \ref{MeminTheorem} holds.
\cref{lem:FV}.\ref{lem:FV-ii} yields that condition \ref{MeminTh-ii} of the aforementioned theorem also holds.
Moreover, from \cref{cor:R2-J1Cont} with $p=2$, we obtain the convergence
\begin{align}\label{eq:JointXXI}
\big( (X^k)^\top, X^{k,{\rm R}_2,I} \big)^\top  
	   	\xrightarrow[]{\hspace*{0.2cm}(\J_1(\mathbb R^{\ell+1}),\mathbb P)\hspace*{0.2cm}}
\big( (X^\infty)^\top, X^{\infty,{\rm R}_2,I} \big)^\top  .
\end{align}
The last convergence in conjunction with conditions \ref{MSI} and \ref{MFilweak}, Remark \ref{EquivalentExtended} and Corollary \ref{cormultiskorokhod}, is equivalent to the convergence 
$(X^{k,{\rm R}_2,I},\mathbb G^k)\extsense (X^{\infty,{\rm R}_2,I},\mathbb G^\infty).$ 
Therefore, condition \ref{MeminTh-iii} of Theorem \ref{MeminTheorem} is also satisfied.

\vspace{.5em} 
Applying now Theorem \ref{MeminTheorem} to the sequence $(X^{k,{\rm R}_2,I})_{k\in\overline\N}$, and keeping in mind the decomposition from \cref{lem:XIisSpecial}, we obtain the convergence
\begin{align}\label{eq:JointXXII}
\big(X^{k,{\rm R}_2,I}, {\rm R}_{2,I}\star\widetilde{\mu}^{(X^{k,d},\mathbb G^k)}\big)^\top
	\xrightarrow[]{\hspace*{0.2cm}(\J_1(\mathbb R^2),\mathbb P)\hspace*{0.2cm}}
\big(X^{\infty,{\rm R}_2,I},{\rm R}_{2,I}\star\widetilde{\mu}^{(X^{\infty,d},\mathbb G^\infty)}\big)^\top.
\end{align}
Using \cref{cormultiskorokhod}, we can combine the convergences in \eqref{eq:JointXXI} and \eqref{eq:JointXXII} to obtain
\begin{align*}
\big((X^{k})^\top, {\rm R}_{2,I}\star\widetilde{\mu}^{(X^{k,d},\mathbb G^k)} \big)^\top
	\xrightarrow[]{\hspace*{0.2cm}(\J_1(\mathbb R^\ell\times \mathbb R),\mathbb P)\hspace*{0.2cm}}
\big((X^{\infty})^\top,{\rm R}_{2,I}\star\widetilde{\mu}^{(X^{\infty,d},\mathbb G^\infty)})^\top.
\end{align*}

\vspace{0.5em}
The last result can be further strengthened to an $\mathbb{L}^2-$convergence in view of the following arguments:
let $\alpha^k: = \big((X^{k})^\top,{\rm R}_{2,I}\star\widetilde{\mu}^{(X^{k,d},\mathbb G^k)}\big)$, then by Vitali's Convergence Theorem, the latter is equivalent to showing that $d^2_{{\rm J}_1}(\alpha^k,\alpha^\infty)$ is uniformly integrable.
Moreover, by the inequality 
$$d^2_{{\rm J}_1}(\alpha^k,\alpha^\infty)
	\le \big(d_{{\rm J}_1}(\alpha^k,0) + d_{{\rm J}_1}(0,\alpha^\infty)\big)^2 
	\le 2d^2_{{\rm J}_1}(\alpha^k,0) + 2d^2_{{\rm J}_1}(0,\alpha^\infty)
	\le 2\Vert\alpha^k\Vert_\infty^2 + 2\Vert \alpha^\infty \Vert_\infty^2, $$
and \cite[Theorem~1.7, Corollary~1.10]{he1992semimartingale}, it suffices to show that $(\Vert\alpha^k\Vert_\infty^2)_{n\in\overline{\mathbb N}}$ is uniformly integrable.

\vspace{.5em}
By Corollary \ref{cor:R2I-Integr}.\ref{cor:R2I-Integr-ii} we know that the sequence $\big([ {\rm R}_{2,I}\star \widetilde{\mu}^{(X^k,\mathbb G^k)}]_{\infty} \big)_{k\in\overline{\mathbb N}}$ is uniformly integrable.
Therefore, using de La Vall\'ee Poussin's criterion, there exists a moderate Young function $\phi$ such that 
\begin{align}\label{cond:dLVP-SqBr}
\sup_{k\in\overline \N}\mathbb E \Big[ \phi\Big( 
\big[ {\rm R}_{2,I}\star \widetilde{\mu}^{(X^k,\mathbb G^k)}\big]_{\infty}\Big)\Big]<\infty.
\end{align}
\cref{prop:Young_functions} yields that the map $\mathbb R_+\ni x\overset{\psi}{\longmapsto} \psi(x):=\phi(\frac12 x^2)$ is again moderate and Young.
We can apply now the Burkholder--Davis--Gundy (BDG) inequality \cite[cf.][Theorem 10.36]{he1992semimartingale} to the sequence of martingales $\big( {\rm R}_{2,I}\star \widetilde{\mu}^{(X^k,\mathbb G^k)}\big)_{k\in\overline{\mathbb N}}$ using the function $\psi$, and we obtain that
\begin{align*}
\sup_{k\in\overline\N} \mathbb E\bigg[\phi\bigg( \frac12 \sup_{s> 0} \big| ({\rm R}_{2,I}\star \widetilde{\mu}^{(X^k,\mathbb G^k)})_s \big|^2 \bigg)\bigg]
	&= 	\sup_{k\in\overline\N} \mathbb E\bigg[\psi\bigg(\sup_{s> 0} \big| ({\rm R}_{2,I}\star \widetilde{\mu}^{(X^k,\mathbb G^k)})_s \big| \bigg)\bigg]\\
	&\overset{\text{BDG}}{\le }
		C_\psi 		\sup_{k\in\overline \N}\mathbb E \Big[ \psi\Big( \big[ {\rm R}_{2,I}\star \widetilde{\mu}^{(X^k,\mathbb G^k)}\big]_{\infty}^{\frac{1}{2}} \Big)\Big]\\
	&=   C_\psi 	\sup_{k\in\overline \N}\mathbb E \bigg[ \phi\bigg( \frac12 \big[ {\rm R}_{2,I}\star \widetilde{\mu}^{(X^k,\mathbb G^k)}\big]_{\infty} \bigg)\bigg]
	\overset{\eqref{cond:dLVP-SqBr}}{<} \infty. \numberthis\label{eq:helpforPUT}
\end{align*}
Hence the sequence $\big( \sup_{s> 0} |{\rm R}_{2,I}\star \widetilde{\mu}^{(X^k,\mathbb G^k)}|^2\big)_{k\in\overline\N}$ is uniformly integrable, again from de La Vall\'ee Poussin's criterion.
Moreover, $({\rm Tr}[[X^k]_{\infty}])_{k\in\overline{\mathbb N}}$ is a uniformly integrable sequence; cf. \cref{lem:Optional-UI}.\ref{lem:Optional-UI-i}.
Using analogous arguments to the above inequality, we can conclude that the sequence $( \sup_{s>0} |X^{k}_s|^2)_{k\in\overline{\mathbb N}}$ is also uniformly integrable.
Hence, the family
$\big( \sup_{s>0} |X^{k}_s|^2+\sup_{s>0}|{\rm R}_{2,I}\star \widetilde{\mu}^{(X^k,\mathbb G^k)} |^2\big)_{k\in\overline\N}$
is uniformly integrable, which allows us to conclude. 
\end{proof}

\begin{lemma}\label{lem:R-2I-PUT}
The sequence  $({\rm R}_{2,I}\star \widetilde{\mu}^{(X^{k,d},\mathbb G^k)})_{k\in\overline{\mathbb N}}$ possesses the $\PUT$ property, consequently
\begin{align}\label{conv:R2I-PUT-QV}
\big({\rm R}_{2,I}\star\widetilde{\mu}^{(X^{k,d},\mathbb G^k)}, [{\rm R}_{2,I}\star\widetilde{\mu}^{(X^{k,d},\mathbb G^k)}]\big)
\sipmulti{2}
\big({\rm R}_{2,I}\star\widetilde{\mu}^{(X^{\infty,d},\mathbb G^\infty)},[{\rm R}_{2,I}\star\widetilde{\mu}^{(X^{\infty,d},\mathbb G^\infty)}]\big).
\end{align}
\end{lemma}

\begin{proof}
By \eqref{eq:helpforPUT}, we obtain that the martingale sequence $({\rm R}_{2,I}\star \widetilde{\mu}^{(X^{\infty,d},\mathbb G^\infty)})_{k\in\overline{\mathbb N}}$ is $\mathbb L^2-$bounded, which allows us further to conclude that 
\begin{align*}
&\hspace{-3em}\sup_{k\in\overline{\mathbb N}} \mathbb E\bigg[\sup_{s>0} \big|\Delta\big({\rm R}_{2,I}\star \widetilde{\mu}^{(X^{\infty,d},\mathbb G^\infty)}\big)_s\big| \bigg]\\
&\overset{\textcolor{white}{\text{The process}}}{\underset{\textcolor{white}{\text{is \cadlag}}}{\le}} \sup_{k\in\mathbb N} \mathbb E\bigg[ \sup_{s\le t} \big|\Delta\big({\rm R}_{2,I}\star \widetilde{\mu}^{(X^{\infty,d},\mathbb G^\infty)}\big)_s\big| 
	+ \sup_{s\le t} \big|\Delta\big({\rm R}_{2,I}\star \widetilde{\mu}^{(X^{\infty,d},\mathbb G^\infty)}\big)_{s-}\big|\bigg]\\ 
&\overset{\text{The process}}{\underset{\text{is \cadlag}}{\le}} 
	2\sup_{k\in\mathbb N} 	\mathbb E\bigg[ \sup_{s\le t} \big|\Delta\big({\rm R}_{2,I}\star \widetilde{\mu}^{(X^{\infty,d},\mathbb G^\infty)}\big)_s\big|\bigg]\\
&\overset{\textcolor{white}{\text{The process}}}{\underset{\textcolor{white}{\text{is \cadlag}}}{\le}}
	2\sup_{k\in\mathbb N} 	\mathbb E\bigg[ \sup_{s\le t} \big|\Delta\big({\rm R}_{2,I}\star \widetilde{\mu}^{(X^{\infty,d},\mathbb G^\infty)}\big)_s\big|^2\bigg]^{\frac{1}{2}}\overset{\mathbb L^2-\text{boundedness}}{<}\infty.
\end{align*}
Now, \eqref{conv:AppMemin} yields that the sequence converges in law, hence \cite[Corollary VI.6.30]{jacod2003limit} allow us to conclude that the sequence $\big({\rm R}_{2,I}\star \widetilde{\mu}^{(X^{k,d},\mathbb G^k)}\big)_{k\in\mathbb N}$ possesses the $\PUT$ property.
Then, we obtain the required convergence from \cite[Theorems VI.6.26 and VI.6.22.(c)]{jacod2003limit}.
\end{proof}

\subsection{\texorpdfstring{Step $3$ is valid}{Step 3 is valid}}\label{subsubsection:Sufficient}

The following result provides a sufficient criterion for showing that a martingale $L$ is orthogonal to the space generated by another martingale $X$. 
We adopt the notation of the previous section.
\begin{proposition}\label{prop:SCO}
Let $X$ be an $\mathbb F-$quasi--left--continuous, $\mathbb R^\ell-$valued $\mathbb F-$martingale and $L$ be a uniformly integrable, $\mathbb R-$valued $\mathbb F-$martingale.
Assume that
\begin{enumerate}[label={\rm(\roman*)},leftmargin=*]
	\item\label{prop:SCO-i} 
		$[X,L]$ is a uniformly integrable $\Fil-$martingale, where 
		$[X,L]:=([X^1,L], \dots, [X^\ell,L])^\top.$
	\item\label{prop:SCO-ii} 
		$\mathcal{I}$ is a family of subsets of $\R{\ell}$ such that $\sigma(\mathcal{I})=\mathcal{B}(\R{\ell})$ and the martingale ${\rm R}_{2,A} \star \mutilde^{(X,\Fil)}$ is well--defined for every $A\in\mathcal I.$
		Moreover, $\big[ {\rm R}_{2,A} \star \mutilde^{(X,\Fil)}, L \big]$ is a uniformly integrable $\Fil-$martingale, for every $A\in \mathcal{I}$.
	\item\label{prop:SCO-iii}
		$|\Delta L|\mu^X\in\widetilde{\mathcal A}_{\sigma}(\mathbb F)$\footnote{For the notation, recall the discussion before Definition \ref{CondFPredProj}.}.
\end{enumerate}
Then, we have that
\begin{align}\label{eq:suffJumps}
\langle X^{c,i}, L^{c} \rangle^\Fil = 0,\text{ for every } i=1,\dots\ell,\text{ and }
M_{\mu^X}[\Delta L | \widetilde{\mathcal P}^\Fil]= 0.
\end{align}
\end{proposition}

\begin{proof}
By condition \ref{prop:SCO-i} we have that $[X^i,L]$ is a process of class {\rm (D)}, see \cite[Definition I.1.46]{jacod2003limit}, for every $i=1,\dots,\ell.$
Hence, by the Doob--Meyer decomposition we obtain that 
\begin{align}\label{SCO:AnBrNull}
\langle X^i, L\rangle = 0, \text{ for every } i=1,\dots,\ell.
\end{align} 

We are going to translate condition \ref{prop:SCO-ii} into the following one
\begin{align}\label{SCO:EqforW}
\mathbb E\bigg[	\int_{(0,\infty)\times \mathbb R^{\ell}} W(s,x) {\rm R}_2(x) \Delta L_s \mu^X(\ds,\dx)  \bigg] = 0, \text{ for every $\mathbb F-$measurable function }W.
\end{align}
Before we proceed, recall that we have assumed $X$ to be an $\mathbb F-$quasi--left--continuous martingale.
Thus, by \cite[Corollary II.1.19]{jacod2003limit} and Remark \ref{rem:PurDisJum} we can conclude that
\begin{align}\label{SCO:JumpEq}
\Delta \big({\rm R}_{2,A}\star\mutilde^{(X,\mathbb F)}\big)_s
	= {\rm R}_2\big( \Delta X^d_s\big)\mathds{1}_A(\Delta X^{d}_s)
	= {\rm R}_2 \big(\Delta X_s\big) \mathds{1}_A(\Delta X_s)
	= \Delta \big({\rm R}_{2,A}*\mu^{X}\big)_s.
\end{align}
By the martingale property of $\big[ {\rm R}_{2,A} \star \mutilde^{(X,\Fil)}, L \big]$, we obtain for every $0\le t<u<\infty$ and every $C\in\mathcal F_t$, that
\begin{align*}
0 	& = \mathbb E\Big[ \mathds{1}_C\, 	\mathbb E \big[\,  [{\rm R}_{2,A}\star \widetilde{\mu}^{(X,\mathbb F)}, L ]_u \big| \mathcal F_t \big]    - \mathds{1}_C\, [{\rm R}_{2,A}\star \widetilde{\mu}^{(X,\mathbb F)}, L ]_t\Big]\\
	& = \mathbb E\Big[ \mathds{1}_C\, 	[{\rm R}_{2,A}\star \widetilde{\mu}^{(X,\mathbb F)}, L ]_u - \mathds{1}_C\, [{\rm R}_{2,A}\star \widetilde{\mu}^{(X,\mathbb F)}, L ]_t\Big] 
	= \mathbb E\bigg[ \mathds{1}_C 	\sum_{t<s\le u} \Delta \big({\rm R}_{2,A}\star \widetilde{\mu}^{(X,\mathbb F)}\big)_s \Delta L_s \bigg]\\
	& \hspace{-0.5em}\overset{\eqref{SCO:JumpEq}}{=} 
		\mathbb E\bigg[ \mathds{1}_C 	\sum_{t<s\le u} \Delta \big({\rm R}_{2,A}* \mu^X\big)_s \Delta L_s \bigg] \ \ 
	\hspace{-0.5em}\overset{\eqref{SCO:JumpEq}}{=} 
	 	\mathbb E\bigg[	\int_{(t,u]\times A} \mathds{1}_C {\rm R}_2(x) \Delta L_s \mu^X(\ds,\dx)  \bigg],
\end{align*}
where in the third equality we have used that ${\rm R}_{2,A}\star \widetilde{\mu}^{(X,\mathbb F)}\in \mathcal H^{2,d}(\mathbb F,\infty;\mathbb R)$ and \cite[Theorem I.4.52]{jacod2003limit}.
Moreover, observe that
\begin{align*}
\widetilde{\mathcal P}^{\mathbb F}=\mathcal P^{\mathbb F}\otimes \mathcal B(\mathbb R^\ell) 
		&	= \sigma\big( P\times A, \text{ where }P\in\mathcal P^{\mathbb F} \text{ and } A\in\mathcal B(\mathbb R^\ell) \big)\\
		& 	\hspace{-2.4em}\overset{\text{\cite[Theorem I.2.2]{jacod2003limit}}}{=} \sigma\big( C\times(t,u]\times A, \text{ where }0\le t<u, C\in\mathcal F_t \text{ and } A\in\mathcal I \big),
\end{align*}
where we have from \ref{prop:SCO-ii} that $\sigma(\mathcal I)=\mathcal B(\mathbb R^\ell).$
Hence, by a monotone class argument we can conclude that condition \eqref{SCO:EqforW} holds. 

\vspace{.5em}
The next observation is that the function $\widetilde{\Omega}\ni (\omega,s,x) \longmapsto {\rm R}_2(x) \in\mathbb R$ is $\mathbb F-$predictable, since it is deterministic and continuous. 
Moreover, the function ${\rm R}_2$ is positive $M_{\mu^X}-a.e.$, recall that $\mu^X$ has been defined using the random set $[\Delta X\neq 0].$
Therefore, we obtain 
\begin{align}\label{SCO:EqforU}
\mathbb E\bigg[	\int_{(0,\infty)\times \mathbb R^{\ell}} U(s,x) \Delta L_s \mu^X(\ds,\dx)  \bigg] = 0, \text{ for every $\mathbb F-$measurable function }U,
\end{align}
by substituting in \eqref{SCO:EqforW} the $\mathbb F-$predictable function $W(\omega,s,x)$ with the $\mathbb F-$predictable function $\frac{U(\omega,s,x)}{{\rm R}_2(x)},$ where $U$ is an arbitrary $\mathbb F-$predi\-ctable function.

\vspace{0.5em}
Now, since we have assumed that $|\Delta L|\mu^X\in \widetilde{\mathcal A}_{\sigma}(\mathbb F)$, in view of \eqref{SCO:EqforU} and Definition \ref{CondFPredProj} we obtain that the conditional $\mathbb F-$predictable projection of $L$ on $\mu^X$ is well--defined and is identically zero, \emph{i.e.}
\begin{align}\label{SCO:RNDerisNull}
M_{\mu^X}\big[\Delta L|\widetilde{\mathcal P}^{\mathbb F}\big] = 0.
\end{align}

\vspace{0.5em}
We proceed now to show the validity of the first $\ell$ conditions in \eqref{eq:suffJumps}.
By \cite[Theorem 13.3.16]{jacod2003limit}, we obtain that\footnote{We abuse notation and denote $\widetilde{\Omega}\ni (\omega,s,x) \overset{\pi^i}{\longmapsto} x^i\in\mathbb R$.}
\begin{align*}
\langle X^{d,i}, L\rangle 
	= \langle \pi^i\star \widetilde{\mu}^{(X,\mathbb F)}, L\rangle
	= \big(\pi^i M_{\mu^X}[\Delta L|\widetilde{\mathcal P}^{\mathbb F}] \big)*\nu^{(X,\mathbb F)} \overset{\text{\eqref{SCO:RNDerisNull}}}{=} 0, \text{ for every }i=1,\dots,\ell.
	\numberthis\label{SCO:PDareOrth}
\end{align*}
The combination of \eqref{SCO:AnBrNull} and \eqref{SCO:PDareOrth} yields
\begin{align*}
\langle X^{c,i}, L^c\rangle
	= \langle X^{i}, L\rangle - \langle X^{d,i}, L^d\rangle^{\mathbb F} 
	\overset{\eqref{SCO:AnBrNull}}{\underset{\eqref{SCO:PDareOrth}}{=}} 0, \text{ for every }i=1,\dots,\ell,
\end{align*}
where for the first equality we have applied \cite[Theorem I.4.52]{jacod2003limit}. 
\end{proof}


\subsection{\texorpdfstring{Step $2$ is valid}{Step 2 is valid}}\label{sec:StepB}

Now that we have obtained a family of converging martingales by Proposition \ref{prop:AppMemin}, we proceed by proving some technical lemmata which are going to be useful in the proof of Theorem \ref{RobMartRep}.

\vspace{0.5em}
Recall that our aim in \ref{StepB} is to prove that 
\begin{align}\label{Cond:XiisZero}
[Y^\infty, \overline N] \text{ is an $\mathbb F-$martingale,} 
\end{align}
for every weak--limit point $\overline N$ of $(N^k)_{k\in\mathbb N}$ and for some filtration $\mathbb F$ which includes $\mathbb G^\infty$, and may depend on $\overline N$. 	
In the next few lines, we are going to explain why this is sufficient (for showing that the limit of $(\langle N^k\rangle)_{k\in\mathbb N}$ equals zero) and how the filtration $\mathbb F$ is going to be determined. 

\vspace{0.5em}
Observe that by the orthogonal decomposition of $Y^k$ with respect to $(X^{k,c},\mu^{X^{k,d}},\mathbb G^k)$, we have for every $k\in\mathbb N$ 
\begin{align}\label{Iden:AngleYN-N}
\langle Y^k,N^k \rangle =\langle N^k \rangle.
\end{align}
This identity is the link with what follows. 
In Lemma \ref{lem:NisTight}, we show that the sequence $(N^k)_{k\in\N}$ is tight, while the sequence $(\langle N^k \rangle)_{k\in\N}$ is $C-$tight.
Thus, there exists a subsequence $(k_l)_{l\in\N}$ such that $(N^{k_l})_{l\in\mathbb N}$ and $(\langle N^{k_l} \rangle)_{l\in\mathbb N}$ converge jointly and we denote the respective weak--limit points by $\overline N$ and $\Xi$, \textit{i.e.}
\begin{align}\label{conv:weaklimit}
(N^{k_l}, \langle N^{k_l}\rangle)\xrightarrow[l\to\infty]{\hspace{0.2cm}\mathcal L \hspace{0.2cm}}(\overline N, \Xi).
\end{align}
The last convergence will enable us to prove, see Corollary \ref{cor:Conv-Theta}, that
\begin{align}\label{conv:QCYN-Xi}
[Y^{k_l},N^{k_l}]-\langle N^{k_l}\rangle \xrightarrow[l\to\infty]{\hspace{0.2cm}\mathcal L \hspace{0.2cm}} [Y^\infty,\overline N] - \Xi.
\end{align}
Using classical arguments, we can prove that the limit process is a uniformly integrable martingale with respect to $\mathbb F^{(Y^\infty,\overline N, \Xi)}$.
Then, we can conclude that $\Xi=0$, if \eqref{Cond:XiisZero} is valid for the filtration $\mathbb F^{(Y^\infty,\overline N, \Xi)}$, since in this case $\Xi$ is an $\mathbb F^{(Y^\infty,\overline N, \Xi)}-$predictable $\mathbb F^{(Y^\infty,\overline N, \Xi)}-$martingale of finite variation, see \cite[Corollary I.3.16]{jacod2003limit}.

\vspace{0.5em}
However, the problem is that the filtration $\mathbb F^{(Y^\infty,\overline N, \Xi)}$ does not necessarily contain $\mathbb G^\infty.$
Before we proceed let us fix an arbitrary $I\in\mathcal J(X^\infty)$ and an arbitrary $G\in\mathcal G_\infty^\infty.$
In order to enlarge the filtration with respect to which the limiting process in \eqref{conv:QCYN-Xi} is a martingale, we will use the sequences $(\Theta^{k,I})_{k\in\mathbb N}$ and $(\Theta^{k,I,G})_{k\in\mathbb N}$, which are defined as
\begin{align}
\Theta^{k,I}&:=
		\big(
						X^{k}, 										[X^{k}]-\langle X^{k}\rangle, 
						Y^k,										N^{k}, 										
						[X^k,N^k],									[{\rm R}_{2,I}\star \widetilde{\mu}^{(X^{k,d},\mathbb G^k)}, N^k], 
						[Y^{k},N^{k}]-\langle N^{k}\rangle
		\big)^\top\label{def:Thetak},\\
\Theta^{k,I,G}&:=
		\big(
				(\Theta^{k,I})^\top, 	\E [\mathds{1}_G|\mathcal G^{k}_\cdot]
		\big)^\top. \label{def:ThetakG}
\end{align} 
We will prove initially that 
\begin{multline}\label{def:Theta-infG}
\Theta^{\infty,I,G}:= \\
	\big(	
			X^\infty, 							[X^\infty]-\langle X^\infty\rangle^{\mathbb G^\infty}, 
			Y^\infty,							\overline N,
			[X^\infty,\overline N],				[{\rm R}_{2,I}\star \widetilde{\mu}^{(X^{\infty,d},\mathbb G^\infty)}, \overline N],
			[Y^\infty,\overline N]- \Xi,
			\E [\mathds{1}_G|\mathcal G^{\infty}_\cdot]
	\big)^\top,
\end{multline}
\emph{i.e.} the weak--limit of $(\Theta^{k_l,I,G})_{l\in\mathbb N}$, is an $\mathbb F^{M^{\infty,G}}-$martingale, where 
\begin{align}\label{def:M-infG}
M^{\infty,G}:=(X^\infty, \langle X^\infty\rangle^{\mathbb G^\infty}, Y^\infty, \overline N, \Xi, \mathbb E[\mathds{1}_G|\mathcal G^\infty_\cdot]).
\end{align}
Since the set $G$ was arbitrarily chosen, we will deduce in Proposition \ref{prop:Fmart} the martingale property of
\begin{align}\label{def:Theta-inf}
\Theta^{\infty,I}:=\big(	X^\infty, 						[X^\infty]-\langle X^\infty\rangle^{\mathbb G^\infty}, 
						Y^\infty,							\overline N, 						
						[X^\infty,\overline N],				[{\rm R}_{2,I}\star \widetilde{\mu}^{(X^{\infty,d},\mathbb G^\infty)}, \overline N],
						[Y^\infty,\overline N]- \Xi
			\big)^\top,
\end{align}
with respect to the filtration 
\begin{align}\label{def:FilF}
\mathbb F:=\mathbb G^\infty\vee \mathbb F^{M^\infty}, \text{ for } M^\infty:=(X^\infty, \langle X^\infty\rangle^{\mathbb G^\infty}, Y^\infty, \overline N, \Xi).
\end{align}
Observe that the filtration $\mathbb F$ depends on $\overline N$ and $\Xi$, but we notationally suppressed this dependence.
Moreover, it does not necessarily hold that $\mathbb F$ is right--continuous. 
However, in view of \cite[Theorem 2.46]{he1992semimartingale} and the \cadlag property of $\Theta^\infty$, we can conclude also that $\Theta^\infty$ is an $\mathbb F_{\text{\tiny $+$}}-$martingale, where 
$\mathbb F_{\text{\tiny $+$}}:=(\mathcal F_{t\text{\tiny $+$}})_{t\in\mathbb R_+}$ and $\mathcal F_{t\text{\tiny $+$}}:=\bigcap_{s>t}\mathcal F_s.$

\subsubsection{\texorpdfstring{$\Theta^{\infty,I}$ is an $\mathbb F$--martingale}{Theta is an F--martingale}}\label{sub-sub-sec:StepB-Prep}

In view of the discussion above, we need to show that the family $(\Theta^{k,I,G})_{k\in\N}$ is tight and uniformly integrable for every $I\in\mathcal J(X^\infty)$ and every $G\in\mathcal G_\infty^\infty$. These results are proved in Lemmata \ref{lem:Theta-i} and \ref{lem:Theta-ii}.
Before that, we present some necessary results. 

\begin{lemma}\label{lem:NisTight}
Assume the setting of Theorem \ref{RobMartRep}.
Then, for the sequence $(N^k)_{k\in\N}$ the following are true
\begin{enumerate}[label={\rm (\roman*)}]
\item\label{lem:N-Tight-i} The sequence $(\langle N^k\rangle)_{k\in\mathbb N^k}$ is $C-$tight in  $\mathbb D(\mathbb R)$.
\item\label{lem:N-Tight-ii} The sequence $(N^k)_{k\in\mathbb N}$ is tight in $\mathbb D(\mathbb R)$.
\item\label{lem:N-Tight-iii} The sequence $\big((N^k,\langle N^k\rangle)\big)_{k\in\mathbb N}$ is tight in $\mathbb D(\mathbb R^2)$.
\item\label{lem:N-Tight-iv} The sequence $(N^k)_{k\in\mathbb N}$ is $\mathbb L^2-$bounded, \emph{i.e.} $\sup_{k\in\N}\mathbb E\Big[ \sup_{t\in[0,\infty]} |N^k_t|^2\Big]<\infty.$
\item\label{lem:N-Tight-v} The sequence $(\langle N^k\rangle_\infty)_{k\in\mathbb N^k}$ is uniformly integrable.
\item\label{lem:N-Tight-vi} The sequence $([ N^k]_\infty)_{k\in\mathbb N^k}$ is $\mathbb L^1-$bounded, \emph{i.e.} $\sup_{k\in\mathbb N}\mathbb E\big[ \,[N^k]_\infty\big]<\infty.$
\end{enumerate}
\end{lemma}
\begin{proof}
\begin{enumerate}[leftmargin=0cm, itemindent=1em, align=left]
\item[\ref{lem:N-Tight-i}]
By the construction of the orthogonal decompositions of $Y^k$ with respect to $(X^{k,c},\mu^{X^{k,d}},\mathbb G^k)$ and Corollary \ref{cor:orthogSpace}, we obtain that 
\begin{align}\label{id:OrthMart}
 \langle Z^k\cdot X^{k,c} + U^k\star \widetilde{\mu}^{(X^{k,d},\mathbb G^k)}, N^k\rangle = 0, \text{ for every }k\in\mathbb N.
\end{align}
By conditions \ref{MFilweak}, \ref{Mfinalrv} and Proposition \ref{MeminProp1}, we have
$$ (Y^k, \mathbb G^k)\sipmulti{}(Y^\infty,\mathbb G^\infty),$$ 
which allows us to apply Theorem \ref{MeminCorollary} and therefore obtain by part \ref{MeminCorollaryJ1P} the convergence
\begin{align}\label{conv:MeminCorY}
 (Y^k, [Y^k], \langle Y^k\rangle)\sipmulti{3}(Y^\infty, [Y^\infty], \langle Y^\infty \rangle).
\end{align}
Hence, by Assumption \ref{MFilqlc} and convergence \eqref{conv:MeminCorY}, we get that the sequence $(\langle Y^k\rangle)_{k\in\mathbb{N}}$ is C$-$tight, see \cite[Definition VI.3.25]{jacod2003limit}. 
Moreover, \eqref{id:OrthMart} implies that
\begin{align}\label{StrongMajor}
\langle Y^k\rangle =\,
		\langle Z^k\cdot X^{k,c} + U^k\star \widetilde{\mu}^{(X^{k,d},\mathbb G^k)} \rangle		+		\langle N^k\rangle, \text{ for every } k\in\N,
\end{align}
which in turn yields that the process $\langle Y^k\rangle$ strongly majorises both $\langle Z^k\cdot X^{k,c} + U^k\star \widetilde{\mu}^{(X^{k,d},\mathbb G^k)}\rangle $ and $\langle N^k\rangle$, for every $k\in\N$, see \cite[Definition VI.3.34]{jacod2003limit}.
We can conclude thus the C$-$tightness of $(\langle N^k\rangle)_{k\in\mathbb{N}}$ by  \cite[Proposition~VI.3.35]{jacod2003limit}.
\vspace{0.5em}
\item[\ref{lem:N-Tight-ii}]
Using \cite[Theorem VI.4.13]{jacod2003limit} to obtain the tightness of $(N^k)_{k\in\N}$, it suffices to show that $(\langle N^k \rangle)_{k\in\N}$ is C$-$tight and $(N_0^k)_{k\in\N}$ is tight.
The first statement follows from \ref{lem:N-Tight-i}, while for the second one we get that $N_0^k=0$, by the definition of the orthogonal decomposition of $Y^k$ with respect to $(X^{k,c},\mu^{X^{k,d}},\mathbb G^k),$ for every $k\in\mathbb N.$ 
\vspace{0.5em}
\item[\ref{lem:N-Tight-iii}]
This is immediate in view of \cite[Corollary VI.3.33]{jacod2003limit} and \ref{lem:N-Tight-i}--\ref{lem:N-Tight-ii}.
\vspace{0.5em}
\item[\ref{lem:N-Tight-iv}]
We have by Doob's $\mathbb{L}^2-$inequality, see \cite[Theorem 5.1.3]{jacod2003limit} and \cite[Theorem 6.8]{he1992semimartingale}, that
\begin{align}\label{ineq:DoobIneq}
\E\bigg[\,\sup_{t\in[0,\infty]}\big|N_t^k\big|^2\bigg]
	\leq 4\E\Big[\big|N_{\infty}^k\big|^2\Big] 
	= 4\E\Big[\langle N^k\rangle_{\infty}\Big]
	\leq 4\E\Big[\langle Y^k\rangle_{\infty}\Big],
\end{align}
by identity \eqref{StrongMajor}.
By convergence \eqref{conv:MeminCorY}, compare Theorem \ref{MeminCorollary}.\ref{BracketConvinL}, we obtain that the sequence $(\langle Y^k\rangle_\infty)_{k\in\N}$ is uniformly integrable and in particular $\mathbb L^1-$bounded.

\vspace{0.5em}
\item[\ref{lem:N-Tight-v}]
Since $(\langle N^k\rangle_\infty)_{k\in\N}$ is strongly majorized by $(\langle Y^k\rangle_\infty)_{k\in\N}$ which is uniformly integrable, see the comments in the proof of \ref{lem:N-Tight-iv}, we can conclude by \cite[Theorem 1.7]{he1992semimartingale}.
\vspace{0.5em}
\item[\ref{lem:N-Tight-vi}]
The identity $\mathbb E\big[ [N^k]_\infty]=\mathbb E\big[\langle N^k\rangle_\infty\big]$ for every $k\in\mathbb N$ and the uniform integrability of the sequence $(\langle N^k\rangle_\infty)_{k\in\mathbb N}$ which implies the $\mathbb L^1-$boundedness of the sequence allow us to conclude.
\end{enumerate}
\end{proof}

\begin{lemma}\label{lem:Xk-Yk-PUT}
Let conditions \ref{MFilqlc}, \ref{MSI}, \ref{MFilweak} and \ref{Mfinalrv} hold. 
Then the sequences $(X^k)_{k\in\mathbb N}$ and $(Y^k)_{k\in\mathbb N}$ possess the $\PUT$ property.
\end{lemma}
\begin{proof}
The proof is analogous to the one of \cref{lem:R-2I-PUT}, in particular we are going to apply \cite[Corollary VI.6.30]{jacod2003limit}.
In view of conditions \ref{MFilqlc}, \ref{MSI}, \ref{MFilweak}, \ref{Mfinalrv}, Remark \ref{EquivalentExtended} and Proposition \ref{MeminProp1}, we have that the sequences $(X^k)_{k\in\overline{\mathbb N}}$ and $(Y^k)_{k\in\overline{\mathbb N}}$ satisfy Theorem \ref{MeminCorollary}, \emph{i.e.}
$$X^k\sipmulti{\ell}X^\infty \text{ and } Y^k\sipmulti{} Y^\infty.$$
Moreover,  $(X^k)_{k\in\mathbb N}$ is $\mathbb L^2-$bounded and for every $k\in\mathbb N$ the process $X^k$ is a \cadlag $\mathbb G^k-$martingale.
Therefore, using similar arguments to the ones in the proof of \cref{lem:R-2I-PUT}, we get that
\begin{align*}
\sup_{k\in\mathbb N} \mathbb E\bigg[ \sup_{s\le t} \big|\Delta X^k_s\big|\bigg] < \infty,
\end{align*}
which allows us to conclude.
The steps for the sequence $(Y^k)_{k\in\mathbb N}$ are completely similar, so that we omit them.
\end{proof}

%
%
%
%
\begin{lemma}\label{lem:Theta-i}
The family of random variables $\big(\big\Vert \Theta^{k,I,G}_t\big\Vert_1\big)_{k\in\N,t\in[0,\infty]}$ is uniformly integrable for every $I\in\mathcal J(X^\infty)$ and $G\in\mathbb G^\infty_\infty$, where, abusing notation, we have defined 
\begin{align*}
\big\Vert \Theta^{k,I,G}_t \big\Vert_1:=&\
		\		\big\Vert X^{k}_t \big\Vert_1 
			+ 	\big\Vert [X^{k}]_t- \langle X^{k}\rangle_t \big\Vert_1
			+ 	\Vert (Y^k_t, N^k_t)\Vert_1
			+	\Vert [X^k,N^k]_t \Vert_1\\
			&+ 	\big\Vert \big( [{\rm R}_{2,I}\star\widetilde{\mu}^{(X^{\infty,d},\mathbb G^\infty)},\overline N]_t,	[Y^{k}, N^{k}]_t-\langle N^{k}\rangle_t, 		\E [\mathds{1}_G|\mathcal G^{k}_t] \big) \big\Vert_1.
\numberthis\label{def:NormTheta-i}
\end{align*}
\end{lemma}
\begin{proof}
Let $I\in\mathcal J(X^\infty)$ and $G\in\mathbb G^\infty_\infty$ be fixed.
We are going to prove that for each summand of $\Vert \Theta^{k,I,G}_t\Vert_1$ the associated family indexed by $\{k\in\N,t\in[0,\infty]\}$ is uniformly integrable. 
Then, by \cite[Corollary 1.10]{he1992semimartingale}, we can conclude also for their sum $(\Vert \Theta^{k,I,G}_t\Vert_1)_{k\in\N,t\in[0,\infty]}$.

\begin{enumerate}[label={\bf ($\mathbf{\Theta}$\arabic*)}, itemindent=1em, align=left, leftmargin=0cm, label={\rm (\roman*)}]
\item\label{item:Xk-UI} Let $i\in\{1,\dots,\ell\}.$ 
Theorem \ref{MeminCorollary} and Burkholder--Davis--Gundy's inequality imply the $\mathbb L^2$--boundedness of $(\sup_{t\in[0,\infty]}|X^{k,i}|)_{k\in\N}.$
By de La Vall\'ee--Poussin's criterion, we obtain that  $(\sup_{t\in[0,\infty]}|X^{k,i}|)_{k\in\N}$ is uniformly integrable.
By the obvious domination $|X^{k,i}_t|\le \sup_{t\in[0,\infty]}|X^{k,i}_t|,$ for every $t\in[0,\infty],$ and \cite[Theorem 1.7]{he1992semimartingale}, we can conclude in particular the uniform integrability of the family $(|X^{k,i}_t|)_{k\in\N,t\in[0,\infty]}.$
By \cite[Corollary 1.10]{he1992semimartingale}, we conclude that $(\Vert X^{k}_t\Vert_1)_{k\in\N,t\in[0,\infty]}$ is uniformly integrable.
\vspace{0.5em}
\item Let $i,j\in\{1,\dots,\ell\}.$ 
By Theorem \ref{MeminCorollary} and Lemma \ref{UIplusL2Bounded}, we obtain the uniform integrability of the sequence $({\rm Var}([X^{k,i},X^{k,j}])_\infty)_{k\in\N}.$ 
Hence the sequence $(|[X^{k,i},X^{k,j}]_t|)_{k\in\N,t\in[0,\infty]}$ is also uniformly integrable, in view of the domination $|[X^{k,i},X^{k,j}]_t|\le {\rm Var}([X^{k,i},X^{k,j}])_\infty.$
Then, $(\Vert [X^k]_t\Vert_1)_{k\in\N,t\in[0,\infty]}$ is also uniformly integrable as the sum of uniformly integrable families.

\vspace{0.5em}
\noindent For the uniform integrability of $(\langle X^{k,i},X^{k,j}\rangle_t)_{k\in\N,t\in[0,\infty]}$, we can use arguments completely analogous to the ones above. 
We have only to mention that Lemma \ref{UIplusL2Bounded} is valid when we substitute the quadratic covariation by the predictable quadratic covariation.
Therefore we can conclude by the inequalities
\begin{align*}
|[X^{k,i},X^{k,j}]_t-\langle X^{k,i},X^{k,j}\rangle_t|
	\le |[X^{k,i},X^{k,j}]_t|+|\langle X^{k,i},X^{k,j}\rangle_t| 
	\le \Vert {\rm \Var}([X^k])_\infty \Vert_1 + \Vert {\rm \Var}(\langle X^k\rangle)_\infty \Vert_1.
\end{align*}
%
\item 
We can conclude the uniform integrability of $(Y^k)_{k\in\mathbb N}$ by arguments analogous to \ref{item:Xk-UI}, since the sequence $(Y^k)_{k\in\overline{\mathbb N}}$ satisfies the assumptions of Theorem \ref{MeminCorollary}. 
\vspace{0.5em}
\item
The uniform integrability of $(N^k_t)_{k\in\N,t\in[0,\infty]}$ is immediate by Lemma \ref{lem:NisTight}.\ref{lem:N-Tight-iv} and de La Vall\'ee Poussin's criterion.
\vspace{0.5em}
\item 
In Corollary \ref{cor:R2I-Integr}.\ref{cor:R2I-Integr-ii} we have obtained that the sequence $([ {\rm R}_{2,I}\star \widetilde{\mu}^{(X^{k,d},\mathbb G^k)}]_{\infty} )_{k\in\overline{\mathbb N}}$ is uniformly integrable.
Moreover, the sequence $( [N^k]_\infty)_{k\in\mathbb N}$ is $\mathbb L^1-$bounded, by Lemma \ref{lem:NisTight}.\ref{lem:N-Tight-vi}.
Now we can conclude the uniform integrability of the sequence $({\rm Var}([{\rm R}_{2,I}\star \widetilde{\mu}^{(X^{k,d},\mathbb G^k)}, N^k])_{\infty})_{k\in\mathbb N}$ by Lemma \ref{UIplusL2Bounded}.
This further implies the uniform integrability of $([{\rm R}_{2,I}\star \widetilde{\mu}^{(X^{k,d},\mathbb G^k)}, N^k]_t)_{k\in\mathbb N, t\in[0,\infty]}$.
\vspace{0.5em}
\item
We prove first the uniform integrability of the family $([Y^k,N^k]_t)_{k\in\N,t\in[0,\infty]}$. 
We can obtain it by applying Lemma \ref{lem:NisTight} to the uniformly integrable sequence $({\rm Tr}[ [Y^k]_\infty])_{k\in\mathbb N}$ and the $\mathbb L^1-$bounded sequence 
$( [N^k]_\infty)_{k\in\mathbb N}$; see Lemma \ref{lem:NisTight}.\ref{lem:N-Tight-vi}. 
The uniform integrability of $(\langle N^k\rangle_\infty)_{k\in\infty}$ is provided by Lemma \ref{lem:NisTight}.\ref{lem:N-Tight-v}.
Then
\begin{align*}
|[Y^k,N^k]_t-\langle N^k\rangle_t| 
	&\le |[Y^k,N^k]_t| + |\langle N^k\rangle_t| \\
	&\le {\rm Var}\big([Y^k,N^k]\big)_t + {\rm Var}\big(\langle N^k\rangle\big)_t 
	 \le {\rm Var}\big([Y^k,N^k]\big)_\infty + {\rm Var}\big(\langle N^k\rangle\big)_\infty,
\end{align*}
and we can now conclude by \cite[Theorem 1.7]{he1992semimartingale}.
\vspace{0.5em}
\item
Observe that for every $p>1$ and for every $k\in\N, t\in[0,\infty]$ it holds
\begin{align*}
\E\big[(\E[\mathds{1}_G|\mathcal G^k_t])^p\big] \overset{\text{Jensen ineq.}}{\le} \E\big[ \E[(\mathds{1}_G)^p|\mathcal G^k_t]\big] = \E[\mathds{1}_G]=\mathbb P(G).
\end{align*}
Then, we can conclude by de La Vall\'ee--Poussin's criterion. 
We could have, more generally, used that $0\le \E[\mathds{1}_G|\mathcal G_t^k]\le 1$, for every ${k\in\N,t\in[0,\infty]}$, see \cite[Theorem 1.18]{he1992semimartingale}. \qedhere
\end{enumerate}
\end{proof}
%
%
%
%
\begin{lemma}\label{lem:Theta-ii}
The sequence $\big(\Theta^{k,I,G}\big)_{k\in\N}$ is tight in 
$\mathbb D(\mathbb R^{(\ell+1)\times\ell}\times \mathbb R^{1\times2}\times\mathbb R^{\ell}\times\mathbb R^{1\times3})$
for every $I\in\mathcal J(X^\infty)$ and for every $G\in\mathcal G_\infty^\infty$.
\end{lemma}

\begin{proof}
Let $I\in\mathcal J(X^\infty)$ and $G\in\mathcal G_\infty^\infty$ be fixed.
We claim that the sequence $\Phi^k:=(X^{k},Y^{k},N^{k})$ is tight in $\mathbb D(\R{\ell}\times\R{1\times2})$, and that the tightness of $(\Phi^k)_{k\in\N}$ is sufficient to show the tightness of $(\Theta^{k,I,G})_{k\in\N}$.


\vspace{.5em}
The space $\mathbb D(\mathbb R^{(\ell+1)\times\ell}\times \mathbb R^{1\times2}\times\mathbb R^{\ell}\times\mathbb R^{1\times3})$ is Polish since it is isometric to $\mathbb D(\mathbb R^{\ell^2+2\cdot\ell+5})$ which is Polish, see \cite[Theorem VI.1.14]{jacod2003limit}.
Hence in this case tightness is equivalent to sequential compactness.
Therefore, it suffices to provide a weakly convergent subsequence $(\Theta^{k_{l_m},I,G})_{m\in\N}$ for every subsequence $(\Theta^{k_l,I,G})_{l\in\N}$.

\vspace{0.5em}
Let us therefore consider a subsequence $(\Theta^{k_l,I,G})_{l\in\N}$. 
Assuming the tightness of $(\Phi^k)_{k\in\N}$, there exists a weakly convergent subsequence $(\Phi^{k_{l_m}})_{m\in\N},$ converging to, say, $(X^\infty,Y^\infty, \widetilde N).$
Moreover, by \cref{lem:Xk-Yk-PUT}, we have that the sequences $(X^k)_{k\in\N}$, $(Y^k)_{k\in\N}$ possess the $\PUT$ property. 
In particular, the subsequences $(X^{k_{l_m}})_{m\in\N}$, $(Y^{k_{l_m}})_{m\in\N}$ possess the $\PUT$ property.
On the other hand, $(N^k)_{k\in\N}$ is $\mathbb L^2-$bounded and by using arguments completely analogous to the ones used in Lemma \ref{lem:R-2I-PUT}, we get that the sequence $(N^{k_{l_m}})_{m\in\N}$ possesses also the $\PUT$ property.
Finally, by Lemma \ref{lem:R-2I-PUT} we have that the sequence $({\rm R}_{2,I}\star\widetilde{\mu}^{(X^{k_{l_m},d},\mathbb G^{k_{l_m}})})_{m\in\mathbb N}$ possesses the $\PUT$ property.
By Theorem \ref{MeminCorollary}, we therefore obtain the convergence
\begin{align*}
(X^k,[X^k],\langle X^k\rangle)
\xrightarrow{\hspace{0.2cm}(\J_1(\mathbb R^{\ell}\times \R{\ell\times\ell}\times\R{\ell\times\ell}),\mathbb P)\hspace{0.2cm}}
(X^\infty,[X^\infty], \langle X^\infty\rangle^{\mathbb G^\infty}).
\end{align*}
Hence, by \cite[Theorem VI.6.26]{jacod2003limit} and by Proposition \ref{prop:AppMemin}, we also get the convergence

\vspace{0.5em}
\hfil$\big(X^{k_{l_m}},[X^{k_{l_m}}]-\langle X^{k_{l_m}}\rangle, Y^{k_{l_m}}, N^{k_{l_m}}, [X^{k_{l_m}},N^{k_{l_m}}], [{\rm R}_{2,I}\star\widetilde{\mu}^{(X^{k_{l_m},d},\mathbb G^{k_{l_m}})},N^{k_{l_m}}], [Y^{k_{l_m}},N^{k_{l_m}}] \big)$

\vspace{0.3em}
{\hfil $\Big\downarrow \mathcal L$}

\vspace{0.3em}
\hfil$\big(X^{\infty},[X^\infty]-\langle X^\infty\rangle^{\mathbb G^\infty}, Y^{\infty}, \widetilde N, [X^\infty,\widetilde N], [{\rm R}_{2,I}\star\widetilde{\mu}^{(X^{\infty,d},\mathbb G^\infty)},\widetilde N], [Y^\infty,\widetilde N] \big).$

\vspace{0.5em}
In view of the $C-$tightness of $(\langle N^k\rangle)_{k\in\N},$ see Lemma \ref{lem:NisTight}, we can pass to a further subsequence $(k_{l_{m_{n}}})_{n\in\N}$ so that 
$$\langle N^{k_{l_{m_{n}}}}\rangle\xrightarrow[n\to\infty]{\hspace{0.2cm}\mathcal L \hspace{0.2cm}} \widetilde \Xi.$$
Hence, we can finally obtain a jointly weakly convergent subsequence
$$\Theta^{k_{l_{m_{n}}},I,G} 
\xrightarrow[n\to\infty]{\hspace{0.2cm}\mathcal L \hspace{0.2cm}} 
\big(X^{\infty},[X^\infty]-\langle X^\infty\rangle^{\mathbb G^\infty}, Y^{\infty}, \widetilde N, [X^\infty,\widetilde N], [{\rm R}_{2,I}\star\widetilde{\mu}^{(X^{\infty,d},\mathbb G^\infty)},\widetilde N], [Y^\infty,\widetilde N]-\widetilde{\Xi}, \mathbb E[\mathds{1}_G|\mathcal G^\infty_\cdot] \big),$$
where we have used that $(X^k, \mathbb E[\mathds{1}_G|\mathcal G^k_\cdot])\xrightarrow[k\to\infty]{\hspace{0.2cm}(\J_1(\mathbb R^{\ell+1}),\mathbb P)\hspace{0.2cm}}(X^\infty, \mathbb E[\mathds{1}_G|\mathcal G^\infty_\cdot])$ to conclude.

In order to prove our initial claim that $(\Phi^k)_{k\in\N}$ is tight, we will apply \cite[Theorem 4.13]{jacod2003limit} to the $\mathbb L^2-$bounded sequences $(X^k)_{k\in\N},(Y^k)_{k\in\N},$ and $(N^k)_{k\in\N}.$
The sequences $(X^k_0)_{k\in\N},(Y^k_0)_{k\in\N},$ and $(N^k_0)_{k\in\N}$ are clearly tight in $\R{\ell}$, $\R{}$ and $\R{}$ respectively, as $\mathbb L^2-$bounded.
The sequences $\big({\rm Tr}\big[\langle X^k\rangle\big]\big)_{k\in\N}$ and $(\langle Y^k\rangle)_{k\in\N}$ are $C-$tight as a consequence of Theorem \ref{MeminCorollary} and the quasi--left--continuity of $X^\infty$ and $Y^\infty.$
The sequence $(\langle N^k\rangle)_{k\in\N}$ is $C-$tight by Lemma \ref{lem:NisTight}.
Finally, the sequence $(\Psi^k)_{k\in\N}$, where 
$$\Psi^k:={\rm Tr}\big[\langle X^k\rangle\big] + \langle Y^k\rangle + \langle N^k\rangle,$$
is $C-$tight as the sum of $C-$tight sequences; see \cite[Corollary VI.3.33]{jacod2003limit}.
This concludes the proof.
\end{proof}

\begin{remark}\label{rem:weakLimit-qlc}
In the proof of the previous lemma, in order to prove the tightness of the sequence $\Phi^k$ we have used \cite[Theorem VI.4.13]{jacod2003limit}, which in turn makes use of \emph{Aldous's criterion for tightness}, see \cite[Section VI.4a]{jacod2003limit}.
This allows us to conclude that every weak--limit point of $(N^k)_{k\in\mathbb N}$ is quasi--left--continuous for its natural filtration.
Recall that by Condition \ref{MFilqlc}, we have in particular that $X^\infty$ and $Y^\infty$ are quasi--left--continuous for their natural filtrations. 
To sum up, for the arbitrary  weak--limit point $\overline N$ of $(N^k)_{k\in\mathbb N}$, it holds 
$$\mathbb P(\Delta X^{\infty,i}_t= 0)=\mathbb P(\Delta Y^\infty_t= 0)=\mathbb P(\Delta \overline N_t= 0)=1,\text{ for every }t\in\mathbb R_+, i=1,\dots,\ell.$$
Observe that the above property is independent of the filtration with respect to which the processes are adapted to.
\end{remark}
%
%
%
%

\vspace{0.5em}
In view of Lemma \ref{lem:NisTight}, we will fix for the rest of this subsection a pair $(\overline N , \Xi)$ such that \eqref{conv:weaklimit} is valid, \emph{i.e.} the
weak--limit point $(\overline N , \Xi)$ is approximated by the subsequence $(N^{k_l},\langle N^{k_l}\rangle)_{l\in\mathbb N}$.
Consequently the subsequence $(k_l)_{l\in\mathbb N}$ will be fixed hereinafter.

\begin{lemma}\label{cor:Nkl-PUT}
The subsequence $(N^{k_l})_{l\in\mathbb N}$ possesses the $\PUT$ property. 
\end{lemma}
\begin{proof}
This is immediate by \cite[Corollary VI.6.30]{jacod2003limit}.
To verify the assumptions of the aforementioned corollary, in view of convergence \eqref{conv:weaklimit}, it is only necessary to prove that
$$\sup_{l\in\mathbb N} \mathbb E\bigg[ \sup_{s\le t} \big|\Delta N^{k_l}_s\big|\bigg]<\infty.$$
This follows by analogous arguments to the ones used in the proof of \cref{lem:R-2I-PUT}.
\end{proof}

\begin{corollary}\label{cor:Conv-Theta}
Assume that conditions {\rm \ref{MFilqlc}, \ref{MSI}, \ref{MFilweak}, \ref{Mfinalrv}} and convergence \eqref{conv:weaklimit} hold. 
Then, for every $I\in\mathcal J(X^\infty),$ $G\in\mathcal G_\infty^\infty,$ we have that 
$\Theta^{k_l,I,G}\xrightarrow[l\to\infty]{\hspace{0.2cm}\mathcal L\hspace{0.2cm}} \Theta^{\infty,I,G}$ 
and that the process $\Theta^{\infty,I,G}$ is a uniformly integrable $\mathbb F^{M^{\infty,G}}-$martingale, where $M^{\infty,G}$ is defined in \eqref{def:M-infG}.
\end{corollary}
\begin{proof}
Let us fix an $I\in\mathcal J(X^\infty)$ and a $G\in\mathcal G^\infty_\infty.$
By Lemma \ref{lem:Theta-ii}, \emph{i.e.} by the tightness of the sequence $(\Theta^{k,I,G})_{k\in\N}$, we obtain that the sequence $(\Theta^{k_l,I,G})_{l\in\N}$ is tight. 
Therefore, it is sufficient to prove the convergence in law of each element of the subsequence $(\Theta^{k_l,I,G})_{l\in\N}$. 

\begin{enumerate}[label={\bf ($\mathbf{\Theta}$\arabic*)}, itemindent=1em, align=left, leftmargin=0cm, label={\rm (\roman*)}]
\item
By conditions \ref{MSI}, \ref{MFilweak} and Proposition \ref{MeminProp1}, the following convergence holds 
						$$X^k\sipmulti{\ell} X^\infty.$$
\item
By the above convergence, conditions \ref{MFilqlc}, \ref{MSI}, \ref{MFilweak} and Theorem \ref{MeminCorollary}, the following convergence holds
						$$[X^k]-\langle X^k\rangle\sipmulti{\ell\times\ell} [X^\infty]-\langle X^\infty\rangle^{\mathbb G^\infty}.$$ 
\item
By conditions \ref{MSI}, \ref{MFilweak}, \ref{Mfinalrv} and Proposition \ref{MeminProp1}, the following convergence holds 
						$$Y^k\sipmulti{}Y^\infty.$$
\item
By Lemma \ref{lem:Xk-Yk-PUT}, the sequences $(X^{k_l})_{l\in\mathbb N}$ and $(Y^{k_l})_{l\in\mathbb N}$ possess the $\PUT$ property.
Moreover, by Lemma \ref{cor:Nkl-PUT} the sequence $(N^{k_l})_{l\in\N}$ has the $\PUT$ property and by Lemma \ref{lem:R-2I-PUT} the sequence $({\rm R}_{2,I}\star\widetilde{\mu}^{(X^{k_l,d},\mathbb G^{k_l})})_{l\in\mathbb N}$ possesses the $\PUT$ property.
Therefore we can obtain by \cite[Theorem VI.6.26]{jacod2003limit} that
\begin{align*}
 [X^{k_l}, N^{k_l}]&\xrightarrow[l\to\infty]{\hspace{0.2cm}\mathcal L\hspace{0.2cm}} [X^\infty, \overline N],\\
[{\rm R}_{2,I}\star\widetilde{\mu}^{(X^{k_l,d},\mathbb G^{k_l})}, N^{k_l}] &\xrightarrow[l\to\infty]{\hspace{0.2cm}\mathcal L\hspace{0.2cm}} [{\rm R}_{2,I}\star\widetilde{\mu}^{(X^{\infty,d},\mathbb G^\infty)},\overline N],\\
 [Y^{k_l}, N^{k_l}] - \langle N^{k_l}\rangle &\xrightarrow[l\to\infty]{\hspace{0.2cm}\mathcal L\hspace{0.2cm}} [Y^\infty, \overline N]- \Xi.
\end{align*}
\item 
By condition \ref{MFilweak}, the following convergence holds
						$$ \E[\mathds{1}_G|\mathcal G^k_\cdot] \sipmulti{} \E[\mathds{1}_G|\mathcal G^\infty_\cdot].$$
\end{enumerate}
In order to show the second statement, we use Corollary \ref{cormultiskorokhod} and that $(X^k)_{k\in\mathbb N}$ is a common element, hence we obtain the convergence 
$$\big(\Theta^{k_l,I,G},M^{k_l,G}\big)
	\xrightarrow[l\to\infty]{\hspace{0.2cm}\mathcal L\hspace{0.2cm}}
	\big(\Theta^{\infty,I,G},M^{\infty,G}\big).$$
Now we can conclude that $\Theta^{\infty,I,G}$ is a uniformly integrable martingale with respect to the filtration generated by $M^{\infty,G}$, which obviously coincides with the filtration generated by $(\Theta^{\infty,I,G},M^{\infty,G})$, by \cref{lem:Theta-i} and \cite[Theorem IX.1.12, Proposition IX.1.10]{jacod2003limit}.
Finally, in order to obtain the martingale property with respect to $\mathbb F^{M^{\infty,G}}$, \emph{i.e.} the usual augmentation of the natural filtration of  $M^{\infty,G},$ we apply \cite[Theorem 2.46]{he1992semimartingale}.
\end{proof}

\begin{remark}\label{rem:ImmFil}
\begin{enumerate}[label={\rm (\roman*)}, itemindent=0.5em, align=left, leftmargin=0cm]
\item\label{rem:ImmFil.i} 
For every $t\in[0,\infty]$, we have $\F_t=\mathcal G_t^\infty\vee\mathcal F^{(\overline{N}, \Xi)^\top}_t.$
\item\label{rem:ImmFil.ii}
The inclusion $\mathbb F^{M^{\infty,G}}\subset \mathbb F$ holds for every fixed $G\in\mathcal G^\infty_\infty$. 
In particular, it also holds that $\mathcal P^{\mathbb F^{M^{\infty,G} } }\subset \mathcal P^{\mathbb F}$, for every fixed $G\in\mathcal G_\infty^\infty$.
%
%
\end{enumerate}
\end{remark}

\begin{proposition}\label{prop:Fmart}
The process $\Theta^{\infty,I}$ is a uniformly integrable $\mathbb F-$martingale, for every $I\in\mathcal J(X^\infty)$.
\end{proposition}

\begin{proof}
Let us fix an $I\in\mathcal J(X^\infty)$.
By applying Corollary \ref{cor:Conv-Theta} for every $G\in\mathcal G_\infty^\infty$ we have that $\Theta^{\infty,I}$ is a uniformly integrable $\mathbb F^{M^{\infty,G}}-$martingale. 
Therefore we have only to prove that $\Theta^{\infty,I}$ is an $\mathbb F-$martingale.
\iftoggle{full}{By Lemma \ref{lem:pi-lambda-signed}, }{By \citet[Lemma~7.2.1]{bogachev2007measure} or \citet[Lemma~A.7]{saplaouras2017backward}, }
it is sufficient to prove that for every $0\leq t<u\leq \infty$ the following condition holds
\begin{align}\label{Equal:onLambda}
\int_{\Lambda} \E[\Theta^{\infty,I}_u|\F_t]\,\dP= \int_{\Lambda} \Theta^{\infty,I}_t\,\dP,
\end{align}
for every
$\Lambda\in\mathcal{A}_t:=
	\big\{\Gamma\cap\Delta,\; \Gamma\in\mathcal G_t^\infty, \Delta\in\mathcal F^{(\overline{N}, \Xi)^\top}_t \big\}.$
Observe that $\mathcal A_t$ is a $\pi-$system\footnote{A $\pi-$system is a non--empty family of sets which is closed under finite intersections.} with $\sigma (\mathcal A_t)=\mathcal F_t.$

\vspace{0.5em}
Let us, therefore, fix $0\leq t<u\leq \infty$ and $\Lambda\in\mathcal{A}_t$, where $\Lambda=\Gamma\cap\Delta$, for some $\Gamma\in\mathcal G_t^\infty,$ $\Delta\in\mathcal F^{(\overline{N}, \Xi)^\top}_t$.
Observe that in particular $\Lambda \in\mathcal F_t^{M^{\infty, \Gamma}}$.
Now we obtain
\begin{align}\label{Property:Fmart}
\int_{\Lambda}\E\big[\Theta^{\infty,I}_u|\F_t\big]\,\dP
	=\int_{\Lambda}\E\big[\E [\Theta^{\infty,I}_u|\F_t]\big|\mathcal F_t^{M^{\infty, \Gamma}}\big]\,\dP
	=\int_{\Lambda}\E \big[\Theta^{\infty,I}_u\big|\mathcal F_t^{M^{\infty, \Gamma}}\big]\,\dP
	=\int_{\Lambda} \Theta^{\infty,I}_t\,\dP,
\end{align}
where the first equality holds because $\Lambda\in\mathcal F_t^{M^{\infty, \Gamma}}$, so that we can use the definition of the conditional expectation with respect to the $\sigma-$algebra $\mathcal F_t^{M^{\infty, \Gamma}}$.
In the second equality, we have used the tower property and Remark \ref{rem:ImmFil}.\ref{rem:ImmFil.ii}, while for the third one we have used that $\Theta^{\infty,I}$ is a $\mathbb F^{M^{\infty,\Gamma}}-$martingale, \textit{i.e.} we applied Corollary \ref{cor:Conv-Theta} for $\Gamma\in\mathcal G_t^\infty \subset \mathcal G^\infty_\infty.$
Therefore, we can conclude that $\Theta^{\infty,I}$ is a uniformly integrable $\mathbb F-$martingale.
\end{proof}

\subsubsection{\texorpdfstring{$\overline N$ is sufficiently integrable}{N is sufficiently integrable}}

The uniform integrability of the $\mathbb F-$martingale $\overline N$ implies neither the integrability of $[\,\overline N\,]^{1/2}_\infty$ nor of $\sup_{t\in[0,\infty)}|\overline N_t|$.
In Lemma \ref{lem:Integrability-barN} we will prove that there exists a suitable function $\Psi$ such that $$\Psi\bigg(\sup_{t\in[0,\infty)]}\big|\overline N_t\big|\bigg)\;  \text{and}\; \Psi\big([\,\overline N\,]^{1/2}_\infty\big)\; \text{are integrable.}$$
This result is crucial in order to show that $M_{\mu^{X^\infty}}\big[|\Delta \overline N|\,\big|\widetilde{\mathcal P}^{\mathbb F}\big]$ is well--defined, see Corollary \ref{cor:Pred-sigma-finite} and Proposition \ref{prop:barN-orthog}.
The way we choose $\Psi$ is given in Definition \ref{def:Psi-Young}.

\vspace{.5em}
We have shown in \cref{lem:Optional-UI} that the sequence $\big({\rm Tr}\big[ [X^k]_\infty \big]\big)_{k\in\overline{\mathbb N}}$ is uniformly integrable.
The orthogonal decomposition of $Y^\infty$ with respect to $(X^{\infty,c},\mu^{X^{\infty,d}},\mathbb G^\infty)$ implies then that the (finite) family 
$\{\langle Z^\infty\cdot X^{\infty,c}\rangle_\infty, [U^\infty\star\widetilde{\mu}^{(X^{\infty,d},\mathbb G^\infty)}]_\infty\}$ is uniformly integrable.
Therefore, the family 
\begin{align}\label{def:Mcal-UI}
\mathcal M:=\big\{ {\rm Tr}\big[ [X^k]_\infty \big], k\in\overline{\mathbb N}\big\} \cup \big\{ \langle Z^\infty\cdot X^{\infty,c}\rangle_\infty, [U^\infty\star\widetilde{\mu}^{(X^{\infty,d},\mathbb G^\infty)}]_\infty \big\},
\end{align}
is uniformly integrable as a finite union of uniformly integrable families of random variables.
By de La Vall\'ee--Poussin's criterion, see \cite[Corollary 2.5.5]{cohen2015stochastic}, we can construct a Young function $\Phi$, see Appendix \ref{Young_functions} for the respective definition, such that 
$\sup_{\Gamma\in\mathcal M} \mathbb E[ \Phi(\Gamma)]<\infty.$
We can improve the last condition by choosing a moderate Young function $\Phi^{\mathcal A}$, see Appendix \ref{Young_functions} for the definition, such that 
\begin{align}\label{moderate-Tr-Xk}
\sup_{\Gamma\in\mathcal M} \mathbb E[ \Phi^{\mathcal A}(\Gamma)]<\infty,
\end{align}
where $\mathcal A:=(\alpha_m)_{m\in\mathbb N}$ is a sequence of non--negative integers such that $\alpha_0=0,$ $\alpha_{2m}\le 2\alpha_m$ and $\lim_{m\to\infty}\alpha_m=\infty,$ see the first Lemma in 
\cite{meyer1978sur} .
\begin{definition}\label{def:Psi-Young}
Let $\Phi^{\mathcal A}$ be a moderate Young function with associated sequence $\mathcal A$, consisting of non--negative integers such that $\alpha_0=0,$ $\alpha_{2m}\le 2\alpha_m$ and $\lim_{m\to\infty}\alpha_m=\infty,$ 
for which condition \eqref{moderate-Tr-Xk} holds.
For $\mathbb R_+\ni x\overset{{\rm quad}}{\longmapsto} 1/2x^2\in\mathbb R_+$, we define the function $\Psi$ to be the Young conjugate of $\Phi^{\mathcal A}\circ{\rm quad}$.
\end{definition}
\begin{remark}
The crucial property of $\Psi$ is that it is moderate, see Proposition \ref{prop:Young_functions}. 
\end{remark}
\begin{lemma}\label{lem:Integrability-barN}
The weak--limit $\overline N$ given by convergence \eqref{conv:weaklimit} and the Young function $\Psi$ from Definition \ref{def:Psi-Young} satisfies
\begin{align*}
\mathbb E\bigg[ \Psi\bigg(\sup_{t\in[0,\infty)} |\overline N_t|\bigg)\bigg]<\infty, 
 \text{ and } 
\mathbb E\Big[ \Psi\Big(\,[\,\overline N\,]_\infty^{1/2}\Big)\Big]<\infty.
\end{align*}
\end{lemma}
\begin{proof}
By convergence \eqref{conv:weaklimit}, Corollary \ref{cor:Conv-Theta}, Remark \ref{rem:weakLimit-qlc} and \cite[Proposition VI.3.14]{jacod2003limit} we have that 
\begin{align*}
N^{k_l}_t		\xrightarrow[l\to\infty]{\hspace{0.2cm}\mathcal L\hspace{0.2cm}} 	\overline N_t,
\text{ for every }t\in\mathbb R_+.
\end{align*}
The function $\Psi$ is continuous and convex, since it admits a representation via a Lebesgue integral with positive and non--decreasing integrand.
By the continuity of $\Psi$ and the above convergence we can also conclude that
\begin{align}\label{conv:Psi-weak}
\Psi\big(|N^{k_l}_t|\big)		\xrightarrow[l\to\infty]{\hspace{0.2cm}\mathcal L\hspace{0.2cm}} 	\Psi(|\overline N_t|),
\text{ for every }t\in\mathbb R_+.
\end{align}

Using now Proposition \ref{prop:Young_functions}.{\rm(iv)}, there exists a Young function $\Upsilon$ and constants $K>0$ and $\overline C>1$ such that 
\begin{align*}
&\sup_{l\in\mathbb N}\mathbb E\bigg[ \Upsilon\bigg(\Psi\bigg( \sup_{t\in[0,\infty)} \big|N^{k_l}_t\big|\bigg)\bigg)\bigg]
=					\sup_{l\in\mathbb N}\mathbb E\bigg[ (\Upsilon\circ\Psi)\bigg( \sup_{t\in[0,\infty)} \big|N^{k_l}_t\big|\bigg)\bigg]\\
&\hspace{1em}\le 	\sup_{l\in\mathbb N}\mathbb E\bigg[ \overline C  \mathds{1}_{[0,\overline C]}\bigg( \sup_{t\in[0,\infty)}\big|N^{k_l}_t\big|\bigg) \bigg] 
+ 	\frac K2 		\sup_{l\in\mathbb N}\mathbb E\bigg[ \sup_{t\in[0,\infty)}\big|N^{k_l}_t\big|^2\mathds{1}_{[\overline C,\infty)}\bigg( \sup_{t\in[0,\infty)}\big|N^{k_l}_t\big|\bigg) \bigg]\\
&\hspace{1em}\le \overline C + \frac K2 
					\sup_{l\in\mathbb N}\mathbb E\bigg[ \sup_{t\in[0,\infty)} |N^{k_l}_t|^2\bigg]\overset{\text{Lem. \ref{lem:NisTight}.\ref{lem:N-Tight-iv}}}{<}\infty.
\numberthis\label{prop:Ml-UI}
\end{align*}
By the above inequality and de La Vall\'ee Poussin's criterion, we obtain the uniform integrability of the family $\big(\Psi(|N^{k_l}_t|)\big)_{l\in\mathbb N, t\in[0,\infty)}$.
On the other hand, convergence \eqref{conv:Psi-weak} and the Dunford--Pettis compactness criterion, see \citet[Chapter II, Theorem 25]{dellacherie1978probabilities}, yield that the set 
\begin{align*}
\mathcal Q := \big\{ \Psi(|N^{k_l}_t|), t\in[0,\infty), l\in\mathbb N\}\cup\{\Psi(|\overline N_t|), t\in[0,\infty) \big\},
\end{align*}
is uniformly integrable, since we augment the relatively weakly compact set $\big(\Psi(|N^{k_l}_t|)\big)_{t\in[0,\infty),l\in\mathbb N}$ merely by aggregating it with the weak--limits $\Psi(|\overline N_t|),$ for $t\in[0,\infty)$.
In particular, the subset $\overline{\mathcal N}_\Psi:=\big( \Psi(|\overline N_t|)\big)_{t\in[0,\infty)}$ is uniformly integrable.
The $\mathbb L^1-$boundedness of $\overline{\mathcal N}_\Psi$ and the $\mathbb F-$martingale property of $\overline N$, see Proposition \ref{prop:Fmart}, imply that the random variable $\overline N_\infty:=\lim_{t\to\infty}\overline N_t$ exists $\mathbb P-a.s.$
Using the uniform integrability of $\overline{\mathcal N}_\Psi$ and the continuity of $\Psi$ once again, we have that 
\begin{align}
\Psi\big(\big|\overline N_t\big|\big)
\xrightarrow[t\to\infty]{\hspace{0.2cm}\mathbb L^1(\Omega,\mathcal G,\mathbb P)\hspace{0.2cm}}
\Psi\big(\big|\overline N_\infty\big|\big),
\end{align}  
\emph{i.e.} $\Psi(|\overline N_\infty|)\in\mathbb L^1(\Omega, \mathcal G, \mathbb P).$
Recall now that the function $\Psi$ is moderate and convex, see Proposition \ref{prop:Young_functions}. 
By the integrability of $\Psi(|\overline N_\infty|)$ we have that $\Vert \overline N_\infty\Vert_\Psi<\infty,$ where $\Vert \Theta\Vert_\Psi:=\inf\big\{\lambda >0, \mathbb E\big[\Psi\big(\frac{|\Theta|}{\lambda}\big)\big]\le 1\big\}$ is the norm of the Orlicz space associated to the Young function $\Psi,$ see \citet[Paragraph 97]{dellacherie1982probabilities}.
Now we are ready to apply Doob's inequality in the form \cite[Inequality 103.1]{dellacherie1982probabilities}, since the Young conjugate of $\Psi$ is the moderate Young function $\Phi^{\mathcal A}\circ{\rm quad}$, with associated constant 
$\overline c_{\Phi^{\mathcal A}\circ{\rm quad}}<\infty$, see Definition \ref{def:Young-constants} for the associated constants of a Young function.
The above yields
\begin{align}\label{Ineq:Doob}
C_\Psi:=\bigg\| \sup_{t\in[0,\infty)}\big|\overline N_t\big|\bigg\|_\Psi \le \overline c_{\Phi^{\mathcal A}\circ{\rm quad}} \Vert \overline N_\infty\Vert_\Psi<\infty.
\end{align}
Inequality \eqref{Ineq:Doob} yields therefore the finiteness of $\Vert \sup_{t\in[0,\infty)}|\overline N_t|\Vert_\Psi$,
which in conjunction with the fact that $\Psi$ is moderate, \emph{i.e.} $\overline c_\Psi<\infty,$ provides also the finiteness of $\mathbb E[\Psi(\sup_{t\in[0,\infty)}|\overline N_t|)].$
The latter can be easily concluded by \citet[Theorem 3.1.1.(b),(d)]{long1993martingale} and the fact that
\begin{align*}
\mathbb E\bigg[ \Psi\bigg(\sup_{t\in[0,\infty)}|\overline N_t|\bigg)\bigg]  
&\le \mathbb E\bigg[ \Psi\bigg(\frac{\sup_{t\in[0,\infty)}|\overline N_t|}{C_\Psi}\bigg)\bigg]\mathds{1}_{\{C_\Psi\le 1\}}+\overline{c}_\Psi^{C_\Psi } \mathbb E\bigg[ \Psi\bigg(\frac{\sup_{t\in[0,\infty)}|\overline N_t|}{C_\Psi}\bigg)\bigg]\mathds{1}_{\{C_\Psi> 1\}},
\end{align*} 
which is finite in any case.
Now, we use the finiteness of $\mathbb E\big[ \Psi(\sup_{t\in[0,\infty)}|\overline N_t|)\big] $, Burkholder--Davis--Gundy's inequality \cite[Theorem 10.36]{he1992semimartingale}, and the fact that $\Psi$ is moderate to conclude that $\mathbb E\big[ \Psi\big(\,[\,\overline N\,]_\infty^{1/2}\big)\big]<\infty$. \qedhere
\end{proof}

\begin{corollary}\label{cor:Pred-sigma-finite}
The weak--limit $\overline N$ in convergence \eqref{conv:weaklimit} satisfies 
\begin{align*}
\mathbb E\bigg[ \int_{(0,\infty)} |\Delta \overline N_s|\, \Vert x\Vert_1\ \mu^{X^{\infty,d}}(\ud s,\ud x)\bigg]<\infty.
\end{align*}
\end{corollary}
\begin{proof}
We have that
\begin{multline*}
\mathbb E\bigg[\int_{(0,\infty)\times \mathbb R^\ell}\big|\Delta \overline N_s\big|\, \Vert x\Vert_1\,\mu^{X^{\infty,d}}(\ud s,\ud x) \bigg]
	 = \sum_{i=1}^\ell \mathbb E \bigg[ \int_{(0,\infty)} \big|\Delta \overline N_s\big|\, |x^i|\mu^{X^{\infty,d}}(\ud s,\ud x)\Big] \\
	 = \sum_{i=1}^\ell \mathbb E\bigg[\sum_{s>0} \big|\Delta \overline N_s\big|\, \big|\Delta X^{{\infty,d},i}_s\big| \Big]
	 \le \sum_{i=1}^\ell \mathbb E\bigg[\bigg(\sum_{s>0} (\Delta \overline N_s)^2\bigg)^{1/2}\,\bigg(\sum_{s>0} (\Delta X^{{\infty,d},i}_s)^2\bigg)^{1/2} \bigg]\\
	 \le \sum_{i=1}^\ell \mathbb E\Big[\, [\,\overline N\,]^{1/2}_\infty\,[X^{\infty,i}]^{1/2}_\infty \Big] \ \ \ \
	 \hspace{-1.3em}\overset{\text{Young Ineq.}}{\underset{\text{Lem. \ref{lem:Integrability-barN}}}{\le}} 
			\sum_{i=1}^\ell \mathbb E\Big[\Psi\big([\,\overline N\,]^\frac12_\infty\big)	+	\Phi^{\mathcal A}\circ{\rm quad}\big([X^{\infty,i}]^{\frac12}_\infty\big) \Big]\\
	 \hspace{-0.3em}\overset{\eqref{moderate-Tr-Xk}}{\le} 	
			 \mathbb E\Big[\Psi\big([\,\overline N\,]^\frac12_\infty\big)	+	\Phi^{\mathcal A}\big({\rm Tr}\big[[X^\infty]_\infty\big]\big) \Big]
	<\infty,
\end{multline*}
where in the last inequality we used also the convexity of $\Phi^\mathcal A$ in order to take out the coefficient $\frac12$, which appears due to the definition of the function ${\rm quad}.$ 
\end{proof}
We conclude this sub--sub--section with the following result, which yields that $(X^{\infty,1} \overline N,\dots,X^{\infty,\ell} \overline N)^\top$ is an uniformly integrable martingale.
\begin{proposition}\label{prop:barN-orthog}
The weak--limit $\overline N$ in convergence \eqref{conv:weaklimit} satisfies
\begin{align}\label{property:PQCareZero}
\langle X^{\infty,c,i}, \overline N^{c} \rangle^\Fil = 0, \text{ for every } i=1,\dots\ell, \text{ and }
M_{\mu^{X^{\infty,d}}}[\Delta \overline N | \widetilde{\mathcal P}^\Fil]= 0.
\end{align} 
\end{proposition}
\begin{proof}
We will apply Proposition \ref{prop:SCO} to the pair of $\mathbb F-$martingales $(X^\infty, \overline N)$.
Firstly, recall Remark \ref{rem:weakLimit-qlc}, \emph{i.e.} that $X^\infty$ and $\overline N$ are $\mathbb F-$quasi--left--continuous $\mathbb F-$martingales.
Now we verify that the aforementioned pair indeed satisfies the requirements of \cref{prop:SCO}.\ref{prop:SCO-i}--\ref{prop:SCO-iii}.
\begin{enumerate}[label=$\bullet$, itemindent=0.0cm, leftmargin=*]
\item By \cref{prop:Fmart} we have that $[X^\infty, \overline N]$ is a uniformly integrable $\mathbb F-$martingale, hence the first condition is satisfied.
\item By the same proposition, we also have that $[{\rm R}_{2,I}\star \widetilde{\mu}^{(X^{\infty,d},\mathbb G^\infty)}, \overline N]$ is a uniformly integrable $\mathbb F-$martingale, for every $I\in\mathcal J(X^\infty).$
Moreover, \iftoggle{full}{by Lemma \ref{IXGeneratesBorel} we have that }{the density of $\mathcal J(X^\infty)$ (see also Lemma A.9 in the extended version \cite{papapantoleon2018stability}) allows one to obtain}
$\sigma\big(\mathcal J(X^\infty)\big) = \mathcal B(\mathbb R^\ell\setminus\{0\}).$
Hence, the second condition is also satisfied.
\item Finally, the property $|\Delta \overline N|\mu^{X^{\infty,d}}\in\widetilde{\mathcal A}_{\sigma}(\mathbb F)$ is equivalent to 
\begin{align}\label{Integr-Cond:V}
\mathbb E\bigg[\int_{(0,\infty)\times \mathbb R^\ell}\big|\Delta \overline N_s\big| V(s,x)\,\mu^{X^{\infty,d}}(\ud s,\ud x) \bigg]<\infty,
\end{align}
for some strictly positive $\mathbb F-$predictable function $V.$
In view of Corollary \ref{cor:Pred-sigma-finite}, \eqref{Integr-Cond:V} holds for $V(\omega, t,x) = \sum_{i=1}^\ell |\pi^i(x)| = \Vert x\Vert_1,$ which is $\mu^{X^{\infty,d}}-a.e.$ strictly positive and $\mathbb F-$predictable as deterministic.
Therefore, \eqref{Integr-Cond:V} is valid, \emph{i.e.} the third condition is satisfied. \qedhere
\end{enumerate}
\end{proof}

\subsubsection{\texorpdfstring{The filtration $\mathbf{\mathbb G^\infty}$ is immersed in the filtration $\mathbb F$}{The filtration G-infty is immersed in the filtration F}}

In this sub--sub--section we use the notation and framework of subsection \ref{sub-sub-sec:StepB-Prep}.
Recall that the filtration $\mathbb F$ has been defined in \eqref{def:FilF}, and that the subsequence $(k_l)_{l\in\mathbb N}$ is fixed and such that convergence \eqref{conv:weaklimit} holds.

\begin{definition}
The filtration $\mathbb G$ is \emph{immersed} in the filtration $\mathbb F$\footnote{A summary of other terms describing the same property can be found in \citet{tsirelson1998within}.} if
\begin{align*}
\mathcal H^2(\mathbb G,\infty;\mathbb R)\subset \mathcal H^2(\mathbb F,\infty;\mathbb R).
\end{align*}
\end{definition}

\begin{lemma}\label{lem:PD-Mart-Indist}
Let $X_0 + X^{\infty,c,\mathbb F} + {\rm Id}\star \widetilde{\mu}^{(X^{\infty,d},\mathbb F)}$ be the canonical representation of $X^\infty$ as an $\mathbb F-$martingale and 
$X_0 + X^{\infty,c} + {\rm Id}\star \widetilde{\mu}^{(X^{\infty,d},\mathbb G^\infty)}$ be the canonical representation of $X^\infty$ as an $\mathbb G^\infty-$martingale.
Then the respective parts are indistinguishable, \emph{i.e.}
\begin{align*}
X^{\infty,c}=X^{\infty,c,\mathbb F}, \text{ and } {\rm Id}\star \widetilde{\mu}^{(X^{\infty,d},\mathbb G^\infty)}={\rm Id}\star \widetilde{\mu}^{(X^{\infty,d},\mathbb F)},
\text{ up to indistinguishability.}
\end{align*}
\end{lemma}

Therefore, we will simply denote the continuous part of the $\mathbb F-$martingale $X^\infty$ by $X^{\infty,c}$ and we will use indifferently ${\rm Id}\star \widetilde{\mu}^{(X^\infty,\mathbb G^\infty)}$ and ${\rm Id}\star \widetilde{\mu}^{(X^\infty,\mathbb F)}$ to denote the discontinuous part of the $\mathbb F-$martingale $X^\infty.$

\begin{proof}
By Proposition \ref{prop:Fmart} the process $X^\infty$ is an $\mathbb F-$martingale whose canonical representation is given by
\begin{align*}
X^\infty  = X_0 + X^{\infty,c,\mathbb F} + {\rm Id}\star \widetilde{\mu}^{(X^{\infty,d},\mathbb F)},
\end{align*}
see  \cite[Corollary II.2.38]{jacod2003limit}.
However $X^\infty$ is $\mathbb G^\infty-$adapted, which in conjunction with \citet[Theorem 2.2]{follmer2011local} implies that the process $X^{\infty,c,\mathbb F} + {\rm Id}\star \widetilde{\mu}^{(X^{\infty,d},\mathbb F)}$ is an $\mathbb G^\infty-$martingale.
On the other hand, the canonical representation of the process $X^\infty$ as a $\mathbb G^\infty-$martingale is 
\begin{align*}
X^\infty  = X_0 + X^{\infty,c} + {\rm Id}\star \widetilde{\mu}^{(X^{\infty,d},\mathbb G^\infty)}.
\end{align*}
Hence, by \cite[Theorem I.4.18]{jacod2003limit} we can conclude the indistinguishability of the respective parts due to the uniqueness of the decomposition. 
\end{proof}
%
%
\begin{lemma}\label{lem:PQC-unaffected}
We have $\langle X^\infty \rangle^{\mathbb G^\infty} = \langle X^\infty \rangle^{\mathbb F}.$
\end{lemma}

Therefore, we will denote the dual predictable projection of $X^\infty$ simply by $\langle X^\infty\rangle.$

\begin{proof}
In Proposition \ref{prop:Fmart} we showed that the process $[X^\infty]-\langle X^\infty \rangle^{\mathbb G^\infty}$ is an $\mathbb F-$martingale. 
On the other hand, $\big({\rm Tr}\big[ [X^\infty]_t\big]\big)_{t\in[0,\infty]}$ is uniformly integrable and in particular of class (D)\footnote{See \cite[Definition I.1.46]{jacod2003limit}.}.
Consequently, by \cite[Theorem I.3.18]{jacod2003limit}, there exists a unique $\mathbb F-$predictable process, say $\langle X^\infty\rangle^{\mathbb F}$, such that $[X^\infty]-\langle X^\infty\rangle^{\mathbb F}$ is a uniformly integrable $\mathbb F-$martingale.
Recall that by definition $\mathbb G^\infty\subset \mathbb F$, therefore $\mathcal{P}^{\mathbb G^\infty}\subset\mathcal{P}^{\mathbb F}$, which allows us to conclude that $\langle X^\infty\rangle^{\mathbb G^\infty} - \langle X^\infty\rangle^{\mathbb F}$ is a uniformly integrable $\mathbb F-$predictable $\mathbb F-$martingale of finite variation.
Therefore, by \cite[Corollary I.3.16]{jacod2003limit}, we obtain that $\langle X^\infty\rangle^{\mathbb G^\infty} - \langle X^\infty\rangle^{\mathbb F}=0$ up to an evanescent set.
\end{proof}
\begin{corollary}\label{cor:Xinf-Fqlc}
The process $X^\infty$ is $\mathbb F-$quasi--left--continuous.
Therefore
\begin{align}\label{property:nu-F}
\nu^{(X^{\infty,d},\mathbb F)}\big(\omega;\{t\}\times\mathbb R^\ell\big)= 0 \ \text{ for every }(\omega,t)\in\Omega\times\mathbb R_+. 
\end{align}
\end{corollary}
\begin{proof}
The $\mathbb F-$quasi--left--continuity of $X^\infty$ is immediate by the continuity of $\langle X^\infty \rangle$ and \cite[Theorem I.4.2]{jacod2003limit}.
By \cite[Corollary II.1.19]{jacod2003limit}, we conclude that \eqref{property:nu-F} holds.
\end{proof}

%
%
%
Now, we can obtain some useful properties about the predictable quadratic covariation of the continuous and the purely discontinuous martingale part of $X^\infty$. 
\begin{lemma}\label{lem:PQcohen2015stochasticqual}
We have $\langle X^{\infty,c}\rangle^{\mathbb F}=\langle X^{\infty,c}\rangle^{\mathbb G^\infty}$ and $\nu^{(X^{\infty,d},\mathbb F)}|_{\mathcal P^{\mathbb G^\infty}}=\nu^{(X^{\infty,d},\mathbb G^\infty)}. $
\end{lemma}

Therefore, we will denote the dual predictable projection of $X^{\infty,c}$ simply by $\langle X^{\infty,c}\rangle.$

\begin{proof}
For the reader's convenience we separate the proof in two parts.
\begin{enumerate}[label={\rm (\roman*)}, itemindent=0.5em, align=left, leftmargin=0cm]
\item First, we prove that $\nu^{(X^{\infty,d},\mathbb F)}|_{\mathcal P^{\mathbb G^\infty}}=\nu^{(X^{\infty,d},\mathbb G^\infty)}$. 
Indeed, recalling \cref{rem:PurDisJum} and that $X^\infty$ is both $\mathbb G^\infty-$ and $\mathbb F-$quasi--left--continuous, it holds for every $t\in[0,\infty)$
\begin{align*}
\int_{\mathbb R^\ell}x \widetilde{\mu}^{(X^{\infty,d},\mathbb G^\infty)}(\{t\}\times\dx)
	\overset{\ref{MFilqlc}}{=} \int_{\mathbb R^\ell}x \mu^{X^\infty}(\{t\}\times\dx) 
 	\overset{\eqref{property:nu-F}}{=} 	\int_{\mathbb R^\ell}x \widetilde{\mu}^{(X^{\infty,d},\mathbb F)}(\{t\}\times\dx), \hspace{0.3em} \mathbb P-a.s.
\numberthis\label{property:IdentJump}
\end{align*}
Consequently, for every non--negative, $\mathbb G^\infty$--predictable function $\theta$, using \cite[Theorem II.1.8]{jacod2003limit}, it holds
\begin{align*}
\E\big[\theta * \nu^{(X^{\infty,d},\mathbb F)}_\infty\big] 
 	\overset{\eqref{property:IdentJump}}{=} \E\bigg[\sum_{s>0}\theta(s,\Delta X^{\infty, d})\mathds{1}_{[\Delta X^{\infty, d}\neq 0]}\bigg] 
 	 = \E\big[\theta * \nu^{(X^{\infty,d},\mathbb G^\infty)}_\infty\big].
\end{align*}
Therefore we can conclude that $\nu^{(X^{\infty,d}, \mathbb F)}\big|_{\mathcal P^{\mathbb G^\infty}} = \nu^{(X^{\infty,d},\mathbb G^\infty)}$.  

\vspace{0.5em}
\item Now we prove that $\langle X^{\infty,c}\rangle^{\mathbb G^\infty} = \langle X^{\infty,c} \rangle^{\mathbb F}.$
We will combine the previous part with Lemma \ref{lem:PQC-unaffected}.
The map {\rm Id} is both an $\mathbb G^\infty-$ and $\mathbb F-$predictable function as deterministic and continuous.
Then 
\begin{align*}
\langle X^{\infty,c} \rangle^\mathbb F 
	&=  \langle X^\infty \rangle - \langle {\rm Id}\star \widetilde{\mu}^{(X^{\infty,d},\mathbb F)}\rangle^{\mathbb F} 
	 \overset{\text{Lem. \ref{lem:PD-Mart-Indist}}}{\underset{\text{Cor. \ref{cor:Xinf-Fqlc}}}{=}} \langle X^{\infty}\rangle - |{\rm Id}|^2 *\nu^{(X^{\infty,d},\mathbb F)} \\
	&\overset{\text{\rm (i)}}{=} \langle X^{\infty}\rangle - |{\rm Id}|^2 *\nu^{(X^{\infty,d},\mathbb G^\infty)}
	 \overset{\ref{MFilqlc}}{=} \langle X^{\infty}\rangle - \langle {\rm Id}\star \widetilde{\mu}^{(X^{\infty,d},\mathbb G^\infty)}\rangle^{\mathbb G^\infty} =  \langle X^{\infty,c} \rangle^{\mathbb G^\infty}. \qedhere
\end{align*} 
\end{enumerate}
\end{proof}

In view of the previous lemmata, we are able to prove that every $\mathbb G^\infty-$stochastic integral with respect to $X^{\infty}$ is also an $\mathbb F-$martingale. 
The exact statement is provided below.

\begin{lemma}\label{lem:Immerse-Cont}
Let $Z\in\mathbb H^2(X^{\infty,c},\mathbb G^\infty,\infty;\mathbb R^\ell)$, then we have $Z\cdot X^{\infty,c}\in\mathcal H^{2,c}(\mathbb F,\infty;\mathbb R)$. 
\end{lemma}
\begin{proof}
By Lemma \ref{lem:PQcohen2015stochasticqual}, we can easily conclude that $\mathbb H^2(X^{\infty,c},\mathbb G^\infty,\infty;\mathbb R^\ell) \subset \mathbb H^2(X^{\infty,c},\mathbb F,\infty;\mathbb R^\ell)$.
We are going to prove the required property initially for simple integrands.
Assume that $\rho, \sigma$ are $\mathbb G^\infty-$stopping times such that $\rho\le \sigma$, $\mathbb P-a.s.$, and that $\psi$ is an $\mathbb R^\ell-$valued, bounded and $\mathcal G_\rho^\infty-$measurable random variable.
Then, see \cite[Theorem 4.5]{jacod2003limit}, the $\mathbb G^\infty-$stochastic integral $(\psi\mathds{1}_{\rrbracket\rho,\sigma\rrbracket})\cdot X^{\infty,c}$
is defined as 
\begin{align}\label{repr:cont}
\psi\mathds{1}_{\rrbracket\rho,\sigma\rrbracket}\cdot X^{\infty,c} = \sum_{i=1}^\ell \int_0^\cdot \psi^i\mathds{1}_{\rrbracket\rho,\sigma\rrbracket} \ud X^{\infty,c,i}.
\end{align}
Treating now $\psi$ as an $\mathcal F_\rho-$measurable variable, since $\mathbb G^\infty\subset \mathbb F$, and using that $X^{\infty,c}$ is an $\mathbb F-$martingale, Proposition \ref{prop:Fmart} yields that the representation in \eqref{repr:cont} is also an $\mathbb F-$martingale.
By \cite[Theorem III.4.5 -- Part a]{jacod2003limit} we can conclude for an arbitrary $Z\in\mathbb H^2(X^{\infty,c},\mathbb G^\infty,\infty;\mathbb R)$, since the process $Z\cdot X^{\infty,c}$ can be approximated by $\mathbb F-$martingales with representation as in \eqref{repr:cont}.
\end{proof}

\begin{lemma}\label{lem:Immerse-PD}
Let $U\in\mathbb H^2(\mu^{X^{\infty,d}},\mathbb G^\infty,\infty;\mathbb R)$, then it holds that $U\star \widetilde{\mu}^{(X^{\infty,d},\mathbb G^\infty)}\in\mathcal H^{2,d}(\mathbb F,\infty;\mathbb R)$.
Moreover, the processes $U\star \widetilde{\mu}^{(X^{\infty,d},\mathbb G^\infty)}$ and $U\star \widetilde{\mu}^{(X^{\infty,d},\mathbb F)}$ are indistinguishable.
\end{lemma}
\begin{proof}
Let $U\in\mathbb H^2(\mu^{X^{\infty,d}},\mathbb G^\infty,\infty;\mathbb R)$.
The inclusion $\mathcal P^{\mathbb G^\infty}\subset \mathcal P^{\mathbb F}$ and the equalities
\begin{align*}
\mathbb E\bigg[ \sum_{s>0} \bigg| \int_{\mathbb R^\ell} U(s,x)\widetilde{\mu}^{(X^{\infty,d},\mathbb F)}(\{s\}\times\dx) \bigg|^2 \bigg]
&	 \overset{\text{Cor. \ref{cor:Xinf-Fqlc}}}{=} 
\mathbb E\bigg[ \sum_{s>0} \bigg| \int_{\mathbb R^\ell} U(s,x)\mu^{X^\infty}(\{s\}\times\dx) \bigg|^2 \bigg]\\
&\hspace{0.8em}	\overset{\ref{MFilqlc}}{=}
\mathbb E\bigg[ \sum_{s>0} \bigg| \int_{\mathbb R^\ell} U(s,x)\widetilde{\mu}^{(X^{\infty,d},\mathbb G^\infty)}(\{s\}\times\dx) \bigg|^2 \bigg]<\infty,
\end{align*}  
yield that $U\in\mathbb H^2(\mu^{X^{\infty,d}},\mathbb F,\infty;\mathbb R)$, hence the $\mathbb F-$martingale $U\star \widetilde{\mu}^{(X^{\infty,d},\mathbb F)}$ is well--defined.

\vspace{0.5em}
Let now $W$ be a positive $\mathbb G^\infty-$predictable function such that $\mathbb E\big[W*\mu^{X^{\infty,d}}_\infty\big]<\infty.$
By \cite[Theorem II.1.8]{jacod2003limit}, we have that $\mathbb E\big[W*\nu^{(X^{\infty,d},\mathbb G^\infty)}_\infty\big]<\infty$, as well as $\mathbb E\big[W*\nu^{(X^{\infty,d},\mathbb F)}_\infty\big]<\infty.$
Then, the property $\nu^{(X^{\infty,d},\mathbb F)}|_{\mathcal P^{\mathbb G^\infty}}=\nu^{(X^{\infty,d},\mathbb G^\infty)}$ from \cref{lem:PQcohen2015stochasticqual} translates into
\begin{align*}
W\star \widetilde{\mu}^{(X^{\infty,d},\mathbb G^\infty)} 
	& =	W* \widetilde{\mu}^{(X^{\infty,d},\mathbb G^\infty)} 
	  =	W* \mu^{X^{\infty,d}}  - W* \nu^{(X^{\infty,d},\mathbb G^\infty)}\\ 
	& \overset{\mathclap{\text{\cref{lem:PQcohen2015stochasticqual}}}}{=}\hspace{1.5em}	
		W* \mu^{X^{\infty,d}}  - W* \nu^{(X^{\infty,d},\mathbb F)}
	  =	W* \widetilde{\mu}^{(X^{\infty,d},\mathbb F)}\\
	&  = W\star \widetilde{\mu}^{(X^{\infty,d},\mathbb G^\infty)} \text{ up to indistinguishability,}
\end{align*}  
where in the first and the last equalities, we used \cite[Proposition II.1.28]{jacod2003limit}, while in the second one as well as in the second to last one, we used the definition of the compensated integer valued measure.
It is immediate that the above equality holds also if $W$ is a real--valued $\mathbb G^\infty-$predictable function such that $\mathbb E\big[|W|*\mu^{X^{\infty,d}}_\infty\big]<\infty,$ \emph{i.e.}
\begin{align}\label{eqlt:FV-Mart}
W\star \widetilde{\mu}^{(X^{\infty,d},\mathbb G^\infty)} =W\star \widetilde{\mu}^{(X^{\infty,d},\mathbb F)},\text{ for $W$ as described above.}
\end{align}
In other words, when the $\mathbb G^\infty-$predictable integrand $W$ is such that $W\star\widetilde{\mu}^{(X^{\infty,d},\mathbb G^\infty)}$ is of finite--variation, then it is indistinguishable from $W\star\widetilde{\mu}^{(X^{\infty,d},\mathbb F)}$.

\vspace{0.5em}
In view of the above discussion, we can conclude that for an arbitrary $U\in\mathbb H^2(\mu^{X^{\infty,d}},\mathbb G^\infty,\infty;\mathbb R)$, we also have $U\star\widetilde{\mu}^{(X^{\infty,d},\mathbb G^\infty)}=U\star\widetilde{\mu}^{(X^{\infty,d},\mathbb F)}$ up to indistinguishability. 
Indeed, let us denote by $(\tau^m)_{m\in\mathbb N}$  the sequence of $\mathbb G^\infty-$totally inaccessible $\mathbb G^\infty-$stopping times 
 which exhausts the thin $\mathbb G^\infty-$optional set $[\Delta X^{\infty,d}\neq 0]$ (see \citep[Proposition I.1.32, Proposition I.2.26]{jacod2003limit}), and fix some $U\in\mathbb H^2(\mu^{X^{\infty,d}},\mathbb G^\infty,\infty;\mathbb R)$.
Then 
\begin{align*}
\Delta(U\star\widetilde{\mu}^{(X^{\infty,d},\mathbb G^\infty)})_{\tau^m} = U(\tau^m,\Delta X^{\infty,d}_{\tau^m})\in\mathbb L^2(\Omega,\mathcal G_{\tau^m}^\infty,\mathbb P;\mathbb R), \text{ for every }m\in\mathbb N.
\end{align*}
By \cite[Theorem~10.2.10]{cohen2015stochastic}, we know that, for every $m\in\mathbb N,$ there exists a continuous and $\mathbb G^\infty-$adapted process $\Pi^{*,m}$ such that
\begin{align*}
M^m:=U\big(\tau^m,\Delta X^{\infty,d}_{\tau^m}\big)\mathds{1}_{\llbracket \tau^m,\infty\llbracket} - \Pi^{*,m}\in\mathcal H^{2,d}\big(\mathbb G^\infty,\infty;\mathbb R\big), \text{ for every }m\in\mathbb N.
\end{align*}

Using \ref{MXWPRP}, for every $m\in\mathbb N$, there exists $U^m\in\mathbb H^2(\mu^{X^{\infty,d}},\mathbb G^\infty,\infty;\mathbb R)$ such that 
\begin{align*}
M^m
	=	U^m\star\widetilde{\mu}^{(X^{\infty,d},\mathbb G^\infty)}
	\overset{\eqref{eqlt:FV-Mart}}{=}	U^m\star\widetilde{\mu}^{(X^{\infty,d},\mathbb F)},
\end{align*}
where the second equality holds because $M^m$ is a process of finite variation for every $m\in\mathbb N$, since it is a single--jump process.
Consequently, $M^m\in\mathcal H^{2,d}(\mathbb F,\infty;\mathbb R)$ for every $m\in\mathbb N.$
Finally, by \cite[Theorem 10.2.14]{cohen2015stochastic} we can approximate (the precise argument is presented in the proof of the aforementioned theorem) both  
$U\star\widetilde{\mu}^{(X^{\infty,d},\mathbb G^\infty)}$ and $U\star\widetilde{\mu}^{(X^{\infty,d},\mathbb F)}$ by the same sequence $\big(\sum_{m=1}^n M^m\big)_{n\in\mathbb N}$.
By the uniqueness of the limit, $U\star\widetilde{\mu}^{(X^{\infty,d},\mathbb G^\infty)}=U\star\widetilde{\mu}^{(X^{\infty,d},\mathbb F)}$ up to indistinguishability.
\end{proof}

Concluding, since every $\mathbb G^\infty$--martingale is also an $\mathbb F$--martingale, the filtration $\mathbb G^\infty$ is immersed in $\mathbb F$.
The following result shows that the orthogonal decomposition remains the same under both filtrations.

\begin{corollary}\label{cor:Y-inf-repr}
The orthogonal decomposition of $Y^\infty$ with respect to $(X^{\infty,c},\mu^{X^{\infty,d}},\mathbb F)$ is given by
\begin{align}\label{repr:Ginf}
Y^\infty 	 = Y_0 + Z^\infty\cdot X^{\infty,c} + U^\infty\star \widetilde{\mu}^{(X^{\infty,d},\mathbb F)},
\end{align} 
where $Z^\infty\in\mathbb H^2(X^{\infty,c},\mathbb G^\infty,\infty;\mathbb R)$ and $U^\infty\in\mathbb H^2(\mu^{X^{\infty,d}},\mathbb G^\infty,\infty;\mathbb R)$ are determined by the orthogonal decomposition of $Y^\infty$ with respect to $(X^{\infty,c},\mu^{X^{\infty,d}},\mathbb G^\infty)$, see \cref{RobMartRep}.
In other words, the orthogonal decomposition of $Y^\infty$ with respect to $(X^{\infty,c},\mu^{X^{\infty,d}},\mathbb F)$ is indistinguishable from 
the orthogonal decomposition of $Y^\infty$ with respect to $(X^{\infty,c},\mu^{X^{\infty,d}},\mathbb G^\infty)$.
\end{corollary}
\begin{proof}
In Proposition \ref{prop:Fmart}, we have proven that $Y^\infty\in\mathcal H^2(\mathbb F,\infty;\mathbb R)$, therefore the orthogonal decomposition of $Y^\infty$ with respect to $(X^{\infty,c},\mu^{X^{\infty,d}},\mathbb F)$ is well--defined.
Assume that
\begin{align*}
Y^\infty = Y_0 + Z^{\infty,\mathbb F}\cdot X^{\infty,c} + U^{\infty,\mathbb F}\star \widetilde{\mu}^{(X^{\infty,d},\mathbb F)} + N^{\infty,\mathbb F},
\end{align*} 
where 	$Z^{\infty,\mathbb F}\in\mathbb H^2(X^{\infty,c},\mathbb F,\infty;\mathbb R)$,
		$U^{\infty,\mathbb F}\in\mathbb H^2(\mu^{X^{\infty,d}},\mathbb F,\infty;\mathbb R)$
	and $N^{\infty,\mathbb F}\in\mathcal H^2(X^\perp, \mathbb F,\infty;\mathbb R).$
On the other hand, by Lemmata \ref{lem:Immerse-Cont} and \ref{lem:Immerse-PD}, we have that $Z^\infty\cdot X^{\infty,c}\in\mathcal H^{2,c}(\mathbb F,\infty;\mathbb R)$ 
 and $U^\infty\star \widetilde{\mu}^{(X^\infty,\mathbb G^\infty)}\in\mathcal H^{2,d}(\mathbb F,\infty;\mathbb R)$, that is to say
 \begin{align*}
 Y^\infty 
 	= Y_0 + Z^\infty\cdot X^{\infty,c} + U^\infty\star\widetilde{\mu}^{(X^{\infty,d},\mathbb G^\infty)}
 	= Y_0 + Z^\infty\cdot X^{\infty,c} + U^\infty\star\widetilde{\mu}^{(X^{\infty,d},\mathbb F)}.
 \end{align*}

Hence, from \cite[Theorem III.4.24]{jacod2003limit}, we get that, up to indistinguishability
\begin{equation*}
Z^{\infty}\cdot X^{\infty,c} = Z^{\infty,\mathbb F}\cdot X^{\infty,c}, \;
	U^\infty\star \widetilde{\mu}^{(X^{\infty,d},\mathbb G^\infty)}= U^{\infty,\mathbb F}\star \widetilde{\mu}^{(X^{\infty,d},\mathbb F)},\;  
	N^{\infty}=0.	\qedhere
\end{equation*}
\end{proof}

\subsection{Proof of the main Theorem}\label{subsec:Proof}

This subsection is devoted to the proof of the main theorem. 
In view of the preparatory results obtained in the previous sections, as well as of the outline of the proof presented in subsection \ref{sec:outline}, the following proof basically amounts to proving that \ref{StepA} is valid.
%
%
\begin{proof}[Proof of Theorem \ref{RobMartRep}]
By Lemma \ref{lem:NisTight}.\ref{lem:N-Tight-iii} we have that the sequence $\big((N^k, \langle N^k\rangle)\big)_{k\in\mathbb N}$ is tight in $\mathbb D(\mathbb R^2)$.
Therefore, an arbitrary subsequence $\big((N^{k_l}, \langle N^{k_l}\rangle)\big)_{l\in\mathbb N}$ has a further subsequence $\big((N^{k_{l_m}}, \langle N^{k_{l_m}}\rangle)\big)_{m\in\mathbb N}$ which converges in law, say to $(\overline N,\Xi)$, \emph{i.e.}
\begin{align}\label{conv:WeakLimit-m}
(N^{k_{l_m}}, \langle N^{k_{l_m}} \rangle)^\top \xrightarrow[m\to\infty]{\hspace{0.2cm}\mathcal{L}\hspace{0.2cm}} (\overline N, \Xi)^\top, 
\end{align}
where  $\overline N$ is a \cadlag process and $\Xi$ is a continuous and increasing process.
The continuity of $\Xi$ follows from Lemma \ref{lem:NisTight}.\ref{lem:N-Tight-i}.
Therefore, we can use the results of subsection \ref{sec:StepB} for the subsequence $\big((N^{k_{l_m}}, \langle N^{k_{l_m}}\rangle)\big)_{m\in\mathbb N}$ and the pair $(\overline N, \Xi)$.

\vspace{0.5em}
By Proposition \ref{prop:barN-orthog} and Corollary \ref{cor:Y-inf-repr} we conclude 
\begin{align*}
\langle Y^\infty,\overline N\rangle 
	& = \langle Z^\infty \cdot X^{\infty,c},\overline N^c\rangle + \langle U^\infty\star \widetilde{\mu}^{(X^{\infty,d},\mathbb G^\infty)}, \overline N^d\rangle\\
	& = \langle Z^\infty \cdot X^{\infty,c},\overline N^c\rangle + \langle U^\infty\star \widetilde{\mu}^{(X^{\infty,d},\mathbb F)}, \overline N^d\rangle\\
	& = Z^\infty \cdot\langle X^{\infty,c},\overline N^c\rangle + \big(U^\infty M_{\mu^{X^{\infty,d}}}[\Delta \overline N|\widetilde{\mathcal P}^{\mathbb F}]\big)*\nu^{(X^{\infty,d},\mathbb F)} 
	  \overset{\eqref{property:PQCareZero}}{=} 0,
	\numberthis\label{YN-Martingale}
\end{align*}
\emph{i.e.} $[Y^\infty,\overline{N}]$ is an $\mathbb F-$martingale.
On the other hand, we have proved in Proposition \ref{prop:Fmart} that $[Y^\infty,\overline{N}]- \Xi\,$ is an $\mathbb F-$martingale as well.
By subtracting the two martingales, we obtain that $\Xi$ is also an $\Fil-$martingale. 
Hence, $\Xi$ is an $\mathbb F-$predictable process of finite variation and a martingale, therefore it has to be constant, see \cite[Corollary I.3.16]{jacod2003limit}. 
Now, we have that
\begin{align*}
\langle N^{k_{l_m}}\rangle \xrightarrow[m\to\infty]{\hspace{0.2cm}\mathcal L\hspace{0.2cm}} \Xi 
\ \text{ implies that } \
\langle N^{k_{l_m}}\rangle_0 \xrightarrow[m\to\infty]{\hspace{0.2cm}\mathcal L\hspace{0.2cm}} \Xi_0.
\end{align*}
Recall that by definition $\langle N^k\rangle_0=0$ for every $k\in\mathbb N,$ hence $\Xi_0=0$. 
Therefore $\Xi=0$ and, since the limit is a deterministic process, the convergence above is equivalent to the following
\[
\langle N^{k_{l_m}}\rangle \xrightarrow[m\to\infty]{\hspace{0.2cm}(\J_1,\mathbb P)\hspace{0.2cm}} 0.
\]


\vspace{0.5em}
Since the limit above is common for every subsequence and $(\mathbb D, \J_1(\mathbb R))$ is Polish, we can conclude from \cite[Theorem 9.2.1]{dudley2002real} that 
$$\langle N^k \rangle \sipmulti{} 0.$$
Using Lemma \ref{lem:NisTight}.\ref{lem:N-Tight-v} and \cite[Theorem 1.11]{he1992semimartingale}, we can strengthen the above convergence to 
\begin{align}\label{conv:PQCNk-L1}
\langle N^k\rangle\silmulti{}{1}0.
\end{align}
Then, we can also conclude that $N^k\silmulti{}{2}0$. 
Indeed, for every $R>0$ and by Doob's $\mathbb{L}^2-$inequality we obtain
\begin{align*}
\E\bigg[ \sup_{t\in[0,R]}|N^k_t|^2\bigg]\leq 4 \E\big[ |N^k_{R}|^2 \big] = 4 \E\big[\langle N^k\rangle_{R}\big]\xrightarrow[\hspace{0.2cm}k\to\infty\hspace{0.2cm}]{}0,
\end{align*}
which implies the convergence $\E[d_{\lu}^2(N^k,0)]\longrightarrow0$. 

\vspace{0.5em}
Using the convergence of $Y^k\silmulti{}{2}Y^\infty$ and the convergence of $(N^k)_{k\in\N}$ to the zero process, which is trivially continuous, we can obtain the joint convergence 
$$(Y^k,N^k)\silmulti{2}{2} (Y^\infty,0).$$
Moreover, using the orthogonal decompositions of $Y^k$ and $Y^\infty$ and the previous results, we obtain 
$$Z^k\cdot X^{k,c}+ U^k\star\mutilde^{(X^{k,d},\mathbb G^k)}=Y^k-N^k - Y^k_0
	\silmulti{}{2}Y^\infty-Y^\infty_0 = Z^\infty\cdot X^{\infty,c}+ U^\infty\star\mutilde^{(X^{\infty,d},\mathbb G^\infty)},$$
which yields then \eqref{RobMartRepi}. Thus, it is only left to prove convergence \eqref{RobMartRepii}. 
Since the sequences $(X^k)_{k\in\overline{\mathbb N}}$ and $(Y^k)_{k\in\overline{\mathbb N}}$ satisfy the conditions of Theorem \ref{MeminCorollary}, we obtain in particular that the sequences 
$(Y^k + X^k)_{k\in\overline{\mathbb N}}$ and $(Y^k - X^k)_{k\in\overline{\mathbb N}}$ also satisfy the conditions of this theorem.
Therefore we can conclude that  
\begin{align*}
\langle Y^k + X^k\rangle \sipmulti{}\langle Y^\infty + X^\infty \rangle,\text{ and }
\langle Y^k - X^k\rangle \sipmulti{}\langle Y^\infty- X^\infty \rangle.
\end{align*}
By the continuity of the limiting processes, recall \ref{MFilqlc} and \cite[Theorem 4.2]{jacod2003limit}, using the identity 
\begin{align*}
\langle Y^k,X^k\rangle = \frac{1}{4}\big(\langle Y^k + X^k\rangle - \langle Y^k - X^k\rangle \big), \text{ for every }k\in\overline{\mathbb N},
\end{align*}
and convergence \eqref{conv:PQCNk-L1}, we can conclude that
\begin{align*}
(\langle Y^k\rangle, \langle Y^k , X^k\rangle, \langle N^k \rangle )
\xrightarrow{ \hspace{.3cm} \left(\J_1(\R{}\times\R\ell\times\R{}),\mathbb L^1\right) \hspace{.3cm} }
(\langle Y^\infty \rangle, \langle Y^\infty, X^\infty \rangle, 0).
\end{align*}
In order to strengthen the last convergence to an $\mathbb L^1-$convergence, we only need to recall that by Theorem \ref{MeminCorollary} the sequences $\big({\rm Tr}\big[\langle X^k\rangle_\infty\big]\big)_{k\in\overline{\mathbb N}}$ and $\big(\langle Y^k\rangle_\infty\big)_{k\in\overline{\mathbb N}}$  are uniformly integrable. 
Now Lemma \ref{UIplusL2Bounded} provides the uniform integrability of $\big(\big\Vert {\rm Var}(\langle Y^k, X^k\rangle)_\infty\big\Vert_1\big)_{k\in\overline{\mathbb N}}$, which allows us to conclude.
\end{proof}

\appendix

\section{Auxiliary results}\label{AppendixAux}

\subsection{\texorpdfstring{Joint convergence for the Skorokhod $\J_1$--topology}{Joint convergence for the Skorokhod J 1-topology}}\label{subsec:JointConv}

This appendix contains some useful results about the joint convergence of sequences in the Skorokhod $\J_1-$topology, which are heavily used in Sub--section \ref{sec:skor} and throughout \cref{sec:RobMartRepSection}.
Let us recall that the spaces $\D^{q}=\D ([0,\infty);\R{q})$ and $\D ([0,\infty);\R{})^{q}$ 
do not coincide topologically, see \cite[Statements VI.1.21-22]{jacod2003limit}. 
Proposition VI.2.2 in \cite{jacod2003limit} describes the relationship between the convergence on the Skorokhod space of the product state space and the convergence on the product of the Skorokhod spaces of one dimensional state spaces; see also \citet[Lemma~3.5]{aldous1981weak} and \citet[Proposition~3.6.5]{ethier1986markov}.
A suitable variation of \cite[Proposition~VI.2.2]{jacod2003limit} is provided in \cite[Lemma~1]{coquet2001weak}, which we state here for convenience.

\begin{lemma}\label{JointSkorokhodConv}
Let $\alpha^n:=(\alpha^{n,1},\ldots,\alpha^{n,q})^{\transp}\in\D^{q}$, for every $n\in\overline{\N}$.
The convergence $\alpha^n\Jconv{q}\alpha^{\infty}$ holds if and only if the following hold
\[\alpha^{n,i} \Jconv{} \alpha^{\infty,i},\; \text{for $i=1,\ldots,q$, and}\;  \sum_{i=1}^p \alpha^{n,i}\Jconv{} \sum_{i=1}^p \alpha^{\infty,i},\; \text{ for $p=1,\ldots,q$.}\]
\end{lemma}
The following corollary is also useful for our purposes.
\begin{corollary}\label{cormultiskorokhod}
Let $\alpha^n, \beta^n,\gamma^n\in\D,$ for every $n\in\overline{\N}.$
If 
$$(\alpha^n,\beta^n)^{\transp}\Jconv{2}(\alpha^{\infty},\beta^{\infty})^{\transp}\text{ and }(\alpha^n,\gamma^n)^{\transp}\Jconv{2}(\alpha^{\infty},\gamma^{\infty})^{\transp}$$ 
then
	$$(\alpha^n, \beta^n, \gamma^n)^{\transp}\Jconv{3}(\alpha^{\infty}, \beta^{\infty}, \gamma^{\infty})^{\transp}.$$
\end{corollary}
\begin{proof}
Consider the sequence $\big((\alpha^n, \beta^n, -\beta^n, \gamma^n)^{\transp}\big)_{n\in\overline{\N}}$. Using that the pairs converge, conditions (i)--(ii) of Lemma \ref{JointSkorokhodConv} are satisfied, the desired result follows since
\[
(\alpha^n, \beta^n, -\beta^n, \gamma^n)^{\transp}\Jconv{4}(\alpha^{\infty},\beta^{\infty},-\beta^{\infty},\gamma^{\infty})^{\transp}.
\]
\end{proof}
\begin{remark}\label{rem:cormultiskorokhod}
The result above does not depend on the dimension of the state spaces, and can be generalized inductively to an arbitrary number of sequences, as long as a common converging sequence exists.
\end{remark}
%
%
%
The following lemma is another convenient tool when we want to conclude joint convergence.

\begin{lemma}\label{lem:ToolforJoint} 
Let $\alpha^n:=(\alpha^{n,1},\ldots,\alpha^{n,q})^{\transp}\in\D^{q}$ for every $n\in\overline{\N}$, and $f:\mathbb R^q\longrightarrow \mathbb R^p$ be a function continuous on $\mathbb R^p$.
If $\alpha^n\Jconv{q} \alpha^\infty$ then $f(\alpha^n)\Jconv{p}f(\alpha^\infty).$ 
\end{lemma}
\begin{proof}
See \cite[Lemma~2.8]{aldous1981weak}.
\end{proof}
%
\subsubsection{\texorpdfstring{$\J_1$--continuous functions}{J 1-continuous functions}}
Let 
\begin{align}\label{def:mathcalI}
\mathcal I:=\big\{ I \subset \mathbb R,\ I\text{ is a subinterval of }(-\infty,0) \text{ or } I\text{ is a subinterval of }(0,\infty)\big\}.
\end{align}
The aim of this sub--sub--section is the following:
given a sequence $\alpha^n\inJmulti{q} \alpha^\infty$ and a function $g:\mathbb R^q\longrightarrow \mathbb R$, we want to define a sequence $(\zeta^n)_{n\in\mathbb N}\subset\D$, where $\zeta^n$ will be constructed using $\alpha^n$ for all $n\in\mathbb N$, such that
\begin{align}\label{conv:toprove}
\zeta^n\inJmulti{} \sum\limits_{0<t\le \cdot} g(\Delta \alpha^\infty_t)\mathds{1}_I(\Delta\alpha^\infty_t), \text{ for suitable } I\in(\mathcal I\cup\{\mathbb R\})^q.
\end{align}
The above convergence will be presented in Proposition \ref{prop:ContFuncInSkor}. 
To this end we introduce the necessary notation and liminary results.
For $\beta\in\D$
we introduce the sets
\begin{align*}
W(\beta)&:= \{u\in\mathbb{R}\setminus\{0\}, \exists t>0 \text{ with }\Delta\beta_t=u \},\;
\mathcal{I}(\beta):= \{(v,w)\subset\mathbb{R}, v<w,\; vw>0 \text{ and } v,w\notin W(\beta)  \}.
\end{align*}
The set $W(\beta)$ collects the heights of the jumps of the function $\beta.$
The set $\mathcal{I}(\beta)$ collects all the open intervals of $\mathbb{R}\setminus\{0\}$ with boundary points of the same sign, which, moreover, do \emph{not} belong to $W(\beta).$
Observe that $\mathcal I(\beta)\subset \mathcal I,$ for every $\beta \in\D.$ For $\alpha\in\mathbb D^q$ we define the set
\begin{align}\label{def:Jalpha}
\mathcal J(\alpha)
	:=\Big\{ \prod_{i=1}^q I_i, \text{ where } I_i\in\mathcal I(\alpha^i)\cup\big\{\mathbb R\big\} \text{ for every }i=1,\dots,q \Big\}\setminus \big\{\mathbb R^q\big\}.
\end{align}
To every $I:=\prod_{i=1}^q I_i\subset\mathbb R^q$, we associate the set of indexes
\begin{align}\label{def:J_I}
J_I:=\big\{i\in\{1,\dots,q\},\ I_i\neq\mathbb R\}.
\end{align}
For a pair $(\alpha,I)$, where $\alpha\in\mathbb D^q$ and $I:=\prod_{i=1}^qI_i\subset\mathbb R^q$, we define, in analogy to \cite[Notation VI.2.6]{jacod2003limit}, the time points 
\begin{align}\label{def:tp}
s^0(\alpha,I):=0,\ s^{n+1}(\alpha,I):=\inf\{s>s^n(\alpha,I), \Delta \alpha^i_s\in I_i \text{ for every } i \in J_I\},\; n\in\mathbb N.
\end{align}
If $\{s>s^n(\alpha,I),  \Delta \alpha^i_s\in I_i \text{ for every } i \in J_I\}=\emptyset$, then we set $s^{n+1}(\alpha,I):=\infty.$
The value of $s^{n}(\alpha,I)$ marks the $n-$th time at which the value of $\Delta \alpha$ lies in the set $I$. The following proposition is the analogon of \cite[Proposition VI.2.7]{jacod2003limit}. 
\begin{proposition}\label{prop:J1Map}
Fix $q\in\mathbb N.$
\begin{enumerate}[label={\rm(\roman*)}, itemindent=0.8cm, leftmargin=0cm]
\item\label{prop:J1Map-i} The function $\mathbb D^q \ni \alpha \longmapsto s^n(\alpha,I)\in\overline{\mathbb R}_+,$
is continuous at each point $(\alpha,I)\in\mathbb D^q\times \mathcal J(\alpha),$ for $n\in\mathbb N.$

\vspace{0.5em}
\item\label{prop:J1Map-ii} If $s^n(\alpha,I)<\infty,$  for some $n\in\mathbb N,$ $I\in\mathcal J(\alpha)$, then the function
$\mathbb D^q \ni \alpha \longmapsto \Delta \alpha_{s^n(\alpha,I)}\in\mathbb D^q,$
is continuous.
\end{enumerate}
\end{proposition}

\iftoggle{full}{
\begin{proof}
Let $(\alpha^k)_{k\in\overline{\mathbb N}}$ be such that $\alpha^k\inJmulti{q} \alpha^\infty$ and $I:=\prod_{i=1}^qI_i\in\mathcal J(\alpha^\infty).$ 
Observe that $J_I\neq \emptyset,$ by definition of $\mathcal J(\alpha^\infty)$.
We define $s^{k,n}:=s^n(\alpha^k,I)$, for $k\in\overline{\mathbb N}$ and $n\in\mathbb N.$

\vspace{0.5em}
\ref{prop:J1Map-i} The convergence $s^{k,0}\xrightarrow[k\to\infty]{}s^{\infty,0}$ holds by definition. 
Assume that the convergence $s^{k,n}\xrightarrow[k\to\infty]{}s^{\infty,n}$ holds for some $n\in\mathbb N$.
We will prove that the convergence $s^{k,n+1}\xrightarrow[k\to\infty]{}s^{\infty,n+1}$ holds as well.

\vspace{0.5em}
Before we proceed, fix a positive number $u[\alpha^\infty,I]$ such that $u[\alpha^\infty,I]\notin \cup_{i\in J_I}\{|u|, u\in W(\alpha^{\infty,i})\}$ and 
\begin{align*}
u[\alpha^\infty,I]\le\min\cup_{i\in J_I}\{|v|, v\in\partial I_i\footnotemark \}.
\end{align*}\footnotetext{Recall that $\partial A$ denotes the $|\cdot|-$boundary of the set $A\subset \mathbb R$.}
\hspace{-0.2em}Observe now that, for $i\in J_I$ and for $U:=(-\infty, -u[\alpha^\infty,I])\cup(u[\alpha^\infty,I],\infty)$, the sequence $\big( s^l(\alpha^{\infty,i},U)\big)_{l\in\mathbb N}$ exhausts the set of times that $\alpha^{\infty,i}$ exhibits a jump of height greater than $u[\alpha^\infty,I],$ \emph{i.e.}
\begin{align}\label{property:exhaustJumps}
\{t\in\mathbb R_+, |\Delta \alpha^{\infty,i}|>u[\alpha^\infty,I] \}\subset \big( s^l(\alpha^{\infty,i},U)\big)_{l\in\mathbb N}, \text{ for every }i\in J_I.
\end{align}

We will distinguish between several different cases now.

\vspace{0.5em}
\hspace{0.5em}{\bf Case $1$:} \hfil $s^{\infty,n+1}<\infty.$

\vspace{0.5em}
By Property \eqref{property:exhaustJumps}, there exist unique $l^{i,n},l^{i,n+1}\in\mathbb N$ with $l^{i,n} <l^{i,n+1}$ such that 
\begin{align}\label{def:li}
s^{l^{i,n}}(\alpha^{\infty,i},U)=s^{\infty,n} \text{ and }s^{l^{i,n+1}}(\alpha^{\infty,i},U)=s^{\infty,n+1} \text{ for every }i\in J_I .
\end{align}
By \cite[Proposition VI.2.7]{jacod2003limit}\footnote{Observe that for $\alpha\in \mathbb D$ the time point $t^p(\alpha,u)$ defined in \cite[Definition VI.2.6]{jacod2003limit} can be rewritten using our notation as 
$t^p(a,u)=s^p(\alpha,(-\infty,-u)\cup(u,\infty)).$} and the above identities, we obtain 
\begin{align}\label{conv:sli-n}
s^{l^{i,n}}(\alpha^{k,i},U)\xrightarrow[k\to\infty]{}
s^{\infty,n} 
\text{ and }
\Delta \alpha^{k,i}_{s^{l^{i,n}}(\alpha^{k,i},U)} \xrightarrow[k\to\infty]{} \Delta\alpha^{\infty,i}_{s^{\infty,n} }
\text{ for every }i\in J_I,
\end{align}
as well as
\begin{align}\label{conv:sli-n+1}
s^{l^{i,n+1}}(\alpha^{k,i},U)\xrightarrow[k\to\infty]{}
s^{\infty,n+1} 
\text{ and }
\Delta \alpha^{k,i}_{s^{l^{i,n+1}}(\alpha^{k,i},U)} \xrightarrow[k\to\infty]{} \Delta\alpha^{\infty,i}_{s^{\infty,n+1} }
\text{ for every }i\in J_I.
\end{align}

In particular, by \cite[Proposition 2.1.b]{jacod2003limit}, the convergence $s^{k,n}\xrightarrow[k\to\infty]{}s^{\infty,n}$, which has been assumed true as the induction hypothesis, and the Convergence \eqref{conv:sli-n}, we can obtain that
\begin{align}\label{conv:c00-n-induction}
\big(s^{l^{i,n}}(\alpha^{k,i},U) - s^{k,n}\big)_{k\in\mathbb N} \in c_{\text{\tiny $00$}}(\mathbb N), \text{ for every }i\in J_I,
\end{align} 
where $c_{\text{\tiny $00$}}(\mathbb N):=\{ (\gamma^m)_{m\in\mathbb N}\subset \mathbb R^{\mathbb N}, \exists m_0\in\mathbb N \text{ such that } \gamma^m=0 \text{ for every }m\ge m_0\}.$
Define $$k_0^{n,i}:=\max\big\{k\in\mathbb N, s^{l^{i,n}}(\alpha^{k,i},U) \neq s^{k,n}\big\}\text{ for }i\in J_I.$$
By Property \eqref{conv:c00-n-induction}, we obtain that $k_0^{n,i}<\infty$ for every $i\in J_I$.
Since $J_I$ is a finite set, the number $\bar k_0^{n}:=\max\big\{k_0^{n,i}, i\in J_I\big\},$ is well--defined and finite, therefore
\begin{align}\label{def:barkn0}
s^{l^{1,n}}(\alpha^{k,1},U)=\dots=s^{l^{\ell,n}}(\alpha^{k,\ell},U) = s^{k,n}, \text{ for every }k> \bar k_0^n.
\end{align}
Now, in view of Convergence \eqref{conv:sli-n+1}, we can conclude the induction step once we prove the analogue to \eqref{conv:c00-n-induction} for $m=n+1$, \emph{i.e.}
\begin{align}\label{conv:c00-n+1-induction}
(s^{l^{i,n+1}}(\alpha^{k,i},U) - s^{k,n+1})_{k\in\mathbb N}\in c_{\text{\tiny $00$}}(\mathbb N), \text{ for every }i\in J_I.
\end{align} 
At this point we further distinguish between two cases.

\vspace{0.5em}
\hspace{1.5em}{\bf Case $1.1$:}\hfil For every $i\in J_I$, we have $l^{i,n+1}=1+l^{i,n}$.

\vspace{0.5em}
By \cite[Proposition VI.2.1.b]{jacod2003limit}, Convergence \eqref{conv:sli-n+1} and the convergence $\alpha^{k}\inJmulti{q}\alpha^\infty$, we can conclude that 
\begin{align}\label{conv:c00-n+1-ij}
\big(s^{l^{i,n+1}}(\alpha^{k,i},U) - s^{l^{j,n+1}}(\alpha^{k,j},U)\big)\in c_{\text{\tiny $00$}}(\mathbb N), \text{ for every }i,j\in J_I.
\end{align}
Therefore, we can fix hereinafter an index from $J_I$ and we will do so for $\mu:=\min J_I$, \emph{i.e.} $\mu$ is the minimum element of $J_I$.
Define $$k_0^{n+1,i}:=\max\big\{k\in\mathbb N:s^{l^{i,n+1}}(\alpha^{k,i},U) \neq s^{l^{\mu,n+1}}(\alpha^{k,\mu},U)\big\},\text{ for }i\in J_I\setminus\{\mu\}.$$
By Property \eqref{conv:c00-n+1-ij}, we obtain that $k_0^{n+1,i}<\infty$ for every $i\in J_I\setminus\{\mu\}$.
Since $J_I\setminus\{\mu\}$ is a finite set, the number $\bar k_0^{n+1}:=\max\big\{k_0^{n+1,i}, i\in J_I\setminus\{\mu\}\big\},$ is well--defined and finite. Observe that 
\begin{align}\label{property:finallylMu}
s^{1+l^{i,n}}(\alpha^{k,i},U)=s^{l^{i,n+1}}(\alpha^{k,i},U) = s^{l^{\mu,n+1}}(\alpha^{k,\mu},U), \text{ for } k> \bar k_0^{n+1} \text{ and for }i\in J_I,
\end{align}
where the first equality holds by assumption.
Moreover, by Convergence \eqref{conv:sli-n+1} we obtain that 
\begin{align}\label{conv:finallyinIi}
1=\prod_{i\in J_I}\mathds{1}_{I_i}\Big(\Delta \alpha^{k,i}_{s^{1+l^{i,n}}(\alpha^{k,i},U)}\Big), \text{ for all but finitely many }k,
\end{align}
since $\Delta \alpha^{\infty,i}_{s^{\infty,n+1}}$ lies in the interior of the open interval $I_i$, for every $i\in J_I.$
For notational convenience, we will assume that the above convergence holds for $k>\bar k_0^{n+1}$.
Therefore, 
\begin{align*}
s^{k,n+1}	&=\inf\{s>s^{k,n}, \Delta\alpha^{k,i}_s\in I_i \text{ for every } i\in J_I \}
			 =\inf\bigcap_{i\in J_I}\{s>s^{k,n}, \Delta\alpha^{k,i}_s\in I_i \}\\
			&\hspace{-0.4em}\overset{\eqref{conv:c00-n-induction}}{=}\inf\bigcap_{i\in J_I}\{s>s^{l^{i,n}}(\alpha^{k,i},U),\Delta\alpha^{k,i}_s\in I_i  \}, \text{ for }k> \bar k_0^n\\
			&=\inf\bigcap_{i\in J_I}\{s\in\big(s^{l^{i,n}}(\alpha^{k,i},U),s^{1+l^{i,n}}(\alpha^{k,i},U)\big], \Delta \alpha^{k,i}_s\in I_i\}\\
			&\hspace{4em}\wedge\inf\bigcap_{i\in J_I}\{s>s^{1+l^{i,n}}(\alpha^{k,i},U),\Delta\alpha^{k,i}_s\in I_i  \}\big), \text{ for }k> \bar k_0^n\\
			&\hspace{-0.6em}\overset{\eqref{conv:finallyinIi}}{=}\bigcap_{i\in J_I}\{s^{1+l^{i,n}}(\alpha^{k,i},U)\}\wedge\inf\bigcap_{i\in J_I}\{s>s^{1+l^{i,n}}(\alpha^{k,i},U),\Delta\alpha^{k,i}_s\in I_i  \}\big),\text{ for }k> \bar k_0^n\vee \bar k_0^{n+1}\\
			&\hspace{-0.6em}\overset{\eqref{property:finallylMu}}{=}s^{l^{\mu,n+1}}(\alpha^{k,i},U), \text{ for }k> \bar k_0^n\vee\bar k_0^{n+1}=s^{l^{i,n+1}}(\alpha^{k,i},U), \text{ for }k> \bar k_0^n\vee\bar k_0^{n+1}, i\in J_I,
\end{align*}
\emph{i.e.} Property \eqref{conv:c00-n+1-induction} holds.

\vspace{0.5em}
\hspace{1.5em}{\bf Case $1.2$:}\hfil There exists $i\in J_I$ for which $l^{i,n+1}>1+l^{i,n}$.

\vspace{0.5em}
Define $\bar J_I:=\{i\in J_I, l^{i,n+1}>1+l^{i,n}\}$
and fix $\xi^i\in(l^{i,n},l^{i,n+1})\cap\mathbb N,$ for every $i\in\bar J_I.$
Recall that by \cite[Proposition VI.2.7]{jacod2003limit}, we have
\begin{align}\label{conv:xi-i}
\lim_{k\to\infty}s^{\xi^i}(\alpha^{k,i},U)=s^{\xi^i}(\alpha^{\infty,i},U) \text{ and } \lim_{k\to\infty}\Delta \alpha^{k,i}_{s^{\xi^i}(\alpha^{k,i},U)}
= \Delta\alpha^{\infty,i}_{s^{\xi^i}(\alpha^{\infty,i},U)}.
\end{align}
We can conclude that Property \eqref{conv:c00-n+1-induction} holds, if for every $\bar s\in(s^{\infty,n},s^{\infty,n+1})$ such that
\begin{align}\label{conv:xi-itoxi}
\lim_{k\to\infty}s^{\xi^i}(\alpha^{k,i},U)=\bar s, \text{ for every }i\in\bar J_I,
\end{align}
holds 
\begin{align}\label{conv:finally0}
0	=	\prod_{i\in J_I\setminus \bar J_I}\mathds{1}_{I_i}\big(\Delta \alpha^{k,i}_{\bar s}\big) 
		\prod_{i\in \bar J_I}\mathds{1}_{I_i}\Big(\Delta \alpha^{k,i}_{s^{\xi^i}(\alpha^{k,i},U)}\Big),\text{ for all but finitely many }k.
\end{align}
However, if we had
\begin{align*}
1	=	\prod_{i\in J_I\setminus \bar J_I}\mathds{1}_{I_i}\big(\Delta \alpha^{k,i}_{\bar s}\big) 
		\prod_{i\in \bar J_I}\mathds{1}_{I_i}\Big(\Delta \alpha^{k,i}_{s^{\xi^i}(\alpha^{k,i},U)}\Big),\text{ for all but finitely many }k,
\end{align*}
then in view of the definition of $s^{\infty,n+1}$, we would have $s^{\infty,n+1}=s^{\xi^i}(\alpha^{\infty,i},U)$ for every $i\in\bar J_I$, which contradicts Property \eqref{def:li}.
The contradiction arises in view of 
\begin{align*}
s^{\xi^i}(\alpha^{\infty,i},U) < s^{l^{i,n+1}}(\alpha^{k,i},U), \text{ since } \xi^i<l^{i,n+1}.
\end{align*}

\vspace{0.5em}
\hspace{0.5em}{\bf Case $2$:} \hfil $s^{\infty,n+1}=\infty.$

\vspace{0.5em}
We distinguish again between two cases.

\vspace{0.5em}
\hspace{1.5em}{\bf Case $2.1$:}\hfil $s^{\infty,n}<\infty$

\vspace{0.5em}
Using the same arguments as the ones used in Property \eqref{def:li} for $s^{\infty,n}$, we can associate to $s^n(\alpha^{\infty,i},U)$ a unique natural number $l^{i,n}$ such that Convergence \eqref{conv:sli-n} holds.
Moreover, by definition of $s^{\infty,n+1}$ we obtain that $$\{s>s^{n}(\alpha^\infty,I),  \Delta \alpha^{\infty,i}_s\in I_i \text{ for every } i \in J_I\}=\emptyset,$$ or, equivalently
\begin{align}\label{PropertyforContr}
\text{for every }s>s^{\infty,n}\text{ there exists an }i\in J_I\text{ such that }\Delta \alpha^{\infty,i}\not\in I_i.
\end{align}

\vspace{0.5em}
Assume now that $$\liminf_{k\to\infty}s^{k,n+1}=\bar s,\text{ for some }\bar s\in(s^{k,n},\infty),$$ 
\emph{i.e.} there exists $(k_l)_{l\in\mathbb N}$ such that 
$s^{k_l,n+1}\xrightarrow[l\to\infty]{}\bar s$. 
Equivalently, for every $i\in J_I$ holds $\Delta \alpha^{k_l,i}_{s^{k_l,n+1}}\in I_i$ for all but finitely many $l.$
By convergence $\alpha^k\inJmulti{q} \alpha^\infty$, \cite[Proposition VI.2.1]{jacod2003limit} and convergence $s^{k_l,n+1}\xrightarrow[l\to\infty]{}\bar s$, we have that $\Delta \alpha^{\infty,i}_{\bar s}\in I_i,$ for every $i\in J_I.$
But this contradicts the assumption $s^{\infty,n+1}=\infty$ in view of the equivalent form \eqref{PropertyforContr}.

\vspace{0.5em}
\hspace{1.5em}{\bf Case $2.2$:}\hfil $s^{\infty,n}=\infty.$

\vspace{0.5em}
By the induction hypothesis it holds that
$s^{k,n}\xrightarrow[k\to\infty]{}s^{\infty,n}$,
and by the definition of $s^{k,n+1}$ it holds that 
$s^{k,n}\le s^{k,n+1}, \text{ for every } k\in\overline{\mathbb N},$ $n\in\mathbb N.$
The previous convergence yields
$s^{k,n+1}\xrightarrow[k\to\infty]{}\infty=s^{\infty,n+1}.$

\vspace{0.5em}
\ref{prop:J1Map-ii} Assume that there exist $n\in\mathbb N$ and $I\in\mathcal J(\alpha)$ such that $s^{\infty,n}<\infty.$
By \ref{prop:J1Map-i}, we have that $s^{k,n}\xrightarrow[k\to\infty]{}s^{\infty,n}$, which in conjunction with the convergence $\alpha^k\inJmulti{q}\alpha^\infty$ and \cite[Proposition VI.2.1]{jacod2003limit} implies that 
\begin{equation*}
\Delta \alpha^{k}_{s^{k,n}} \xrightarrow[k\to\infty]{}\Delta \alpha^{\infty}_{s^{\infty,n}}. \qedhere
\end{equation*}   
\end{proof}
}
{
\begin{proof}
See the proof of Proposition A.5 in the extended version \cite{papapantoleon2018stability}.	
\end{proof}
}

\begin{corollary}\label{cor:J1Map-2}
Let $\alpha^k\inJmulti{q}\alpha^\infty$ and $I\in\mathcal J(\alpha^\infty)$.
Define
$$\hat n:=
\begin{cases}
\max\{n\in\mathbb N, s^n(\alpha^\infty, I)<\infty\}, \text{ if }\{n\in\mathbb N, s^n(\alpha^\infty, I)<\infty\} \text{ is non--empty and finite},\\
\infty, \text{ otherwise}.
\end{cases}$$
Then for every function $g:\mathbb R^q\longrightarrow \mathbb R$ which is continuous on 
$$C:=\prod_{i=1}^q A_i, \text{ where } A_i:=
\begin{cases}
W(\alpha^{\infty,i}), 			\text{ if }i\in J_I,\\
W(\alpha^{\infty,i})\cup\{0\}, 	\text{ if }i\in\{1,\dots,q\}\setminus J_I,
\end{cases}$$ 
and for every $0\le n\le \hat n$ holds
$$g(\Delta \alpha^k_{s^n(\alpha^k,I)}) \xrightarrow[k\to\infty]{}g(\Delta\alpha^\infty_{s^n(\alpha^\infty, I)}).$$
In particular, for the \cadlag functions
\begin{align*}
\beta^k_{\cdot}:=g(\Delta \alpha^k_{s^n(\alpha^k,I)}) \mathds{1}_{[s^n(\alpha^k,I),\infty)}(\cdot), \text{ for }k\in\overline{\mathbb N},
\end{align*}
the convergence $\beta^k\inJmulti{}\beta^\infty$ holds.
\end{corollary}

\begin{proof}
\iftoggle{full}{
Fix an $n\in\mathbb N$ such that $n\le \hat n.$
By Proposition \ref{prop:J1Map}.\ref{prop:J1Map-ii} holds 
\begin{align*}
\Delta \alpha^k_{s^n(\alpha^k,I)}\xrightarrow[k\to\infty]{}\Delta\alpha^\infty_{s^n(\alpha,\infty,I)},
\text{ where }
\Delta \alpha^{\infty,i}_{s^n(\alpha^\infty,I)}\in 
\begin{cases}
W(\alpha^{\infty,i}), 			\text{ if } i\in J_I,\\
W(\alpha^{\infty,i})\cup\{0\},	\text{ if } i\in\{1,\dots,q\}\setminus J_I.
\end{cases}
\end{align*}
Therefore, by definition of the time $s^n(\alpha^\infty,I)$ and of the set $C$ holds
\begin{align}\label{conv:J1Map-2-ii}
g(\Delta \alpha^k_{s^n(\alpha^k,I)})\xrightarrow[k\to\infty]{}g(\Delta\alpha^\infty_{s^n(\alpha^\infty,I)}).
\end{align}
By Proposition \ref{prop:J1Map}.\ref{prop:J1Map-i}, the above convergence and \cite[Example VI.1.19]{jacod2003limit} we obtain the convergence 
\begin{equation*}
\beta^k\inJmulti{}\beta^\infty.	\qedhere
\end{equation*}
}
{
	See the proof of Corollary A.6 in the extended version \cite{papapantoleon2018stability}.
}
\end{proof}
\iftoggle{full}{The following simple counterexample shows that in Corollary \ref{cor:J1Map-2} the convergence $\beta^k\inJmulti{}\beta^\infty$ does not necessarily hold for a function $g$ which is not continuous on the set $C$.
\begin{example}
Let $\big((t_k,x_k)\big)_{k\in\overline{\N}}$ be such that $(t_k,x_k)\in\Rp\times\R{}$ for every $k\in\overline{\N}$, $t_k\longrightarrow t_{\infty}\in\R{}$ and $x_k\searrow x_{\infty}\in\R{}\setminus\{0\}.$ 
Define $\gamma^k_{\cdot}:=x_k\mathds{1}_{[t_k,\infty)}(\cdot)$, for every $k\in\overline{\N}$.
By {\rm \cite[Example VI.1.19.ii)]{jacod2003limit}}, we have
$\gamma^k\Jconv{} \gamma^{\infty}.$ 
On the other hand, for $I:=(\frac{1}{2}x_\infty,\frac{3}{2}x_\infty)$ holds $s^1(\gamma^k,I)=t_k$ for all but finitely many $k$ and $s^1(\gamma^\infty,I)=t_\infty$, \emph{i.e.} $s^1(\gamma^k,I)\xrightarrow[k\to\infty]{}s^1(\gamma^\infty,I)$.
Moreover, for $w>x_{\infty}$ we also have
$$\Delta\gamma^k_{t_k}\mathds{1}_{(x_{\infty},w)}(\Delta\gamma^k_{t_k})= x_k\mathds{1}_{(x_{\infty},w)}(x_k)=x_k\ \text{ for all but finitely many }k\in\N,$$
and
$$\Delta\gamma^{\infty}_{t_{\infty}}\mathds{1}_{(x_{\infty},w)}(\Delta\gamma^{\infty}_{t_{\infty}})=x_{\infty}\mathds{1}_{(x_{\infty},w)}(x_{\infty})=0.$$
Therefore, for $\mathbb R\ni x\overset g \longmapsto x \mathds{1}_{(x_\infty,w)}(x)$
\begin{align*}
g(\Delta \gamma^k_{s^1(\gamma^k,I)})=\Delta\gamma^k_{t_k}\mathds{1}_{(x_{\infty},w)}(\Delta\gamma^k_{t_k})\xrightarrow[n\to\infty]{} x_{\infty}\neq 0 = 
\Delta\gamma^{\infty}_{t_{\infty}}\mathds{1}_{(x_{\infty},w)}(\Delta\gamma^{\infty}_{t_{\infty}})=g(\Delta\gamma^{\infty}_{s^1(\alpha^\infty)}),
\end{align*}
and for this reason we \emph{cannot} obtain the convergence
$$\Delta\gamma^k_{t_k}\mathds{1}_{(x_{\infty},w)}(\Delta\gamma^n_{t_k})\mathds{1}_{(x_{\infty},w)}(\cdot)\Jconv{}
\Delta\gamma^{\infty}_{t_{\infty}}\mathds{1}_{(x_{\infty},w)}(\Delta\gamma^{\infty}_{t_{\infty}})\mathds{1}_{(x_{\infty},w)}(\cdot).$$
\end{example}}{}
%
%
%
\begin{proposition}\label{prop:ContFuncInSkor}
Fix some subset $I:=\prod_{i=1}^qI_i$ of $\mathbb R^q$ and a function $g:\mathbb R^q\rightarrow\mathbb R$.
Define the map 
$$\mathbb D(\mathbb R^q)\ni\alpha\longmapsto \hat\alpha[g,I]:=(\alpha^1,\dots,\alpha^q,\alpha^{g,I})^\top\in\mathbb D(\R{q+1}),$$
where 
\begin{align}\label{def:hatAlpha-gI}
\alpha^{g,I}_{\cdot} := \sum_{0< t\leq\cdot}g(\Delta \alpha_t)\mathds{1}_{I}(\Delta \alpha_t). 
\end{align} 
Then, the map $\hat \alpha[g,I]$ is $\J_1-$continuous at each point $\alpha$ for which $I\in\mathcal J(\alpha)$ and for each function $g$ which is continuous on the set
$$C:=\prod_{i=1}^q A_i, \text{ where } A_i:=
\begin{cases}
W(\alpha^{i}), 			\text{ if }i\in J_I,\\
W(\alpha^{i})\cup\{0\}, 	 \text{ if }i\in\{1,\dots,q\}\setminus J_I.
\end{cases}$$  
\end{proposition}
\begin{proof}
The arguments are similar to those in the proof of Corollary VI.2.8 in \cite{jacod2003limit}, therefore they are omitted for the sake of brevity.
The interested reader can also consult \citet[Proposition~I.134]{saplaouras2017backward} for the corresponding proof.
\end{proof}

\iftoggle{full}
{
\subsubsection{\texorpdfstring{Technical proofs}{Technical proofs}}

We conclude this subsection with the proofs that were omitted from Section \ref{sec:skor}.
\begin{proof}[Proof of Proposition \ref{MeminProp1}]
The first step is to show that the $\mathbb{L}^1$ convergence of $(M^k_\infty)_k$ together with the weak convergence of the filtrations, imply the convergence of the martingales in the $\J_1(\R q)-$topology.
Let $\varepsilon>0$ and $\mathbb F^k\weakFil \mathbb F^\infty$, then
\begin{align*}
\mathbb P\Big( d_{\J_1(\R{q})}(M^k,M^\infty)>\varepsilon\Big) 
&\leq 				\mathbb P 				\Big( d_{\J_1(\R{q})}\big(M^k,\mathbb E[M^\infty_\infty|\mathcal F_{\cdot}^{k}]\big)>\frac{\varepsilon}{2}\Big)+		\Pm\Big( d_{\J_1(\R{q})}\big(\mathbb E[M^\infty_\infty|\mathcal F_{\cdot}^{k}], M^\infty\big)>\frac{\varepsilon}{2}\Big)\\
&\leq				\Pm\Big( \sup_{t\in[0,\infty)}\big|M^k_t-\mathbb E[M^\infty_\infty|\mathcal F_{\cdot}^{k}]\big|>\frac{\varepsilon}{2}\Big)+
					\Pm\Big( d_{\J_1(\R{q})}\big(\mathbb E[M^\infty_\infty|\mathcal F_{\cdot}^{k}], M^\infty\big)>\frac{\varepsilon}{2}\Big)\\
&\leq	\frac{2}{\varepsilon}
					\E \big[ |M_\infty^k - M^\infty_\infty|\big] +		\Pm\Big( d_{\J_1(\R{q})}\big(\mathbb E[M^\infty_\infty|\mathcal F_{\cdot}^{k}], \mathbb E[M^\infty_\infty|\mathcal F_\cdot^\infty]\big)>\frac{\varepsilon}{2}\Big)
\numberthis \label{ConvinSiP}
\xrightarrow{\hspace{0.1cm} k\to\infty\hspace{0.1cm}} 0,
\end{align*}
where the first summand converges to $0$ by assumption and the second one by the weak convergence of the filtrations.
Let us point out that for the second inequality we have used that for $\alpha, \beta\in\D^{q}$ holds 
\begin{align*}
d_{\J_1(\R{q})}(\alpha,\beta)\leq d_{\lu}(\alpha,\beta) \leq d_{\norm{\cdot}_\infty}(\alpha,\beta), 
\end{align*}
by the definition of the metrics, while for the third inequality we used Doob's martingale inequality.

\vspace{0.5em}
The next step is to apply Lemma \ref{JointSkorokhodConv} to $(M^k, \E[\xi|\F^k_\cdot])_k^\top=:(N^k)_k$, for $\xi\in\mathbb{L}^1(\Omega,\mathcal F^\infty_\infty,\Pm;\mathbb R)$, in order to obtain the convergence in the extended sense. 
The $\J_1(\R{})-$convergence of each $(N^{k,i})_k$, for $i=1,\dots,q$, and of the partial sums $(\sum_{i=1}^p N^{k,i})_k$, for $p=1,\dots,q$ follows from the previous step and Lemma \ref{JointSkorokhodConv}. 
Moreover, the $\J_1(\R{})-$convergence of $N^{k,q+1}$ follows from the definition of the weak convergence of filtrations.
Hence, we just have to show the $\J_1(\R{})-$convergence of $(\sum_{i=1}^{q+1} N^{k,i})_k$.

\vspace{0.5em}
By assumption, we have 
$
\sum_{m=1}^{q}M_\infty^{k,m} + \xi\xrightarrow[k\to\infty]{\hspace{0.2cm}\mathbb{L}^1(\Omega,\mathcal G,\Pm;\mathbb R)\hspace{0.2cm}} 
\sum_{m=1}^{q}M_\infty^{\infty,m} + \xi,$
and arguing as in \eqref{ConvinSiP} we obtain  
\begin{align*}
\mathbb E\bigg[\sum_{m=1}^{q}M_\infty^{k,m} + \xi \bigg|\mathcal F_{\cdot}^{k}\bigg] 
\underset{k\to\infty}{\sip} 
\mathbb E\bigg[\sum_{m=1}^{q}M_\infty^{\infty,m} + \xi\bigg|\mathcal F_{\cdot}^\infty\bigg].
\end{align*}
Moreover, by the linearity of conditional expectations, we get that
\[
	\sum_{m=1}^{q} \mathbb E\big[M_\infty^{k,m} \big| \mathcal F_{\cdot}^{k}\big]
	+ \E\big[ \xi \big| \mathcal F_{\cdot}^{k}\big]
	\underset{k\to\infty}{\sip} 
	\sum_{m=1}^{q} \mathbb E\big[M_\infty^{m} \big| \mathcal F_{\cdot}^\infty\big]
	+ \E\big[ \xi \big| \mathcal F_{\cdot}^{\infty}\big].	
\]

\vspace{0.5em}
The converse statement is trivial.
\end{proof}

\vspace{0.5em}
\begin{proof}[Proof of Theorem \ref{MeminCorollary}]
{\color{black}\ref{MeminCorollaryJ1P}} By \cite[Corollary 12]{memin2003stability}, we obtain for every $i=1,\ldots,q$ the convergence 
\[
	(M^{k,i}, [M^k]^{ii} )^{\transp}	
			\sipmulti{2}
	(M^{\infty,i}, [M^\infty]^{ii} )^{\transp},
\]
which in conjunction with Corollary \ref{cormultiskorokhod} and the convergence $M^k\sipmulti{q} M^\infty$ implies
\begin{align}\label{DiagQCConv}
\left(
\begin{array}{c c}
M^{k,1} & [M^k]^{11} \\
M^{k,2} & [M^k]^{22} \\
\vdots&\vdots\\
M^{k,q} & [M^k]^{qq} \\
\end{array}
\right)
\sipmulti{q\times 2}
\left(
\begin{array}{c c c}
M^{\infty,1} & [M^\infty]^{11} \\
M^{\infty,2} & [M^\infty]^{22} \\
\vdots&\vdots\\
M^{\infty,q} & [M^\infty]^{qq} \\
\end{array}
\right).
\end{align}
On the other hand, let $i,j\in\{1,\ldots,q\}$ with $i\neq j$. 
Using that the sequence of square integrable martingales $(M^{k,l})$, for every $l=1,\dots,q$, is $\mathbb L^2-$bounded, Doob's maximal inequality and \cite[Corollary VI.6.30]{jacod2003limit}, we get that the sequences $(M^{k,i})_{k\in\overline{\N}}$ and $(M^{k,j})_{k\in\overline{\N}}$ possess the $\PUT$ property; see \cite[Section VI.6]{jacod2003limit}.
Therefore, by \cite[Theorem VI.6.22]{jacod2003limit}, we obtain
\begin{multline}\label{QCijPUT}
(M^{k,i}, M^{k,j}, M^{k,i}_{-}\cdot M^{k,j}, M^{k,j}_{-}\cdot M^{k,i})^{\transp}
	\sipmulti{4}
(M^{\infty,i}, M^{\infty,j}, M^{\infty,i}_{-}\cdot M^{\infty,j}, M^{\infty,j}_{-}\cdot M^{\infty,i})^{\transp}.
\end{multline}

Now, in order to show that the quadratic variation of $M^{k,i}$ and $M^{k,j}$ converges, we just need to show that the product $M^{k,i}M^{k,j}$ converges, by the definition of the quadratic variation, see \cite[Definition I.4.45]{jacod2003limit}. 

By 
the convergence $(M^{k,i}, M^{k,j})^{\transp}\sipmulti{2} (M^{\infty,i}, M^{\infty,j})^{\transp}$, recall \eqref{DiagQCConv}, we obtain the convergence 
\begin{align}\label{QCijProd}
(M^{k,i}, M^{k,j}, M^{k,i}M^{k,j})^{\transp}\sipmulti{3} (M^{\infty,i}, M^{\infty,j}, M^{\infty,i}M^{\infty,j})^{\transp},
\end{align} 
where we have applied Lemma \ref{lem:ToolforJoint} for the continuous function $\mathbb R^2\ni (x^1,x^2)^\top\longmapsto (x^1,x^2,x^1x^2)^\top\in\mathbb R^3$.
Then \eqref{QCijPUT}--\eqref{QCijProd} imply that
\begin{align}\label{OffDiagQCConv}
(M^{k,i}, M^{k,j},[M^k]^{ij})^{\transp}
\sipmulti{3}
(M^{\infty,i}, M^{\infty,j},[M^\infty]^{ij})^{\transp},
\end{align}
while \eqref{DiagQCConv} and \eqref{OffDiagQCConv}, in conjunction with Remark \ref{rem:cormultiskorokhod}, yield that 
\begin{align}\label{QCConv}
(M^k, [M^k]) 
\xrightarrow{\hspace{0.2cm} \left(\J_1(\R{q}\times\R{q\times q}),\Pm\right)\hspace{0.2cm}} 
(M^\infty, [M^\infty]).
\end{align}

\vspace{0.5em}
Let us now show that the predictable quadratic variations converge as well.
By \cite[Corollary 12]{memin2003stability}, we have for $i=1,\ldots,q$ 
\begin{align}\label{DiagPredQCConv}
\langle M^k\rangle^{ii} \sipmulti{} \langle M^\infty \rangle^{ii}. 
\end{align}
Moreover the convergence $M^k_\infty \silmulti{2}{2} M^\infty_\infty$ implies in particular, for every $i,j=1,2,\ldots,q$ with $i\neq j$, that  
\begin{align}\label{FVijConv}
M^{k,i}_\infty+M^{k,j}_\infty \silmulti{2}{2} M^{\infty,i}_\infty+M^{\infty,j}_\infty.
\end{align}
In view of \eqref{DiagQCConv} and \eqref{FVijConv}, we can apply \cite[Corollary 12]{memin2003stability} to $(M^{k,i}+M^{k,j})_{k\in\overline{\N}}$ and $(M^{k,i}-M^{k,j})_{k\in\overline{\N}}$, for every $i,j=1,2,\ldots,q$ with $i\neq j$.
Therefore we get that
\begin{align*}
\langle M^{k,i}+M^{k,j}\rangle \sipmulti{} \langle M^{\infty,i}+M^{\infty,j}\rangle,\\
\text{ and }\ \langle M^{k,i}-M^{k,j}\rangle \sipmulti{} \langle M^{\infty,i}-M^{\infty,j}\rangle.
\end{align*} 
Now recall that $M^\infty$ is quasi--left--continuous, which implies that the processes $\langle M^{\infty,i}+M^{\infty,j} \rangle$ and $\langle M^{\infty,i}-M^{\infty,j} \rangle$ are continuous for every $i,j=1,2,\ldots,q$ with $i\neq j$. 
Therefore, by \cite[Theorem I.4.2, Proposition VI.2.2]{jacod2003limit} and the last results we obtain 
\begin{multline}
\langle{M^k}\rangle^{ij}=\frac{1}{4}(\langle M^{k,i}+M^{k,j}\rangle - \langle M^{k,i}-M^{k,j}\rangle) 
\sipmulti{}\\
\frac{1}{4}(\langle M^{\infty,i}+M^{\infty,j}\rangle - \langle M^{\infty,i}-M^{\infty,j}\rangle)=\langle{M^\infty}\rangle^{ij}.
\numberthis \label{OffDiagPredQCConv}
\end{multline}
Concluding, by \eqref{QCConv}, \eqref{DiagPredQCConv}, \eqref{OffDiagPredQCConv} and due to the continuity of $\langle M^\infty\rangle$ we have 
\begin{align*}
(M^k, [M^k], \langle M^k\rangle)	
\xrightarrow{\hspace{0.2cm} \left(\J_1(\R{q}\times\R{q\times q}\times\R{q\times q}),\Pm\right)\hspace{0.2cm}} 
(M^\infty, [M^\infty], \langle M^\infty\rangle).
\end{align*}

\vspace{0.3em}
\ref{BracketConvinL} Let $i=1,\ldots,q.$
	The sequence $(M^{k,i})_{k\in\overline{\N}}$ satisfies the conditions of \cite[Corollary 12]{memin2003stability}. 
	In the middle of the proof of the aforementioned corollary, we can find the required convergences.
\end{proof}

\vspace{0.5em}
\begin{proof}[Proof of Lemma \ref{UIplusL2Bounded}]
By \cite[Corollary 1.10]{he1992semimartingale}, it is enough to prove that the sequence $\big({\rm \Var}\big([L^{k,i},N^{k,j}]\big)_\infty\big)_{k\in\N}$ is uniformly integrable, for every $i=1,\dots,p$ and $j=1,\dots,q.$
We will use \cite[Theorem 1.9]{he1992semimartingale} in order to prove it. 
Let $i,j$ be arbitrary but fixed. 
The $\mathbb{L}^1$--boundedness of the sequence $({\rm Var}([L^{k,i},N^{k,j}])_\infty)_{k\in\N}$ is obtained using the Kunita--Watanabe inequality in the form \cite[Corollary 6.34]{he1992semimartingale} 
and the $\mathbb L^1-$boundedness of the sequences $\big([L^{k,i}]_\infty \big)_{k\in\N}$ and $\big([N^{k,j}]_\infty \big)_{k\in\N}$, due to the uniform integrability of the former; see \cite[Theorem 1.7.1)]{he1992semimartingale}. 
By the Kunita--Watanabe inequality again, but now in the form \cite[Theorem 6.33]{he1992semimartingale}, and the Cauchy--Schwarz inequality, we obtain for any $A\in\mathcal G$
\begin{align*}
\int_A \Var\big([L^{k,i},N^{k,j}]\big)_\infty\ \dP 
&\hspace{0cm}\overset{\textrm{K-W ineq.}}{\leq}
\int_A [L^{k,i}]_\infty^{\frac{1}{2}} [N^{k,j}]_\infty^{\frac{1}{2}} \dP\\
&\hspace{0.07cm}\overset{\textrm{C-S ineq.}}{\leq}
\bigg(\int_A [L^{k,i}]_\infty\ \dP\bigg)^{\frac{1}{2}} \bigg(\int_A [N^{k,j}]_\infty\ \dP\bigg)^{\frac{1}{2}}\\
&\overset{\hspace{1.18cm}}{\leq}  
 \bigg(\int_A [L^{k,i}]_\infty\ \dP\bigg)^{\frac{1}{2}}
 \bigg(\sup_{k\in\N}\int_{\Omega}[N^{k,j}]_\infty\ \dP \bigg)^{\frac{1}{2}}\\
&\overset{\hspace{1.18cm}}{\leq}
 C  \bigg(\int_A [L^{k,i}]_\infty\ \dP\bigg)^{\frac{1}{2}},
 \numberthis \label{UIVarBound}
 \end{align*}
where $C^2:=\underset{k\in\N}{\sup}\ \E\Big[{\rm Tr}\big[[N^k]_\infty\big]\Big]$.
Note that for the third inequality we used that $[N^{k,j}]_\infty\ge 0$ for all $k\in\N$.

\vspace{0.5em}
Now we use \cite[Theorem 1.9]{he1992semimartingale} for the uniformly integrable sequence $\big([L^{k,i}]_\infty \big)_{k\in\N}$. For every $\varepsilon>0$, there exists $\delta>0$ such that, whenever $\Pm(A)<\delta$, it holds
$$\underset{k\in\N}{\sup}\ \int_A [L^{k,i}]_\infty\ \dP<\frac{\varepsilon^2}{C^2}
\Longleftrightarrow
\underset{k\in\N}{\sup}\bigg(\int_A [L^{k,i}]_\infty\ \dP\bigg)^{1/2}<\frac{\varepsilon}{C}.$$ 
This further implies
$\underset{k\in\N}{\sup}\,\E \big[\mathds{1}_A\Var\big([L^{k,i},N^{k,j}]_\infty\big)\big]\overset{\eqref{UIVarBound}}{<}\varepsilon$, which is the required condition.
\end{proof}

}{}

\iftoggle{full}{

\subsection{Measure Theory}

This appendix contains some necessary results of measure theoretic nature.
\medskip
Let us endow the Euclidean space $\R{\ell}$ with the usual metric $d_{|\cdot|},$ where we have suppressed the dimension in the notation. 
%
Let $D$ be a countable and dense subset of $\R{}.$ 
Then every open $U\subset\R{}$ is the countable union of pairwise disjoint intervals with endpoints in $D,$ \textit{i.e.} there exists a sequence of intervals $\big((a_k,b_k)\big)_{k\in\N}$ with $a_k,b_k\in D$, for every $k\in\N,$ such that $U=\cup_{k\in\N}(a_k,b_k).$  
In particular, these intervals can be selected from the set $\mathcal I(X^{\infty,i})$, for every $i=1,\dots,\ell$.

\begin{lemma}\label{IXGeneratesBorel}
The following equality holds $\sigma\big(\mathcal{J}(X^\infty)\big)=\borel{\mathbb R^{\ell}},$ where $\mathcal{J}(X^\infty)$ has been introduced in Subsection \ref{subsec:StepD} and $\sigma\big(\mathcal{J}(X^\infty)\big)$ denotes the $\sigma-$algebra on $\mathbb R^\ell$ generated by the family $\mathcal{J}(X^\infty).$
\end{lemma}
\begin{proof}
The space $\mathbb R^\ell$ is a finite product of the second countable metric space $\mathbb R$.
Therefore, it holds that $$\mathcal B(\mathbb R^\ell)=\underbrace{\mathcal B(\mathbb R)\otimes\dots\otimes\mathcal B(\mathbb R)}_{\ell-\text{times}},$$ so that it is sufficient to prove that $\sigma(\mathcal I(X^{\infty,i}))=\mathcal B(\mathbb R)$ for every $i=1,\dots,\ell.$

\vspace{0.5em}
Using the comment at the beginning of this appendix, we have the following chain of equalities
\begin{align*}
\sigma(\mathcal I(X^{\infty,i})) 
& = \sigma(\{U\subset \mathbb R, U \text{ is open and subset of } \R{} \setminus\{0\}\} )
  = \sigma(\{U\subset \mathbb R\setminus\{0\}, U \text{ is open}\}\cup\{\mathbb R\}\cup\{0\} )\\
& = \sigma(\{U\subset \mathbb R, U \text{ is open}\} ) = \mathcal B(\mathbb R),
\end{align*} 
which allows us to conclude.
\end{proof}

\begin{lemma}\label{lem:pi-lambda-signed}
Let $(\Sigma, \mathcal S)$ be a measurable space and $\varrho_1,$ $\varrho_2$ be two finite signed measures on this space. 
Let $\mathcal A$ be a family of sets with the following properties:
\begin{enumerate}[label={\rm (\roman*)}]
\item $\mathcal A$ is a $\pi-$system, i.e. if $A,B\in \mathcal A$, then $A\cap B\in\mathcal A.$
\item We have $\sigma(\mathcal A) = \mathcal S$.
\item $\varrho_1(A)=\varrho_2(A)$ for every $A\in\mathcal A.$
\end{enumerate}
Then it holds $\varrho_1=\varrho_2$.
\end{lemma}

\begin{proof}
See \citet[Lemma~7.2.1]{bogachev2007measure}, or \citet[Lemma~A.7]{saplaouras2017backward}.
\end{proof}

} {}

\subsection{Young functions}\label{Young_functions}\label{App:Young}

This appendix contains a brief overview as well as some new results on moderate Young functions.
In particular, in subsection \ref{sec:StepB}, we are interested in using the Burkholder--Davis--Gundy inequality for moderate Young functions and Doob's maximal inequality for functions whose Young conjugate is moderate.

\medskip
We will use and adapt the terminology of \cite[Chapter VI, Section 3, Paragraph 97]{dellacherie1982probabilities}. 
A function $\Phi:\mathbb R_+\longrightarrow\mathbb R_+$ is called \textit{Young} if it is increasing, convex and satisfies 
\begin{align*}
\Phi(0)=0\hspace{0.3em} \text{ and }\hspace{0.3em}\lim_{x\to \infty} \frac{\Phi(x)}x=\infty.
\end{align*}
We can write every Young function $\Phi$ in the form $
\Phi(x)=\int_0^x \phi(t)\ud t,$
where $\phi:\mathbb R_+\longrightarrow \mathbb R_+$ is an increasing and c\`adl\`ag function.
Due to the growth condition, it is immediate to check that $\lim_{x\to\infty}\phi(x)=\infty,$ which implies that the c\`adl\`ag inverse of $\phi$, which is defined by
\begin{align}\label{def:psi}
\psi(s):=\inf\{t\ge 0, \phi(t)>s \}, \ s\in\mathbb R_+,
\end{align}
is real--valued and unbounded as well.

\begin{definition}\label{def:Young-constants}
Every Young function $\Phi$ is associated to the constants
\begin{align*}
\underline c_\Phi:=\inf_{x>0} \frac{x\phi(x)}{\Phi(x)}, \text{ and }
\overline c_\Phi:=\sup_{x>0} \frac{x\phi(x)}{\Phi(x)}.
\end{align*}
If $\overline c_\Phi<\infty$, then $\Phi$ is called moderate.
\end{definition}
Observe that by the immediate inequality $\Phi(x)\le x\phi(x)$, we have that $\overline c_\Phi\ge \underline c_\Phi\ge 1.$
A characterization of moderate Young functions is given in \cite[Theorem 3.1.1]{long1993martingale}.
In other words, a Young function $\Phi$ is moderate if and only if $\Phi(\lambda x)\le \lambda^{\overline c_\Phi} \Phi(x)$ for every $x\in\mathbb R_+$ and for every $\lambda \ge 1$.
However, for a Young function to be moderate, it turns out that we actually only need to prove the property for some $\lambda>1$, \emph{e.g.} for $\lambda=2$, see \cite[Definition 10.32, Lemma 10.33.2)]{he1992semimartingale}. The Young conjugate of $\Phi$ is the Young function $\Psi:\mathbb R_+\longrightarrow\mathbb R_+$ defined as
\begin{align*}
\Psi(x):=\int_0^x\psi(s)\ud s,
\end{align*}
where $\psi$ is the c\`adl\`ag inverse of $\phi$ defined in \eqref{def:psi}.
By \cite[Theorem 3.1.1]{long1993martingale}, we have that $\underline c_\Psi$ is the conjugate index of $\overline c_\Phi$ and $\overline c_\Psi$ is the conjugate index of $\underline c_\Phi$, \emph{i.e.}
\begin{align*}
\underline c_\Psi = \begin{cases}\frac{\overline c_\Phi}{\overline c_\Phi - 1}, \text{ if }\overline c_\Phi>1,\\ \infty, \text{ if }\overline c_\Phi=1,
\end{cases}
\hspace{-1em} \text{ and } 
\overline c_\Psi = \begin{cases}\frac{\underline c_\Phi}{\underline c_\Phi - 1}, \text{ if }\underline c_\Phi>1,\\ \infty,  \text{ if }\underline c_\Phi=1. \end{cases}
\end{align*}
Therefore, the Young conjugate of $\Phi$ is moderate if $\underline c_\Phi>1.$
In the following, to every sequence $\mathcal A:=(\alpha^k)_{k\in\mathbb N}$ such that $\alpha^0=0$ and $\alpha^k\le \alpha^{k+1}$, for every $k\in\mathbb N$, and $\lim_{k\to\infty}\alpha_k=\infty$, we will associate the Young function $\Phi^{\mathcal A}$ with $
\Phi^{\mathcal A}(x):=\int_0^x \sum_{k=1}^\infty\alpha^k\mathds{1}_{[k,k+1)}(t)\ud t.$ For convenience, we define $\phi^{\mathcal A}:\mathbb R_+\longrightarrow\mathbb R_+$ by $\phi^{\mathcal A}(t):=\sum_{k=1}^\infty\alpha^k\mathds{1}_{[k,k+1)}(t).$
If, moreover, $\alpha^{2k}\le 2\alpha^k$ for every $k\in\mathbb N,$ then the Young function $\Phi^{\mathcal A}$ is moderate, as an immediate consequence of the comments after \cref{def:Young-constants}.

\vspace{0.5em}
\begin{proposition}\label{prop:Young_functions}
Let $\mathcal A=(\alpha_k)_{k\in\mathbb N}$ be an increasing and unbounded sequence of positive integers for which holds $\alpha^k\le \alpha^{k+1}$ and $\alpha^{2k}\le 2\alpha^k$, for every $k\in\mathbb N$.
Let now $\Phi_{\mathcal A}$ be the moderate Young function associated to the sequence $\mathcal A$.
Define ${\rm quad}:\mathbb R_+\longrightarrow\mathbb R_+$ by ${\rm quad}(x)=\frac{1}{2} x^2$, and let $\Psi$ be the Young conjugate of $\Phi^{\mathcal A,{\rm quad}}:=\Phi^{\mathcal A}\circ {\rm quad},$ with associated right derivative $\psi$.
Then
\begin{enumerate}
\item\label{prop:Young_functions-1} 
	$\Phi_{\mathcal A,\textup{quad}},\Psi$ are moderate Young functions.
\item\label{prop:Young_functions-2} 
	$\psi$ is continuous and can be written as a concatenation of linear and constant functions defined on intervals.
		Besides, the slopes of the linear parts constitute a non--increasing sequence converging to $0$.
\item\label{prop:Young_functions-3}
	We have $\Psi(x)\le\textup{quad}(x)$, where the equality holds on a compact neighborhood of $0$, and 
\begin{align*}
\lim_{x\uparrow\infty} \big\{\textup{quad}(x)-\Psi(x) \big\}=\infty.
\end{align*}
\item\label{prop:Young_functions-4}
	There exists a Young function $\Upsilon$  such that $\Upsilon\circ \Psi = \textup{quad}$.
\end{enumerate}

\end{proposition}
\begin{proof}
\iftoggle{full}{
\begin{enumerate}
\item[\ref{prop:Young_functions-1}]
We will prove initially that $\Phi_{\mathcal A,\textup{quad}}$ is a Young function.
In view of the comments before and after \cref{def:Young-constants}, it is sufficient to prove that it can be written as a Lebesgue integral whose integrand is a \cadlag, increasing and unbounded function.

\vspace{0.5em}
\noindent For every $x\in\mathbb R_+$, we have by definition 
\begin{align}\label{right-der:comp}
\Phi_{\mathcal A,\textup{quad}}(x)
	= \Phi_{\mathcal A}\Big(\frac{1}{2}x^2\Big) 
	= \int_{[0,\frac{1}{2}x^2]} \phi_{\!_{\mathcal A}}(z)\,\textup{d} z 
	= \int_{[0,x]} t \phi_{\!_{\mathcal A}}\Big(\frac12 t^2\Big) \,\textup{d} t.
\end{align}
We define $\phi_{\!_{\mathcal A,\textup{quad}}}:\mathbb R_+\longrightarrow\mathbb R_+$ by
\begin{align}\label{def:phi-A-quad}
\phi_{\!_{\mathcal A,\textup{quad}}}(t):=t \phi_{\!_{\mathcal A}}\Big(\frac12 t^2\Big) = t \mathds{1}_{[0, \sqrt{2})}(t) + t \sum_{k=1}^\infty \alpha^k \mathds{1}_{[\sqrt{2k}, \sqrt{2k+2})}(t),
\end{align}
\emph{i.e.} $\phi_{\!_{\mathcal A,\textup{quad}}}$ is \cadlag and piecewise-linear.
Observe, moreover, that 

\vspace{0.5em}
\begin{enumerate}[label={\tiny\raisebox{1pt}{$\blacksquare$}},itemindent=0.0cm,leftmargin=*]
	\item $\Delta \phi_{\!_{\mathcal A,\textup{quad}}}(\sqrt{2k+2})=(\alpha^{k+1}-\alpha^k)\sqrt{2k+2}\ge 0$, for every $k\in\mathbb N$.
	\item $\phi_{\!_{\mathcal A,\textup{quad}}}$ has increasing slopes; the value of the slope of the linear part defined on the interval $[\sqrt{2k}, \sqrt{2k+2})$ is determined by the value of the respective element $\alpha^k\ge 1$, for every $k\in\mathbb N$.
	\item $\lim_{s\to\infty}\phi_{\!_{\mathcal A,\textup{quad}}}(s)=\infty.$
\end{enumerate}

\vspace{0.5em}
\noindent 
Therefore, $\Phi_{\mathcal A,\textup{quad}}$ is a Young function and its conjugate $\Psi$ is also a Young function.

\vspace{0.5em}
\noindent We will prove now that both $\Phi_{\mathcal A,\textup{quad}}$ and $\Psi$ are moderate. 
We have directly that $\underline c_\textup{quad} =\overline c_\textup{quad} =2$. 
Moreover, by the property $\alpha^{2k}\le 2\alpha^k$ we have that $\Phi_{\mathcal A}$ is moderate, hence $\overline c_{\Phi_{\!_\mathcal A}}<\infty.$ 
Now we obtain
\begin{align*}
\underline c_{\Phi_{\mathcal A,\textup{quad}}}
	& = \inf_{x>0} \frac{x \phi_{\!_{\mathcal A,\textup{quad}}}(x)}{\Phi_{\mathcal A,\textup{quad}}(x)} 
	  = \inf_{x>0} \frac{x^2 \phi_{\!_{\mathcal A}}(\frac12 x^2)}{\Phi_{\mathcal A}(\frac12 x^2)}
	  = 2 \inf_{x>0} \frac{\frac12x^2 \phi_{\!_{\mathcal A}}(\frac12 x^2)}{\Phi_{\mathcal A}(\frac12 x^2)}
	  = 2 \inf_{u>0} \frac{u \phi_{\!_{\mathcal A}}(u)}{\Phi_{\mathcal A}(u)} 
	  = 2 \underline c_{\Phi_{\mathcal A}} \ge 2,
\end{align*}
because for every Young function $\Upsilon$ holds $\underline c_{_\Upsilon}\ge 1$.
In addition, for $\overline c_{\Phi_{\mathcal A,\textup{quad}}}$ we have 
\begin{align*}
\overline c_{\Phi_{\mathcal A,\textup{quad}}} 
	& = \sup_{x>0} \frac{x \phi_{\!_{\mathcal A,\textup{quad}}}(x)}{\Phi_{\mathcal A,\textup{quad}}(x)} 
	  = 2 \sup_{u>0} \frac{u \phi_{\!_{\mathcal A}}(u)}{\Phi_{\mathcal A}(u)} = 2 \overline c_{\Phi_{\mathcal A}} <\infty.
\end{align*}
Hence, $\Phi_{\mathcal A,\textup{quad}}$ is a moderate Young function.
Besides, since $\underline c_{\Phi_{\mathcal A,\textup{quad}}}>1$, we have from \cite[Theorem 3.1.1 (f)]{long1993martingale} that $\overline c_\Psi<\infty$.
Therefore, $\Psi$ is also moderate.

\vspace{0.7em}
\item[\ref{prop:Young_functions-2}]
For the rest of the proof, \emph{i.e.} for parts (ii)--(iv), we will simplify the notation and simply write $\phi$ for the function $\phi_{_{\mathcal A,\textup{quad}}}$.

\vspace{0.5em}
\noindent Firstly, let us observe that $\psi$ is real--valued, resp. unbounded, when $\phi$ is unbounded, resp. real--valued.\footnote{The reader who is not familiar with generalized inverses may find \citet{embrechts2013note} helpful, especially the comments after \citep[Remark 2.2]{embrechts2013note}. }
In order to determine the value $\psi(s)$ for $s\in(0,\infty)$, let us define two sequences of subsets of $\mathbb R_+$, $(C^k)_{k\in\mathbb N\cup\{0\}}$, $(J^k)_{k\in\mathbb N\cup\{0\}}$, by
\begin{align}\label{def:Ck-Jk}
C^k:=\big[\phi(\sqrt{2k}), \lim_{t\uparrow \sqrt{2k+2}}\phi(t) \big) 
\text{ and }
J^{k+1}:=\big[\lim_{t\uparrow \sqrt{2k+2}}\phi(t), \phi(\sqrt{2k+2}) \big), \text{ for }k\in\mathbb N\cup\{0\}.
\end{align} 
Observe that 

\vspace{0.5em}
\begin{enumerate}[label={\tiny\raisebox{1pt}{$\blacksquare$}},leftmargin=*]
	\item $C^k=\phi\big([\sqrt{2k},\sqrt{2k+2})\big)\neq \emptyset$, $\forall k\in\mathbb N\cup\{0\}$, since $\phi$ is continuous and increasing on $[\sqrt{2k},\sqrt{2k+2}).$
	\item $J^k=\emptyset$ if and only if $\phi$ is continuous at $\sqrt{2k+2}$, which is further equivalent to $\alpha^k=\alpha^{k+1}.$ 
\end{enumerate}

\vspace{0.5em}
\noindent For convenience, let us define two sequences $(s^k)_{k\in\mathbb N\cup\{0\}}$, $(s^{k+1}_{-})_{k\in\mathbb N\cup\{0\}}$ of positive numbers as follows
\begin{equation}\label{def:sk-sk--}
\begin{split}
s^0&:=0,\hspace{.3em}s^k:=\phi(\sqrt{2k})=\alpha^k \sqrt{2k}, \text{ for } k\in\mathbb N \text{ and }\\
s^{k+1}_{-}&:=\lim_{x\uparrow\sqrt{2k+2}}\phi(x)=\alpha^k\sqrt{2k+2}, \text{ for }k\in\mathbb N\cup\{0\}.
\end{split}
\end{equation}
The introduction of the last notation allows us to rewrite $C^k=[s^k, s^{k+1}_{-})$ and $J^{k+1}=[s^{k+1}_{-}, s^{k+1})$, for $k\in\mathbb N\cup\{0\}$.
Now we are ready to determine the values of $\psi$ on $(0,\infty)$. 
The reader should keep in mind that the function $\phi$ is increasing and right--continuous. 

\vspace{0.5em}
\begin{enumerate}[label=$\bullet$,leftmargin=*,itemindent=0.0cm]
 	\item  Let $s\in C^0=\big[\phi(0),\phi(\sqrt2-)\big)=\phi\big([0,\sqrt2)\big)$, then
 			\begin{align*}
 				\psi(s)	& \overset{\phantom{\eqref{def:phi-A-quad}}}{=}	
 						\inf\{t\in \mathbb R_+, \phi(t)>s \}
 					 \overset{s\in\phi\left([0,\sqrt2)\right)}{=}	\inf\big\{t\in [0,\sqrt2), \phi(t)>s\big\}\\
 					& \overset{\eqref{def:phi-A-quad}}{=}
 						\inf\{t\in \mathbb R_+, \textup{Id}(t)\mathds{1}_{[0,\sqrt2)}(t)>s \} = s,
 			\end{align*} 
		where the second equality is valid because $\phi$ is continuous on $[0,\sqrt2)$ with $\phi\big([0,\sqrt2)\big)=[0,s^{1}_{-}).$
		To sum up, we have proven that $\psi\mathds{1}_{[s^0,s^{1}_{-})} = \textup{Id}\mathds{1}_{[s^0,s^{1}_{-})}.$
\vspace{0.5em}
	\item Let $s\in J^1=\big[\phi(\sqrt2-), \phi(\sqrt2)\big)=[s^{1}_{-}, s^1)$. 
			If $J^1=\emptyset$, which amounts to $\alpha^1=1,$ there is nothing to prove.
			On the other hand, if $J^1\neq\emptyset$, then $\phi(\sqrt2\,-) \le s< \phi(\sqrt2)$ and consequently
			\begin{align*}
				\psi(s)	 =	\inf\{t\in\mathbb R_+, \phi(t)> s\} = \sqrt{2} \text{ for every }s\in J^1.
			\end{align*}
		To sum up, we have proven that $\psi\mathds{1}_{[s^{1}_{-},s^1)} = \sqrt2\mathds{1}_{[s^{1}_{-},s^1)}.$ 
 \end{enumerate} 
\vspace{0.5em}
\noindent For the general case let us fix a $k\in\mathbb N$.
We will distinguish between the cases $s\in C^k$ and $s\in J^k$.
For the latter we can argue exactly as in the case $s\in J^1$, but for the sake of completeness we will provide the proof.

\vspace{0.5em}
\begin{enumerate}[label=$\bullet$,leftmargin=*,itemindent=0.0cm]
	\item Let $s\in C^k=\big[\phi(\sqrt{2k}),\phi(\sqrt{2k+2}\,-)\big)=\phi\big([\sqrt{2k},\sqrt{2k+2})\big)$.
			Since $C^k$ is the image of $[\sqrt{2k},\sqrt{2k+2})$ through $\mathbb R_+\ni t\overset{\chi}{\longmapsto} \alpha^{k}t \in\mathbb R_+$, 
			then $\psi$ has to coincide with $\mathbb R_+\ni s \xmapsto{\hspace{.3em}\chi^{-1,r}\hspace{.3em}} \frac{1}{\alpha^{k}}s\in\mathbb R_+$ on $C^k=[s^k, s^{k+1}_{-}).$
			To sum up, we have proved that $\psi(s)\mathds{1}_{[s^k, s^{k+1}_{-})} = \frac{1}{\alpha^{k}}\textup{Id}\mathds{1}_{[s^k, s^{k+1}_{-})}$.
			
			\vspace{0.5em}
	\item Let $s\in J^{k+1}=\big[\phi(\sqrt{2k+2}\,-), \phi(\sqrt{2k+2})\big)=[s^{k+1}_{-}, s^{k+1})$. 
			If $J^{k+1}=\emptyset$, which amounts to $\alpha^{k+1}=\alpha^k,$ there is nothing to prove.
			On the other hand, if $J^{k+1}\neq\emptyset$, then $\phi(\sqrt{2k+2}\,-) \le s< \phi(\sqrt{2k+2})$ and consequently
		\begin{align*}
			\psi(s)	 =	\inf\{t\in\mathbb R_+, \phi(t)> s\} = \sqrt{2k+2}, \text{ for every }s\in J^{k+1}.
		\end{align*}
		To sum up, we have proven that $\psi\mathds{1}_{[s^{k+1}_{-},s^{k+1})} = \sqrt{2k+2}\mathds{1}_{[s^{k+1}_{-},s^{k+1})}.$ 
\end{enumerate}

\vspace{0.5em}
\noindent Overall, we have that the right derivative of $\Psi$ can be written as a concatenation of linear and constant functions defined on intervals, \textit{i.e.}
\begin{align}\label{def:psiAquad}
\psi(s)& =\textup{Id}\mathds{1}_{[s^0,s^1_{-})} + \sum_{k=1}^\infty \frac{1}{\alpha^{k}}\textup{Id}\mathds{1}_{[s^k, s^{k+1}_{-})}+ \sum_{k=0}^\infty \sqrt{2k+2} \mathds{1}_{[s^{k+1}_{-}, s^{k+1})},
\end{align}
where
\begin{align*}
s^0:=0, \; s^k:=\alpha^k\sqrt{2k}, \text{ for }k\in\mathbb N, \text{ and } s^{k+1}_{-}:=\alpha^k\sqrt{2k+2}, \text{ for }k\in\mathbb N\cup\{0\}. 
\end{align*}
Recall now that $\alpha^k\ge 1$, for every $k\in\mathbb N$, therefore we have that the slopes of $\psi$ are smaller than $1$.
Moreover, since $\alpha^k\le \alpha^{k+1}$ and $\lim_{k\to\infty}\alpha^k=\infty,$ we have that $\frac{1}{\alpha^{k+1}}\le \frac{1}{\alpha^k}$ and $\lim_{k\to\infty}\frac{1}{\alpha^k}=0.$ 
Finally, as it can be easily checked, $\psi$ is continuous.
This causes no surprise, since $\phi$ is strictly increasing, see \citet[Proposition 2.3.(7)]{embrechts2013note}\footnote{For the following we keep our notation. In \citep{embrechts2013note} the presented results is for the left--continuous generalized inverse of a function. However, $\phi^{-1,l}$ is the \cadlag version of $\psi$, therefore we can directly conclude that $\psi$ has to be also continuous.}.

\vspace{0.7em}
\item[\ref{prop:Young_functions-3}]
Let us consider the function $\zeta:\mathbb R_+\to\mathbb R$ defined by $\zeta:=\textup{Id}-\psi,$ \emph{i.e.} $\zeta$ is also continuous.
Moreover, $\zeta$ is differentiable on a superset of $\mathbb R_+\setminus\{(s^k)_{k\in\mathbb N\cup\{0\}}\cup(s^{k+1}_{-})_{k\in\mathbb N\cup\{0\}}\}$, which is clearly an open and dense subset of $\mathbb R_+$, since there is no accumulation point in the sequence $(s^k)_{k\in\mathbb N\cup\{0\}}\cup(s^{k+1}_{-})_{k\in\mathbb N\cup\{0\}}$.
Obviously
\begin{align*}
\zeta'(x)=
\begin{cases}
1-\frac{1}{\alpha^{k}}, &\text{for } x\in(s^k,s^{k+1}_{-}) \text{ for }k\in\mathbb N\cup\{0\},\\
1,						&\text{for } x\in(s^{k+1}_{-}, s^{k+1}) \text{ for }k\in\mathbb N\cup\{0\}, \text{ whenever }(s^{k+1}_{-}, s^{k+1}) \neq \emptyset.
\end{cases}
\end{align*}
Define ${M}:=\min\{k\in\mathbb N, \alpha^k>1\},$ which is a well-defined positive integer, since $\alpha^k\xrightarrow{\hspace{1em}} \infty.$
Recall now that $\alpha^k\ge 1$ for $k\in\mathbb N$ and we can conclude that $\zeta'(x)>0$ almost everywhere on $[s^{{M}},\infty)$.

\vspace{0.5em}
\noindent We prove now that \textup{quad} and $\psi$ coincide only on a compact neighborhood of $0.$
By the definition of ${M}$ we have that $\alpha^k=1$ for $k\in\{1,\dots,{M}-1\}$ and $\alpha^{{M}}>1$, 
therefore $\textup{Id}\mathds{1}_{[0,s^{{M}}]}=\psi\mathds{1}_{[0,s^{{M}}]}$ 
and $x>\psi(x)$ for every $x\in(s^{{M}},\infty).$

\vspace{0.5em}
\noindent Finally, it is left to prove that $\lim_{x\uparrow\infty}(\textup{Id}(x)-\psi(x))=\infty.$
Recall that $(s^{k+1}_{-}, s^{k+1})\neq \emptyset$ whenever $k$ is such that $\alpha^k<\alpha^{k+1}$, \emph{i.e.} there are infinitely many non-trivial intervals $(s^{k+1}_{-}, s^{k+1})$.
But these intervals correspond to the intervals where $\psi$ is constant.
Since \textup{Id} is increasing, we can conclude that $\zeta$ is unbounded and that the desired result holds.

\vspace{0.7em}
\item[\ref{prop:Young_functions-4}]
For the following recall \eqref{def:sk-sk--}, \eqref{def:psiAquad} and the definition of $M.$
Let us start with the introduction of the auxiliary function $\eta:\cup_{k=M}^\infty\big[\Psi(s^{k+1}_{-}),\Psi(s^{k+1})\big)\xrightarrow{\hspace{1em}}(0,1]$ defined by $\eta(z):=\frac{\Psi(s^{k+1})-z}{\Psi(s^{k+1})-\Psi(s^{k+1}_{-})}$\footnote{In fact, $\eta(z)$ is the unique number in $[0,1)$ for which $z$ can be written as convex combination of $\Psi(s^{k+1}_{-})$ and $\Psi(s^{k+1})$.}.
Recall now that $\Psi$ is continuous and increasing, which allows us to define the function $\upsilon:\mathbb R_+\longrightarrow\mathbb R_+$ by
\begin{align}\label{def:upsilon-proof}
\begin{split}
\upsilon(z):  = 	\mathds{1}_{[0, \Psi(s^M))}(z)
			&  +		\sum_{k=M}^\infty 	\alpha^{k}\mathds{1}_{[\Psi(s^k), \Psi(s^{k+1}_{-}))}(z) \\
			&  +  	\sum_{k=M}^\infty	\left(\eta(z) \alpha^{k} +(1-\eta(z)) \alpha^{k+1}\right) \mathds{1}_{[\Psi(s^{k+1}_{-}), \Psi(s^{k+1}))}(z).
\end{split}
\end{align}
We can directly check that $\upsilon$ is indeed well--defined, non--negative, non--decreasing and unbounded.
Therefore, the function $\upsilon$ is the right derivative of a Young function, say $\Upsilon$.

\vspace{0.5em}
\noindent
We intend to prove that $\Upsilon\circ\Psi = \textup{quad}$, which is equivalent to proving that the right--derivative of $\Upsilon\circ\Psi$ equals \textup{Id}.
The following simple calculations allow us to evaluate the right derivative of $\Upsilon\circ\Psi$ in terms of $\upsilon,$ $\psi$, and $\Psi.$
For any $x\in\mathbb R_+\to\mathbb R_+$, we have
\begin{align*}
\Upsilon\circ\Psi(x)
	=\int_{[0,\Psi(x)]} \upsilon(t)\,\textup{d}t 
	= \int_{[0,x]} \psi(z)\upsilon\big(\Psi(z)\big)\,\textup{d}z.
\end{align*}
Now we can compare the right derivative of $\Upsilon\circ \Psi$, which is the function $\psi (\upsilon\circ\Psi):\mathbb R_+\longrightarrow\mathbb R_+$, with the identity function \textup{Id}.	
To this end, we will consider the behaviour of $\psi (\upsilon\circ\Psi)$ on the intervals $[0,s^M)$, $[s^k, s^{k+1}_{-})$, for $k\ge M$ and $[s^{k+1}_{-},s^{k+1})$ for $k\ge M,$ which form a partition of $\mathbb R_+.$
Before we proceed, let us evaluate the function $\upsilon$ at $\Psi(s)$ for $s\in\mathbb R_+$ 
\begin{align}\label{def:upsilon-at-Psi}
\nonumber\upsilon\big(\Psi(s)\big) 
			 	= &\	\mathds{1}_{[0, \Psi(s^M))}(\Psi(s))
			  	+		\sum_{k=M}^\infty 	\alpha^{k}\mathds{1}_{[\Psi(s^k), \Psi(s^{k+1}_{-}))}(\Psi(s))\\
\nonumber	&	
			  	+  	\sum_{k=M}^\infty	\left\{\eta\big({\Psi(s)}\big) \alpha^{k} +\big[1-\eta\big({\Psi(s)}\big)\big] \alpha^{k+1}\right\} \mathds{1}_{[\Psi(s^{k+1}_{-}), \Psi(s^{k+1}))}(\Psi(s))\\
\nonumber	 	=&\		\mathds{1}_{[0,s^M)}(s) 
				+  	\sum_{k=M}^\infty 	\alpha^{k}\mathds{1}_{[s^k, s^{k+1}_{-})}(s)\\
			&	
				+ 	\sum_{k=M}^\infty	\left\{\eta\big({\Psi(s)}\big) \alpha^{k} +\big[1-\eta\big({\Psi(s)}\big)\big] \alpha^{k+1}\right\}  \mathds{1}_{[s^{k+1}_{-}, s^{k+1})}(s) 
\end{align}
because $\Psi$ is continuous and increasing.

\vspace{0.5em}
\begin{enumerate}[label=$\bullet$,leftmargin=*,itemindent=0.0cm]
\item Let $s\in[0,s^M)$. 
		At the end of the proof of \ref{prop:Young_functions-3} we obtained that $\textup{Id}\mathds{1}_{[0,s^{M})}=\psi\mathds{1}_{[0,s^{M})}$.
		Therefore, we can conclude that $
			\psi(s)(\upsilon\circ \Psi)(s) 
				= s\upsilon\big(\Psi(s)\big)
				= s.$ 
				
				\vspace{0.5em}
\item Let $s\in [s^k, s^{k+1}_{-})$, for some $k\ge M.$ 
		Then, 
		\begin{align*}
			\psi(s)\upsilon(\Psi(s)) \overset{\eqref{def:psiAquad}}{\underset{\eqref{def:upsilon-at-Psi}}{=}} \frac{1}{\alpha^k} \alpha^k s= s. 
		\end{align*}

\item Let $s\in[s^{k+1}_{-}, s^{k+1})$ for some $k\ge M.$
		Then, for the chosen $s$ there exists (unique) $\mu_s\in[0,1)$ such that 
		\begin{align}\label{prf-prpYF-defss}
		s=\mu_s s^{k+1}_{-} + (1-\mu_s) s^{k+1}.
		\end{align}
		However, $\Psi$ is linear on $[s^{k+1}_{-}, s^{k+1})$, \emph{i.e.}
		\begin{align}\label{eqlt:Psi-s}
			\Psi(s) = \mu_s \Psi(s^{k+1}_{-}) + (1-\mu_s) \Psi(s^{k+1})
		\end{align}
		and $\Psi(s)\in\left[\Psi(s^{k+1}_{-}), \Psi(s^{k+1})\right)$.
		Therefore, by definition of $\eta$ 
		\begin{align*}
			\eta(\Psi(s)) 
				&	= \frac{\Psi(s^{k+1}) - \Psi(s) }{\Psi(s^{k+1})-\Psi(s^{k+1}_{-})} 
					\overset{\eqref{eqlt:Psi-s}}{=} \frac{\Psi(s^{k+1}) - \mu_s \Psi(s^{k+1}_{-}) - (1-\mu_s) \Psi(s^{k+1}) }{\Psi(s^{k+1})-\Psi(s^{k+1}_{-})}\\
				& 	=  \frac{\cancel{\Psi(s^{k+1})} - \mu_s \Psi(s^{k+1}_{-}) - \cancel{\Psi(s^{k+1})} + \mu_s \Psi(s^{k+1})}{\Psi(s^{k+1})-\Psi(s^{k+1}_{-})} = \mu_s 
				\numberthis\label{evaluate:eta}
		\end{align*}
and finally, we can conclude in view of
		\begin{align*}
		\psi(s) (\upsilon\circ\Psi)(s) 
		& = \sqrt{2k +2}\, \upsilon(\Psi(s)) 
		    \overset{\eqref{def:upsilon-at-Psi}}{=}\sqrt{2k+2}\left\{\eta\big({\Psi(s)}\big) \alpha^{k} +\big[1-\eta\big({\Psi(s)}\big)\big] \alpha^{k+1}\right\}  \\
		&   \overset{\eqref{evaluate:eta}}{=} \sqrt{2k +2} (\mu_s \alpha^k + (1-\mu_s)\alpha^{k+1})\\
		&	\overset{\eqref{def:sk-sk--}}{=} \mu_s s^{k+1}_{-} + (1-\mu_s)s^{k+1}
			\overset{\eqref{prf-prpYF-defss}}{=}s.
		\end{align*}
\end{enumerate}
\qedhere
\end{enumerate}
}
{
	See the proof of Proposition A.12 in the extended version \cite{papapantoleon2018stability}.
}
\end{proof}

\iftoggle{full}{

\subsection{Proof of \cref{corr:PII}}\label{Proof_Corr_main_thm}

\begin{proof}
Let us initially define the random variables
\begin{equation*}
D^k:=(X^k_\infty, \xi^k, X^k)\text{ for }k\in\overline{\mathbb N}.
\end{equation*}
Observe that the state space of $D^k$ is $\mathbb R^{\ell + 1}\times\mathbb D(\mathbb R_+;\mathbb R^{\ell})$, which is clearly Polish as a finite Cartesian product of Polish spaces.
Moreover, observe that $
	D^k\xrightarrow{\hspace{0.2em}\mathcal L\hspace{0.2em}} D^\infty,$
due to Assumptions \ref{Weak-MSI} and \ref{Weak-Mfinalrv}.
We are going to use the Skorokhod Representation Theorem, see \citet[Theorem~6.7]{billingsley1999convergence}, in order to obtain a probability space $(\overline{\Omega}, \overline{\mathcal F}, \overline{\mathbb P})$ and a sequence $(\overline{D}^k)_{k\in\overline{\mathbb N}}$ of random variables defined on $(\overline{\Omega}, \overline{\mathcal F}, \overline{\mathbb P})$ such that
\begin{enumerate}
	\item\label{Dk:Equally-Distributed} $\mathcal L(\overline{D}^k) = \mathcal L(D^k)$, for every $k\in\overline{\mathbb N}$,
	\item \label{Skorokhod-Conv} $\overline{D}^k\xrightarrow{\hspace{0.2em}\delta_\Pi\hspace{0.2em}}\overline{D}^\infty$, $\overline{\mathbb P}-$almost surely, 
	where $\delta_\Pi:\big(\mathbb R^{\ell+1}\times\mathbb D(\mathbb R_+;\mathbb R^\ell)\big)\times \big(\mathbb R^{\ell+1}\times\mathbb D(\mathbb R_+;\mathbb R^\ell)\big)\longrightarrow \mathbb R_+ $ with
	\begin{equation*}
		\delta_\Pi\big((x,\alpha), (y,\beta)\big):= |x-y| + \delta_{\textup{J}_1(\mathbb R^\ell)}(\alpha,\beta) \quad \text{ for }x,y\in\mathbb R^{\ell+1} \text{ and }\alpha,\beta\in\mathbb D(\mathbb R_+; \mathbb R^\ell).
	\end{equation*}
\end{enumerate}
The expectation under the measure $\overline{\mathbb P}$ will be denoted by $\overline{\mathbb E}[\cdot]$ and the conditional expectation of a random variable $Z$ with respect to a $\sigma-$algebra $\mathcal H$ under the measure $\overline{\mathbb P}$ will be denoted by $\overline{\mathbb E}[Z|\mathcal H].$ 
In view of the above, $\overline{D}^k$ can be written as $(\overline{X}^k_\infty, \overline{\xi}^k, \overline{X}^k)$ for some $\mathbb R^\ell-$, resp. $\mathbb R-$valued, random variables $\overline{X}^k_\infty, \overline{\xi}^k$ and some $\mathbb R^\ell-$valued process $\overline{X}^k$, for every $k\in\overline{\mathbb N}.$
The next step is to construct the stochastic basis $(\overline{\Omega}, \overline{\mathcal F}, \mathbb F^{\overline{X}^k},\overline{\mathbb P})$ for every $k\in\overline{\mathbb N}$ .
In order to make clear the correspondence with the conditions of \cref{RobMartRep}, we define $\mathbb G^k:=\mathbb F^{\overline{X}^k}$ for every $k\in\overline{\mathbb N}.$ 
Then, we can translate Conditions \ref{Weak-MFilqlc}, \ref{Weak-MSI}, \ref{Weak-MXWPRP} and \ref{Weak-Mfinalrv} into Conditions \ref{MFilqlc}, \ref{MSI}, \ref{MXWPRP} and \ref{Mfinalrv} under the probability space $(\overline{\Omega}, \overline{\mathcal F}, \overline{\mathbb P})$. 
Moreover, in view of Condition \ref{Weak-MFilweak}, we can conclude that $\overline{X}^k$ is a process with independent increments. 
This property in conjunction with the convergence obtained by \ref{Skorokhod-Conv} implies in particular that $\mathbb F^{\overline{X}^k}\xrightarrow{\hspace{0.2em}\textup{w}\hspace{0.2em}}\mathbb F^{\overline{X}^\infty}.$
Our last claim is verified by \citet[Proposition 2]{coquet2001weak}.
Therefore, Condition \ref{MFilweak} is also satisfied when we work on the probability space $(\overline{\Omega}, \overline{\mathcal F},\overline{\mathbb P})$.

\vspace{0.4em}
Finally, we need to prove also that $\mathcal L\big(\mathbb E^k[\xi^k|\mathcal F^{X^k}_\cdot]\big) = \mathcal L\big(\overline{\mathbb E}[\,\overline{\xi}^k|\mathcal F^{\overline{X}^k}_\cdot]\big)$ for every $k\in\overline{\mathbb N}$, in order to be able to transfer the results from $(\overline{\Omega}, \overline{\mathcal F},\overline{\mathbb P})$ to the original space.
We underline that the last statement needs some special care, since the Skorokhod representation theorem does not deal with the associated filtrations.
However, we will see in the next lines that, if we work with the natural filtrations, then we can assume that the laws of the corresponding optional projections (or simply of the conditional expectations) are unaffected.
For the following we will use \citet[Theorem VI.1.14 c)]{jacod2003limit}. 
In other words, we will use the fact that for every $t\in\mathbb R_+$ the $\sigma-$algebra $\mathscr D^0_t(\mathbb R^\ell)$ 
is the $\sigma-$algebra generated by all maps 
$\mathbb D(\mathbb R_+;\mathbb R^\ell)\ni\alpha\longmapsto \alpha(u)\in\mathbb R^\ell$ for $u\le s$, 
and coincides with the Borel $\sigma-$algebra associated to the Polish space $\big(\mathbb D([0,t];\mathbb R^\ell), \delta_{\textup{J}_1(\mathbb R^\ell)}\big)$.

\vspace{0.5em}
\noindent Define now $\mathscr D(\mathbb R^\ell):=\bigvee_{t\in\mathbb R_+}\mathscr D_t^0(\mathbb R^\ell)$.
In order to prove our claim that $\mathcal L\big(\mathbb E^k[\xi^k|\mathcal F^{X^k}_\cdot]\big) = \mathcal L\big(\overline{\mathbb E}[\,\overline{\xi}^k|\mathcal F^{\overline{X}^k}_\cdot]\big)$, it is sufficient by Kolomogorov's extension theorem to prove that for every finite subset $\mathcal T$ of $\mathbb R_+$, the finite dimensional distributions associated to the processes $\mathbb E^k[\xi^k|\mathcal F^{X^k}_\cdot]$ and $\overline{\mathbb E}[\,\overline{\xi}^k|\mathcal F^{\overline{X}^k}_\cdot]$ coincide.

\vspace{0.5em}
We will initially prove that for every $t\in\mathbb R_+$ and for every $k\in\overline{\mathbb N}$ holds
\begin{equation}\label{CE:fixed-t}
 	\mathbb E^k[\xi^k |\mathcal F^{X^k}_t] = g^k_t(X^k)
 	\ \ \text{ as well as }\ \ 
 	\overline{\mathbb E}[\overline{\xi}^k |\mathcal F^{\overline{X}^k}_t] = g^k_t(\overline{X}^k),
\end{equation} 
where $g^k_t:\big(\mathbb D(\mathbb R_+;\mathbb R^\ell), \mathscr D_t(\mathbb R^\ell)\big)\longrightarrow \big(\mathbb R,\mathcal B(\mathbb R)\big)$.

\vspace{0.5em}
Fix $(t,k)\in\mathbb R_+\times\overline{\mathbb N}$. Doob's lemma (\citet[Lemma 1.13]{kallenberg2002foundations}), the $\mathcal F^{X^k}_t-$measurability of $\mathbb E^k[\xi^k |\mathcal F^{X^k}_t]$, 
and the $\mathcal F^{\overline{X}^k}_t-$measurability of $\overline{\mathbb E}[\overline{\xi}^k |\mathcal F^{\overline{X}^k}_t]$,
imply that there exist $g^k_t:\big(\mathbb D(\mathbb R_+;\mathbb R^\ell), \mathscr D_t(\mathbb R^\ell)\big)\longrightarrow \big(\mathbb R,\mathcal B(\mathbb R)\big)$, and $\bar g^k_t:\big(\mathbb D(\mathbb R_+;\mathbb R^\ell), \mathscr D_t(\mathbb R^\ell)\big)\longrightarrow \big(\mathbb R,\mathcal B(\mathbb R)\big)$,
such that $\mathbb E^k[\xi^k |\mathcal F^{X^k}_t] = g^k_t(X^k)$, 
and $\overline{\mathbb E}[\overline{\xi}^k |\mathcal F^{\overline{X}^k}_t] = \bar g^k_t(\overline{X}^k)$.
Our aim, now, is to prove that 
\begin{equation}\label{g-barg-are-Equal}
	g_t^k=\bar g^k_t,\;  \mathbb P_{X^k}-a.s. \; \Big(\text{hence also }\mathbb P_{\overline{X}^k}-a.s.\Big),
\end{equation}
where $\mathbb P_{X^k}$, resp. $\mathbb P_{\overline{X}^k}$, denotes the push forward measure $\mathbb P_{X^k}:\big( \mathbb D(\mathbb R_+;\mathbb R^\ell), \mathscr D(\mathbb R^\ell)\big)\longrightarrow \big(\mathbb R_+, \mathcal B(\mathbb R_+)\big)$ defined by $\mathbb P_{X^k}(A):=\mathbb P^k\big((X^k)^{-1}(A)\big)$.
The interpretation of $\mathbb P_{\overline{X}^k}$ is analogous.
At this point recall \ref{Dk:Equally-Distributed} and the definition of the conditional expectation, in order to derive the following equality
\begin{equation*}
	\mathbb E^k\big[ g_t^k(X^k) \mathds{1}_{(X^k)^{-1}(A)}\big]
		=\mathbb E^k\big[ \xi^k \mathds{1}_{(X^k)^{-1}(A)}\big]
		=\overline{\mathbb E}\big[\, \overline{\xi}^k \mathds{1}_{(\overline{X}^k)^{-1}(A)}\big]
		=\overline{\mathbb E}\Big[\bar g_t^k(\overline{X}^k) \mathds{1}_{(\overline{X}^k)^{-1}(A)}\Big]
\end{equation*}
for every $A\in\mathcal B\big(\mathbb R^\ell)\big).$
Using classical approximation arguments, we can prove using the above equality that 
\begin{equation*}
\mathbb E^k\big[ g_t^k(X^k) f\big(X^k\big)\big] 
	= \overline{\mathbb E}\big[ \, \bar g_t^k (\overline{X}^k)f\big(\overline{X}^k\big)\big],
\end{equation*}
for every bounded $f:\big( \mathbb D(\mathbb R_+;\mathbb R^\ell), \mathscr D(\mathbb R^\ell)\big) \longrightarrow \big(\mathbb R,\mathcal B(\mathbb R^\ell)\big)$.
In particular, using the above equality for every bounded $h_t:\big( \mathbb D(\mathbb R_+;\mathbb R^\ell), \mathscr D_t(\mathbb R^\ell)\big) \longrightarrow \big(\mathbb R,\mathcal B(\mathbb R^\ell)\big)$, of the fact that both $g_t^k$, and $\bar g_t^k$ are $\mathscr D_t(\mathbb R^\ell)-$measurable, and that the measures $\mathbb P_{X^k}$ and $\mathbb P_{\overline X^k}$ are equal, we can conclude that Equality \eqref{g-barg-are-Equal} holds.  
Now, we can easily conclude that the associated finite dimensional distributions are equal in view of
\begin{align*}
&\mathbb P^k\big[\big\{\omega^k\in\Omega^k: \mathbb E^k[\xi^k|\mathcal F^{X^k}_{t_1}](\omega^k)\in A_1, \dots, \mathbb E^k[\xi^k|\mathcal F^{X^k}_{t_n}](\omega^k)\in A_n\big\}\big]\\
&\hspace{2em}=\mathbb P_{X^k}\bigg(\bigcap_{m=1}^n \big[(g^k_{t_1})^{-1}(A_1)] \bigg)
=\mathbb P_{\overline{X}^k}\bigg(\bigcap_{m=1}^n \big[(g^k_{t_1})^{-1}(A_1)] \bigg)\\
&\hspace{2em}=
\overline{\mathbb P}\Big[\Big\{\overline{\omega}\in\overline{\Omega}: 
\overline{\mathbb E}[\, \overline{\xi}^k|\mathcal F^{\overline{X}^k}_{t_1}](\overline{\omega})\in A_1, \dots, \overline{\mathbb E}[\,\overline{\xi}^k|\mathcal F^{\overline{X}^k}_{t_n}](\overline{\omega})\in A_n\Big\}\Big],
\end{align*}
for every $t_1, \dots,t_n\in\mathbb R_+$, $A_1,\dots,A_n\in\mathcal B(\mathbb R^\ell)$ and every $n\in\mathbb N.$ 
\end{proof}

} {}

\bibliographystyle{abbrvnat}
\small
\bibliography{bibliographyDylan}

\end{document}